\newif\ifJOLT 
\JOLTfalse

\ifJOLT
\documentclass{artjlt}
\else
\documentclass{amsart}
\fi
\usepackage{amsmath,amsfonts,amssymb}

\newif\ifLONG 
\LONGtrue 

\newif\ifEXTRAPROOFS 
\EXTRAPROOFStrue 

\newcommand{\n}{\noindent}

\newcommand{\bc}{\begin{center}}
\newcommand{\ec}{\end{center}}
\newcommand{\ds}{\displaystyle}
\newcommand{\ba}{$$\begin{aligned}}
\newcommand{\ea}{\end{aligned}$$}
\newcommand{\lp}{\left(}
\newcommand{\rp}{\right)}
\newcommand{\lb}{\left[}
\newcommand{\rb}{\right]}

\newcommand{\eml}{\end{multline}} 

\newcommand{\bml}{\begin{multline}}
\newcommand{\bmlst}{\begin{multline*}}

\newcommand{\bsp}{\begin{equation}\begin{split}}
\newcommand{\espe}{\end{equation}}
\newcommand{\esps}{\end{split}} 
\newcommand{\esp}{\end{split}\end{equation}} 
\newcommand{\bspst}{\begin{equation*}\begin{split}}
\newcommand{\espst}{\end{split}\end{equation*}} 

\newcommand\cP{\mathcal{P}}
\newcommand\cF{\mathcal{F}}

\ifJOLT
\else
\usepackage{amsthm}
\newtheorem{props}{Prop}[section] 
\newtheorem{thms}[props]{Theorem}

\newtheorem{cors}[props]{Corollary}

\newtheorem{lemmas}[props]{Lemma}
\fi

\title{Askey-Wilson Polynomials and Branching Laws 
}                                     
\author{Allen Back \\ Bent \O{}rsted \\ Siddhartha Sahi \\  Birgit Speh }                 
\ifJOLT
\lastname{Back, \O{}rsted, Sahi, Speh}  


\fi

\address{%
Allen Back \\               
Department Of Mathematics \\ Malott Hall, Cornell University \\ Ithaca, NY 14853 \\ United States \\            
ahb2@cornell.edu                
}
\address{%
Bent \O{}rsted\\               
Department of Mathematics \\
Aarhus University \\ 
Ny Munkegade 118 \\
Building 1530, 431 \\
8000 Aarhus C \\
Denmark \\
orsted@math.au.dk                
}
\address{%
Siddhartha Sahi\\               
Department of Mathematics \\ Rutgers University \\
110 Frelinghhuysen Rd \\
New Brunswick, NJ 08854-8019 \\  United States \\  
sahi@math.rutgers.edu                
}
\address{%
Birgit Speh \\               
Department Of Mathematics \\ Malott Hall, Cornell University \\ Ithaca, NY 14853 \\   United States \\            
bes12@cornell.edu                
}

\begin{document}

\ifJOLT
\maketitle
\fi 

\begin{abstract}
Connection coefficient formulas for special functions describe change of basis matrices under a parameter change, for bases formed by
the special functions.  Such formulas are related to branching questions in representation theory.
The Askey-Wilson polynomials are one of the most general $1$-variable special functions.  Our main results are connection
coefficient formulas for shifting one of the parameters of the nonsymmetric Askey-Wilson polynomials. We also show how one of these results can be used to re-prove an old result of Askey and Wilson in the symmetric case. The method of proof combines establishing a simpler special case of shifting 
one parameter by a factor of $q$ with using a co-cycle condition property of the transition matrices involved. Supporting computations use the
Noumi representation and are based on simple formulas for how some basic Hecke algebra elements act on natural  `almost symmetric' 
Laurent polynomials.
\end{abstract}

\ifJOLT
\else
\maketitle
\tableofcontents
\fi

\ifJOLT
\bc \small {\em It is a pleasure to have this paper be part of a tribute to Karl Heinrich Hofmann.} \ec
\else
\bc \small {\em It is a pleasure to dedicate this paper to Karl Heinrich Hofmann.} \ec
\fi
\section{Introduction \label{INTRO}}

The Askey scheme \cite{KLS} is a hierarchy of hypergeometric and $q$-hypergeometric  orthogonal poynomials, which includes many families of classical polynomials, including the Hermite, Laguerre, Gegenbauer, and Jacobi polynomials. Each family depends on a number of auxiliary parameters, and lower families in the hierarchy can be obtained from higher families by a suitable specialization of parameters. 


If $\cP=\{P_0,P_1,\ldots\}$ and $\cP'=\{P_0',P_1',\ldots\}$ are two sequences of polynomials satisfying $\deg(P_n)=\deg(P'_n)=n$, then one has an expansion of the form
\begin{equation}
\label{cnk}
P_n=\sum\nolimits_{m\le n} c_{m,n}P_m'.
\end{equation}
(The coefficients $c_{m,n}$ depend on the ordered bases $\cP$ and $\cP'$ and so might be more fully denoted by $ c_{m,n}^{\cP,\cP'}.$
However to keep our notation more spare, we usually make clear the bases by context and omit the superscripts.)

Now suppose $\cP$ and $\cP'$ are from the \emph{same} family in the Askey-Wilson scheme and differ in only \emph{one} auxiliary parameter. In this case one can sometimes obtain an explicit formula for the ``connection' coefficients'' 
$c_{m,n}$ 
in \eqref{cnk}  in terms of Pochammer symbols and their $q$-analogs, which are defined as follows:
\begin{align}
(a)_n = \prod\nolimits_{k=0}^{n-1} (a+k),\quad 
& (a_1,\ldots,a_r)_n= \prod\nolimits_{i=1}^{r} (a_i)_n, \\
(a|q)_n = \prod\nolimits_{k=0}^{n-1} (1-aq^k),\quad 
& (a_1,\ldots,a_r|q)_n= \prod\nolimits_{i=1}^{r} (a_i|q)_n.
\end{align} 
For example, one has the following classical result: 
\ifJOLT \begin{Theorem} \else \begin{thms} \fi\label{JACOBICONTHM}( 3.40  in \cite{INTREP})
\n If $\cP=\{P_n^{(\gamma,\beta)}(x)\}$ and $\cP'=\{P_n^{(\alpha,\beta)}(x)\}$ are two sequences of Jacobi polynomials differing in one parameter, then one has 
%
\bml
\label{JACOBICON}
P_n^{(\gamma,\beta)}(x)=\\
\sum_{k=0}^n \frac{ (\beta+k+1 )_{ n - k} ( \gamma-\alpha)_{ n - k}  (\beta+\gamma+n+1 )_{ k } ( 2k+\alpha+\beta+1) } {(\alpha+\beta+k+1 )_{n+1 } (n-k)! } P_k^{(\alpha,\beta)}(x).
\end{multline}
\ifJOLT \end{Theorem} \else \end{thms} \fi
Similarly equation (3.42) in  \cite{INTREP} tells us for Gegenbauer polynomials
\begin{equation}
\label{GEGENCON}
C_n^{\lambda}(x) =
 \sum_{k=0}^{\lb \frac{n}{2} \rb} \frac{ (\lambda - \nu )_{ k} ( \lambda)_{ n - k } ( n - 2k + \nu ) } { ( \nu )_{ n - k + 1} k!} C_{n-2k}^{\nu}(x). 
\end{equation}
(Equations (3.40) and (3.42) in \cite{INTREP} and elsewhere are often expressed in terms of Gamma functions, but the identity
$$
\Gamma(k+\alpha) = (\alpha)_k \Gamma(\alpha) \hspace{5mm} \text{ for } k \in \mathbb{N}
$$
readily leads to formulas \eqref{JACOBICON} and \eqref{GEGENCON}.)

A similar result holds for the Askey-Wilson polynomials $P_n(z;a,b,c,d|q)$, which are a 4 parameter $q$-hypergeometric family at the top of the Askey hierarchy.

\ifJOLT \begin{Theorem} \else \begin{thms} \fi\label{PTCFRM}(\cite{MEMOIRS},Askey-Wilson)
If $\cP= \{P_n(z;a,b,c,d|q)\}=\{P_n\}
\text{ and }$ \\
$\cP'= \{P_n(z;e,b,c,d|q)\}=\{P_n'\}$ are two sequences of Askey-Wilson polynomials differing in one parameter, then one has
\begin{eqnarray}
\label{AWCONSUM} P_n &=&\sum_{m\le n}  c_{m,n} (a, e ; b, c, d) P_m'  \hspace{10mm} \text{ where }\\ 
\label{AWCON} c_{m,n}(a, e ; b, c, d)&=&  \frac{ (q^{n-m+1} | q )_m(bc q^m, bd q^m, cd q^m, ae^{-1}  | q )_{n-m}   }  { (q  | q )_m (abcdq^{n+m-1}, bcdeq^{2m}  | q )_{n-m} } e^{n-m}. 
\end{eqnarray}
\ifJOLT \end{Theorem} \else \end{thms} \fi
Here we follow the notation of \cite{ASKVOL} and regard the Askey-Wilson polynomials as Laurent polynomials in $z$, symmetric under the inversion $z\mapsto z^{-1}$. This is related to the ordinary polynomial variable $x$ of Askey-Wilson \cite{MEMOIRS} by $2x= z +z^{-1}$. Our normalization of $P_n$ as a monic Laurent polynomial in $z$ is also different from \cite{MEMOIRS}. 

As explained in \cite{ASKVOL}, the Askey-Wilson polynomials admit nonsymmetric analogs $E_r=E_r(z; a,b,c,d|q)$, defined for any integer $r$, which are eigenfunctions of certain $q$-difference operators. While $P_n=P_n(z; a,b,c,d|q)$ can be obtained from $E_{\pm n}$ by a suitable symmetrization operator, the $E_r$ themselves are not symmetric under $z\mapsto z^{-1}$.  Also, while the $P_n$ depend symmetrically on $a,b,c \text{ and }d$, the $E_r$ only have the symmetries $a \leftrightarrow b$ and $c \leftrightarrow d$.

We define the ``zig-zag'' order $\prec$ on integers as follows: 
\begin{equation}
0 \prec -1 \prec 1 \prec -2 \prec 2 \prec \cdots.\label{zigzag}
\end{equation}
We say a Laurent polynomial $F(z)$ has degree 
an integer 
$r$ if there is a constant $c\ne0$ such that $F(z)-cz^r$ is in the span of $\{z^s, s\prec r\}$.
Such a Laurent polynomial is called monic if that leading coefficient $c$ is $1.$
 If $\cF=\{F_0,F_{-1},F_{1},\ldots \}$ and $\cF'=\{F'_0,F'_{-1},F'_{1},\ldots \}$  are two families of Laurent polynomials satisfying $\deg(F_r)=\deg(F_r')=r$, then one can again consider connection coefficients (e.g. denoted by $b_{r,s}$) such that
\begin{equation}
\label{crsE}
F_s=\sum\nolimits_{r\preceq s} b_{r,s}F_r'.
\end{equation}

The nonsymmetric Askey-Wilson polynomials satisfy the degree condition, and 
our main results are 
formulas 
of this form
giving the connection coefficients 
for two sequences of such polynomials differing in one parameter. 
In view of the $a \leftrightarrow b$ and $c \leftrightarrow d$ symmetries, there are only two distinct cases to consider: (1) the parameter $a$ is replaced by $e$, say, and (2) the parameter $c$ is replaced by $g$, say.


The first case
is a change of basis relationship  from basis $\{E_r(z; a,b,c,d|q)\}$ 
to basis $\{E_r(z; e,b,c,d|q)\}.$ 
We often refer to the 
corresponding matrix transforming components relative to these bases 
as a {\em transition matrix.} 

\ifJOLT \begin{Theorem} \else \begin{thms} \fi
\label{ETC}
Let $\cF=\{E_r(z;a,b,c,d|q)\}=\{E_r\}$ and $\cF'= \{E_r(z;e,b,c,d|q)\}=\{E_r'\}$ be sequences of nonsymmetric Askey-Wilson polynomials. Then 
$$
\ds E_s =\sum_{r \preceq s} d_{r,s}c_{|r|,|s|} E_r'
$$
where $c_{m,n}$ is the symmetric connection coefficient $c_{m,n}(a, e; b,c,d)$ as in \eqref{AWCON}, and
\[ \quad  d_{r,s}= \begin{cases}  
\ds  \frac{q^{s-r} (abcdq^{s+r-1}| q)_1}{(abcdq^{2s-1}| q)_1} & \text{ if }   r \ge 0, s\ge 0 \\ 
\ds   \frac{ (q^{-(r+s) }|q)_{ 1 } } {( q^{-s}, cdq^{-(s+1) }|q)_{ 1 }} & \text{ if }   r \ge 0, s< 0 \\ 
\ds  \frac{bcdeq^{s-(r+1)} (q^{-r}, cdq^{-(r+1)}, ae^{-1}q^{s+r }|q)_{1 } } { (abcdq^{2s-1}, bcdeq^{-(2r+1) }|q)_{ 1}}  & \text{ if }   r < 0, s\ge 0 \\ 
\ds  \frac{ (q^{-r}, cdq^{-(r+1)}, bcdeq^{-(r+s+1) }|q)_{1 } } { (q^{-s}, cdq^{-(s+1)}, bcdeq^{-(2r+1) }|q)_{ 1} }  & \text{ if }   r < 0, s< 0  \\
\end{cases} \]
\ifJOLT \end{Theorem} \else \end{thms} \fi

\n {\bf Remark:} It 
is 
worth noting that if one were to replace  $E_r(z; a,b,c,d|q)$ above by $\psi(r) E_r(z; a,b,c,d|q)$   then
there is some simplification in the resultant $d_{r,s}$ formulas when $\psi(r)$ is defined to be 
\begin{equation} \label{psi}
\psi(r)= \begin{cases}  
\hfil 1 & \text{if }   r \ge 0\\
(1-q^{-r})(1- cdq^{-(r+1)}) & \text{if }  r < 0. 
\end{cases}
\end{equation}
(This corresponds to  a change of normalization of the $E_r$ for $r<0$.)

We now consider the case where the parameter $c$ is replaced by $g$. 
\ifJOLT \begin{Theorem} \else \begin{thms} \fi
\label{ETC2}
Let $\cF=\{E_r(z;a,b,c,d|q)\}=\{E_r\}$ and $\cF''= \{E_r(z;a,b,g,d|q)\}=\{E_r''\}$ be sequences of nonsymmetric Askey-Wilson polynomials.
Then
$$
E_s =\sum_{r \preceq s}  d_{r,s}'' c_{|r|,|s|}'' E_r''  
$$
where
 \begin{enumerate}
\item $c_{m,n}''$ is the symmetric connection coefficient $c_{m,n}(c, g; a, b, d)$
\item $d_{r,s}''$ is $d^c_{r,s}(c, g; a,b,d )$ and
\end{enumerate}
\[ \quad d^c_{r,s}(c, g; a, b, d )= \begin{cases}  
\ds  \frac{ ( abq^s, abcdq^{r+s-1} | q)_1  }  { ( abq^r, abcdq^{2s-1} | q)_1 }  & \text{ if }   r \ge 0, s\ge 0 \\ 
\ds  -\frac{  abq^r ( q^{-r -s}  | q)_1 }  { ( q^{-s}, abq^r | q)_1  }   & \text{ if }   r \ge 0, s< 0 \\ 
\ds  -  \frac{ dq^{-r-1}g (q^{-r},  abq^s, cg^{-1} q^{s+r}  | q)_1 }  { ( abcd q^{2s-1}, abdgq^{-2r-1} | q)_1 }  & \text{ if }   r < 0, s\ge 0 \\  
\ds   \frac{ ( q^{-r} , abdg q^{ -r-s-1} | q)_1  }  { ( q^{-s} , abdg q^{-2r-1}  | q)_1 }  & \text{ if }   r < 0, s< 0.  \\
\end{cases} \]
\ifJOLT \end{Theorem} \else \end{thms} \fi

\n {\bf Acknowledgement:} The research of SS was partially supported by NSF grants DMS-$1939600$ and $2001537$, as well as Simons Foundation grant $509766.$ The research of BS was partially supported by Simons Foundation grant $89086.$

\bc \bf Motivation and Applications \ec 
For certain values of the parameters, formulas such as (\ref{JACOBICON}) in Theorem \ref{JACOBICONTHM} and (\ref{GEGENCON})  can be interpreted as describing  
aspects  
of 
 branching theorems for spherical
representations of  compact classical Lie groups. 
Spherical representations of compact groups are isomorphic to irreducible subrepresentations of 
 \[L^2(G/K) = \oplus_ \mu W_\mu\] for a compact symmetric space $G/K.$ 
Here $G$ acts on the left. (See  \cite{HELG2} and chapters 1 and 2 of \cite{MODSPHFCNS}.) Each $W_\mu$ contains a 
$K$-bi-invariant function 
 $f_\mu \in C^\infty(K,\ G /K)$. In  representation theory  $ f_\mu$ is referred  to as a spherical function for $G/K,$ but sometimes  in the literature the term zonal spherical function is used. Denoting by $T$ the one dimensional maximal  torus of a rank $1$ symmetric space, we have $G/K=KT$ by the maximal torus theorem for compact symmetric spaces and thus we can consider $f_\mu$ as a function on $T.$ 
As described on page $59$ of \cite{STOCHASTIC} (also e.g. page $65$ of \cite{ASKEYSIAM}, chapter $3$ of \cite{MODSPHFCNS}, or Table $1$ in \cite{JACCROSS}),
Gegenbauer and Jacobi polynomials can be viewed (up to scaling) as spherical functions on the rank $1$ symmetric spaces 
$S^n$ and 
$\mathbb{C}P^n.$

These observations extend to the other compact simply connected  rank $1$ symmetric spaces: the quaternionic projective space
and the Cayley plane. The classical polynomial parameter values associated with these geometric examples are also presented in
these references.
For the sphere $S^n,$ the parameter $\nu$ in the Gegenbauer polynomial $C_k^{\nu}(x)$ is $\ds \frac{n-1}{2}.$ For the complex projective space
$\mathbb{C}P^n,$ 
the parameters in the Jacobi polynomial  $P_k^{(\alpha,\beta)}$ are given by $\alpha=n-1$ and $\beta=0.$ 

\medskip
Recall how  branching of spherical representations from $G$ to $G'$ can give rise to connection coefficient formulas (\ref{GEGENCON}) and
(\ref{JACOBICON}). The $G' \subset G$ situations will be $SO(2n) \subset SO(2n+1)$ and $SU(n-1) \subset SU(n).$
We start with a spherical representation $\tau $ of $G,$ realized as a subrepresentation of $L^2(G/K)$.
We then want to consider  $G' \subset G$ with an associated rank $1$ symmetric space $G'/K'.$ Since they are both rank $1$ symmetric spaces,
their maximal tori in the sense of  symmetric spaces may be identified, so we may view 
\[G'/K'=K'T.\] 
Let $f_{\tau}$ denote the spherical function associated to $\tau.$
When we restrict $\tau$ to $G'$ we obtain $\oplus \tau_{\alpha}'.$
For each summand $\tau_\alpha'$ we have the corresponding spherical function $f_\alpha'$ and 
a relationship
\[ f_\tau(x) = \sum c_{\alpha} f_{\alpha}'(x) \]
where the $c_{\alpha}$ are non-negative constants. 
(Since spherical functions with respect to $G/K$ restricted to $T$ are multiples of matrix elements $<v,\tau(g)v>$ for $v$ a $K$-invariant vector of $\tau,$
the restricted functions  $f_{\alpha}'(x)$ are spherical functions with respect to $G'/K'.)$

For $SO(2n) \subset SO(2n+1)$ and $SU(n-1) \subset SU(n),$ this translates to the formulas (\ref{GEGENCON}) and (\ref{JACOBICON})  for  
Gegenbauer and Jacobi functions respectively.  So  $c_{\alpha} $ becomes 
a  connection coefficient in a formula  like
(\ref{GEGENCON}) and
(\ref{JACOBICON}).

The 
connection coefficient relations (\ref{GEGENCON}) and (\ref{JACOBICON}) can also be viewed as being related to calculating integrals
of a product of two such special polynomials with different parameters; as such one might for example use
Rodriguez' formula and integration by parts to determine them. 

In this paper we shall generalize these classical relations by purely 
combinatorial methods. Although there could also be relations to deformations of algebraic structures
such as double affine Hecke algebras via their generators and eigenfunctions, and also possibly to the
theory of spherical functions for certain quantum 
symmetric spaces,
we shall not go into these.

\bc \bf  Structure of the proof of Theorem \ref{ETC}\ec 
One natural approach to proving this theorem would be to use the vector-valued reformulation of the $E_r$ in \cite{VECVAL}, which extends earlier results in \cite{ZHED}. 

Our proof
here uses an interesting alternative. We start with proving the special case of $a$ and $e$ differing by a factor of $q.$ This is a $q$-shift analogous to a shift by $1$ in one of the classical integer parameters. Since transition functions describing change of basis matrices satisfy a natural co-cycle condition,
we can  establish the case of $a$ and $e$ differing by an integral power of $q$ by showing that our asserted expressions for the transition
functions also satisfy the co-cycle condition. And then, by observing that everything involved is given by rational functions agreeing at infinitely many values,
we obtain the theorem for arbitrary $a$ and $e.$

\pagebreak
\bc \bf The Co-cycle Condition in Detail and Proof Plan A \ec
The partition of powers of $z$ into non-negative vs. negative (corresponding to the same notions on the root system $C_1$ sitting inside the
affine root system $\tilde{C}_1$) gives a direct sum decomposition of Laurent polynomials
$$\mathcal{R} = \mathcal{R}^0 \oplus \mathcal{R}^1$$
where we can choose ordered bases
\ba
 E_0,E_1,E_2,\ldots  \text{ for } &  \mathcal{R}^0 \\
E_{-1},E_{-2},\ldots \text{ for } & \mathcal{R}^1.\\
\ea

We view connection coefficient relations like those given in Theorem \ref{ETC} as describing a change of basis in the space of Laurent polynomials, perhaps truncated
in degree, so as to reference a finite dimensional subspace. 

Unless otherwise specified, we will henceforth treat the common parameters $b,c,d,\text{ and } q$ as unchanged and drop them from the
argument lists of the $E_r$ (and associated connection coefficients.)
 For example $E_n(a)$ is a shorthand for $E_n(z; a,b,c,d|q).$

We use $\boldsymbol{\mathcal{T}}(a,e)$ for the `true' transition matrix from components relative to
the $\{E_r(a)\}$ basis to components relative the $\{E_r(e)\}$ basis. The notation $T(a,e)$ will refer to the transition matrix specified by the formulas
in Theorem \ref{ETC}.
Thus proving Theorem \ref{ETC} amounts to showing $\boldsymbol{\mathcal{T}}(a,e)=T(a,e).$

As described fully in Appendix A, this block decomposition and choice of ordered bases corresponds to $T(a,e)$ having the block decomposition
$$
T=\begin{bmatrix}
T ^{00} & T^{01} \\
T^{10} & T^{11} \\
\end{bmatrix}
$$
with

\ba
T^{00}=&\begin{bmatrix}
\tau_{0,0}  & \tau_{0,1} & \tau_{0,2}  &  \tau_{0,3}  & \ldots\\
0  & \tau_{1,1} & \tau_{1,2}  &  \tau_{1,3}  & \ldots \\
0  & 0 & \tau_{2,2}  & \tau_{2,3} & \ldots \\
0  & 0 & 0  & \ddots & \vdots \\
\end{bmatrix}
& 
T^{01}=&\begin{bmatrix}
\sigma_{0, -1}  & \sigma_{0,-2} & \sigma_{0,-3}  &  \sigma_{0,-4}  & \ldots\\
0  & \sigma_{1,-2} & \sigma_{1,-3}  &  \sigma_{1,-4}  & \ldots \\
0  & 0 & \sigma_{2,-3}  & \sigma_{2,-4} & \ldots \\
0  & 0 & 0  & \ddots & \vdots \\
\end{bmatrix}
\\
T^{10}=&
\begin{bmatrix}
0  & \sigma_{-1, 1} & \sigma_{-1, 2}  &  \sigma_{-1, 3}  & \ldots\\
0  & 0 & \sigma_{-2,2}  &  \sigma_{-2,3}  & \ldots \\
0  & 0 & 0  & \sigma_{-3,3} & \ldots \\
0  & 0 & 0  & \ddots & \vdots \\
\end{bmatrix}
& 
T^{11}=&\begin{bmatrix}
\tau_{-1, -1}  & \tau_{-1,-2} & \tau_{-1,-3}  &  \tau_{-1,-4}  & \ldots \\
0  & \tau_{-2,-2} & \tau_{-2,-3}  &  \tau_{-2,-4}  & \ldots \\
0  & 0 & \tau_{-3,-3}  & \tau_{-3,-4} & \ldots \\
0  & 0 & 0  & \ddots & \vdots \\
\end{bmatrix}.
\ea
Here $\tau_{r,s}$ and $\sigma_{r,s}$ (zero unless $r \preceq s$) are the products of the $c$'s and $d$'s defined by
\begin{align}
\label{TAUDEF} \tau_{r,s}&=&d_{r,s}c_{|r|,|s|} & \hspace{10mm} & \text{if } (r \ge 0 \text{ and } s \ge 0)  \text{ or } (r <  0 \text{ and } s < 0)\\
\label{SIGMADEF} \sigma_{r,s}&=& d_{r,s}c_{|r|,|s|} & \hspace{10mm} & \text{if } (r \ge 0 \text{ and } s < 0)  \text{ or } (r <  0 \text{ and } s > 0). 
\end{align}
We think of this transition matrix as acting on the left on column vectors of components relative to one basis and producing
a column vector of components relative to the other basis.

The `true' transition function $\boldsymbol{\mathcal{T}}(a,e)$ satisfies the co-cycle condition
\begin{equation}
\label{COCYC}
\boldsymbol{\mathcal{T}}(a,e) = \boldsymbol{\mathcal{T}}(f,e)\boldsymbol{\mathcal{T}}(a,f)
\end{equation}
since each side describes a valid way to go from  $a$-coordinates to $e$-coordinates. In particular this means a `discrete' co-cycle condition; namely for any non-negative integer $p:$
\begin{equation}
\label{COCYCAQ}
\boldsymbol{\mathcal{T}}(a,aq^{p+1}) = \boldsymbol{\mathcal{T}}(aq^p,aq^{p+1})\boldsymbol{\mathcal{T}}(a,aq^p).
\end{equation}

Our  proof of Theorem \ref{ETC}  has three steps which we refer to as
\begin{equation} 
\label{PLANA}
\text{ \em (PROOF PLAN A).  }
\end{equation}

\begin{enumerate}
\item Show that the entries of both $T(a,e)$ and  $\boldsymbol{\mathcal{T}}(a,e) $ are rational functions of $e$ with coefficients in the filed $\mathbb{Q}(a,b,c,d,q).$
\item Show $\boldsymbol{\mathcal{T}}(a,aq)=T(a,aq).$
\item Show $T$ also satisfies the discrete co-cycle condition $$T(a,aq^{p+1}) = T(aq^p,aq^{p+1})T(a,aq^p)$$ for any $p \in \mathbb{N}.$
\end{enumerate}
Since both $T(a,a)$ and $\boldsymbol{\mathcal{T}}(a,a)$ are the identity, and  equation (\ref{COCYCAQ}) says $\boldsymbol{\mathcal{T}}$ satisfies the discrete
co-cycle condition, it is immediate from 
parts $2$ and $3$ above 
that  $\boldsymbol{\mathcal{T}}(a,e)$ and $T(a,e)$ agree whenever $e=aq^p$ for non-negative $p.$  Now using 
part $1$ of Proof Plan A, 
we see that each entry of the two matrices is a rational function of $e$ agreeing with the other at infinitely many points. So they must agree
$\big($as rational functions with coefficients in $\mathbb{Q}(a,b,c,d,q) \big)$ for all $e.$

\bc \bf Structure of the Paper \ec 


We have already stated our main results above, explained how these kind of results can relate to branching theorems in
representation theory, and described the basic approach of the proof. 

The heart of the proof of 
Theorem \ref{ETC} are the last two steps of \eqref{PLANA}, Proof Plan A. These  are carried out in sections 
\ref{TAQ} and \ref{COCYCPRF}. We have written out the verification of these two steps in a detailed step-by-step way,
and so these two sections constitute about one third of the main body of our paper.

Section \ref{PREL} briefly recalls the basic double affine Hecke algebra (DAHA) and Noumi representation points of view about
the nonsymmetric Askey-Wilson polynomials. Much of what we use is based on the approach of \cite{ANNALS}, which was
specialized to the one variable case in \cite{ASKVOL}, and further enhanced (including some notational adjustments) in
\cite{SIGMA}. At certain points, we need a little more detail than was recorded in the statements of the theorems proved there,
and so explain how those come from this earlier work. 

The zig-zag order leads to a filtration of the Laurent polynomials and Section \ref{ALMSYM} exploits some aspects of this.
Some of the results are conveniently expressed in terms of what we call `almost symmetric' basis elements. These are based on Laurent 
polynomials which are either symmetric or skew-symmetric under one of the involutions $z \mapsto z^{-1}$  or  $z \mapsto qz^{-1}.$ 

The recursive description of the $\{E_r\},$ equations \eqref{SRCNX}, 
plays a key role. 
It immediately gives us the first step of \eqref{PLANA}, Proof Plan A.
The filtration properties of this `zig-zag recursion'
allow us to obtain explicit formulas for the three highest zig-zag degree terms of each $E_r.$ The coefficients of these terms with
respect to the appropriate almost symmetric basis elements are determined in the later parts of Section \ref{ALMSYM}. Our determination of the last of these, Theorem \ref{CDG2COEFFS}, is a little involved, and so we have also included a detailed step-by-step presentation
of the argument.

Properties of the filtration allow us to quickly express $\boldsymbol{\mathcal{T}}(a,aq)$ in terms of the  three highest zig-zag degree terms of the $E_r.$ Thus the later parts of Section \ref{ALMSYM} are exactly what we need to
carry out the second step of \eqref{PLANA}, Proof Plan A, in Section \ref{TAQ}.

It is natural to compare our nonsymmetric Askey-Wilson connection coefficient results (Theorems \ref{ETC} and \ref{ETC2})  to
the corresponding long known symmetric case, Theorem \ref{PTCFRM}. That is why those former statements are in terms of products 
like $d_{r,s}c_{|r|,|s|}.$ However, for carrying out the proof here, it is cumbersome to be constantly writing out these products
and so an alternate notation for these products is introduced in Section \ref{TAUSIGMA}.

Because of the natural $2 \times 2$ block structure of our transition matrices, verifying the co-cycle condition in Section \ref{COCYCPRF}
has four somewhat similar pieces. It also turns out that those verifications can be carried out  more simply by first simplifying
some ratios, and that is done at the beginning of Section \ref{COCYCPRF}.

Using our result, we also give a re-proof of Askey and Wilson's Theorem \ref{PTCFRM} in 
Section \ref{REPROOF} 
of our paper. Here the $d_{r,s}c_{|r|,|s|}$ representation of the connection coefficients in Theorem \ref{ETC} greatly simplifies the exposition.

While the considerations are very elementary, making explicit how our conventions lead to the precise matrices we use is
important for being able to verify correctness of our arguments. This is done in Section \ref{APPA}, Appendix A.

Section \ref{HAT}, Appendix B, a summary table, includes a few DAHA related formulas that we do not use in this paper, but which 
are of the same nature as ones we fully justify in the main body.

\ifLONG

The technique for proving Theorem \ref{ETC2} is extremely close to the one for Theorem \ref{ETC}, so we omit including the details
within the main body of this paper. 

However, in this longer version of our paper, the details (about a third of this longer form)  are available in Appendices C1 (Section \ref{SHIFTC1}), C2 (Section \ref{SHIFTC2}), and C3 (Section \ref{COCYCPRFC})
located just before the bibliography here.
\fi

We mention
one striking feature, however. In the proof of the discrete co-cycle condition for the shift-c case (Propositions \ref{T00PROP}, \ref{T01PROP}, \ref{T10PROP}, \text{ and }\ref{T11PROP} below in the shift-a case), the polynomials $p_2,$ whose vanishing in the last quarter of the proofs we are demonstrating, can in a natural
way be chosen to be {\em identical } to the $p_2$ of the shift-a case.

\section{Preliminaries \label{PREL}}


First we recall the double affine Hecke algebra (DAHA) point of view and the Noumi representation on Laurent polynomials.

Let $\mathcal{R}$ denote the Laurent polynomials in one variable $z$ with coefficients in a field such as $\mathbb{Q}(a,b,c,d,q).$ 

Let $\mathbb{F}$ be the field $\mathbb{Q}(q^{\frac{1}{2}}, t_0^{\frac{1}{2}}, t_1^{\frac{1}{2}},u_0^{\frac{1}{2}},u_1^{\frac{1}{2}})$. As described in  \cite{ANNALS} and \cite{ASKVOL}, the Noumi representation \cite{NOUMI} is a faithful representation
of a double affine Hecke algebra (DAHA)  $\mathcal{H}$ with coefficients in $\mathbb{F}.$  Here $\mathcal{H}$ is the $\mathbb{F}$-algebra with generators $T_0,T_1,U_0,U_1$
and relations
$$T_0 \sim t_0,\ T_1 \sim t_1,\ U_0 \sim u_0,\ U_1 \sim u_1, \text{ and } T_1T_0U_0U_1=q^{-\frac{1}{2}}$$
where the meaning of $F \sim f$ is $F-F^{-1}=f^{\frac{1}{2}} - f^{-\frac{1}{2}}.$ The scalars $q,a,b,c,d,$ $t_0,t_1,u_0,u_1$ are related by

\begin{equation}
\label{SCALARS}
\begin{array}{rrrr}
 \ds t_0= - \frac{cd}{q} & t_1=-ab & \ds u_0= - \frac{c}{d} & \ds u_1=-\frac{a}{b} \\
\ds a=t_1^{\frac{1}{2}} u_1^{\frac{1}{2}}  &\ds b = - t_1^{\frac{1}{2}} u_1^{-\frac{1}{2}}  & 
\ds c=q^{\frac{1}{2}} t_0^{\frac{1}{2}} u_0^{\frac{1}{2}}  &\ds d = -q^{\frac{1}{2}} t_0 ^{\frac{1}{2}} u_0 ^{-\frac{1}{2}} \\
\end{array}
\end{equation}
For our purposes, we mostly want to reason about Laurent polynomials with {\em rational} coefficients $\mathbb{Q}(a,b,c,d,q).$ So we 
work with elements $\widetilde{T}_0,\widetilde{T}_1, \widetilde{U}_0, \text{ and } \widetilde{U}_1$ (part of this described in \cite{SIGMA}) that are multiples
of the usual $T_0,T_1, U_0, \text{ and } U_1$:
$$
\widetilde{T}_0= t_0^{\frac{1}{2}}T_0 \hspace{10mm} 
\widetilde{T}_1= t_1^{\frac{1}{2}}T_1 \hspace{10mm} 
\widetilde{U}_0= (qt_0)^{\frac{1}{2}}U_0. \hspace{10mm} 
\widetilde{U}_1= (t_1)^{\frac{1}{2}}U_1. %
$$
These adjusted elements now satisfy the product relation $\widetilde{T}_1\widetilde{T}_0\widetilde{U}_0\widetilde{U}_1=t_0t_1$
and the  usual nicely symmetric DAHA inversion formulas, e.g.   
$$T_1^{-1}= T_1 -t_1^{\frac{1}{2}} +t_1^{-\frac{1}{2}},$$
translate (as in \cite{SIGMA}) to ones like 
$$
\lp\widetilde{T_1}+1\rp   \lp \widetilde{T_1}-t_1 \rp=0; \hspace{10mm} \text{i.e. }  \widetilde{T_1}^{-1} = \frac{\widetilde{T}_1}{t_1} -1+ \frac{1}{t_1}.
$$

With these adjustments, the Noumi representation of the DAHA $\mathcal{H}$ on a Laurent polynomial $f \in \mathcal{R}$ looks like:
\begin{eqnarray} 
\label{T0DEF} \lb \widetilde{T}_0 f\rb(z)&=& t_0 f(z)+ \ds \frac{(z-c)(z-d)}{(z^2-q)}\lb f\big(\frac{q}{z}\big)-f(z) \rb \\
\label{T1DEF}  \lb\widetilde{T}_1 f\rb(z)&=& t_1 f(z)+\ds \frac{(1-az)(1-bz)}{(1-z^2)} \lb f \big( \frac{1}{z}\big)-f(z) \rb \\
\label{U0DEF} \widetilde{U}_0 f &=&t_0    \widetilde{T}_0^{-1}X f   \ =  \lp \widetilde{T}_0X+(1-t_0)X\rp f
\end{eqnarray}
since $t_0 \widetilde{T}_0^{-1}=\widetilde{T}_0 + (1-t_0) $
and we use the multiplication operator
$$ \lb Xf \rb(z)= z\big( f(z) \big). $$
We call the operator $X$ because these DAHA preliminaries are often expressed in terms of a {\em Laurent series} variable $x$ while
we consistently use $z$ for the Laurent variable and its related $x=\frac{1}{2}(z+z^{-1})$ for the ordinary polynomial variable.

The nonsymmetric Askey-Wilson polynomials are eigenvectors of 
$$ \widetilde{Y}= \widetilde{T}_1 \widetilde{T}_0 \text{ (as well as of } Y=T_1T_0)$$
with $\widetilde{Y} E_r=\widetilde{\mu}_r E_r$ and 
$$
\widetilde{\mu}_r = \begin{cases} (t_0t_1)q^r &\ \text{if }r\ge 0 \\  q^r &\ \text{if }r<0.\end{cases}\\
$$
$\big(\widetilde{\mu}_r =(t_0t_1)^{\frac{1}{2}}q^{\overline{r}}\big),$ where 
\begin{equation}
\label{QBARDEF} q^{\overline{r}}= \begin{cases} (t_0t_1)^{\frac{1}{2}} q^r &\ \text{if }r\ge 0 \\  (t_0t_1)^{-\frac{1}{2}} q^r &\ \text{if }r<0\end{cases}
\end{equation}
is the notation used in \cite{ASKVOL} for the eigenvalues of $Y=T_1T_0.$ 


The definition of the $E_r$ on page 278 of \cite{ANNALS} is via creation operators $\mathcal{S}_0 \text{ and } \mathcal{S}_1.$ Up to normalization,
for $n \ge 0,$ $\mathcal{S}_0$ takes $E_n$ to $E_{-(n+1)}$ and $\mathcal{S}_1$ takes $E_{-(n+1)}$ to $E_{n+1}.$ 
This is the original recursive description of the nonsymmetric Askey-Wilson polynomials.

Theorem 4.1 on page 402 of
\cite{ASKVOL} makes this (in arbitrary rank)  more explicit. The proof of that theorem introduces a variant  $\mathcal{S'}_0$ of $\mathcal{S}_0$ which also takes
$E_{-(n+1)}$ to a multiple of $E_{n+1}.$ It is not explicitly noted in the theorem, but the proof makes clear that the formulas stated are simply the result of
applying the operators $\mathcal{S'}_0$ and $\mathcal{S}_1.$

Using our usual $\widetilde{U_0} \text{ and } \widetilde{T_1},$ the corresponding creation operators are
\begin{eqnarray*}
\widetilde{\mathcal{S'}}_0=& [\widetilde{Y},\widetilde{U}_0]  \hspace{10mm}  \widetilde{\mathcal{S}}_1=& [\widetilde{T}_1,\widetilde{Y}]    
\end{eqnarray*}
with $\widetilde{\mathcal{S'}}_0=t_0(qt_1)^{\frac{1}{2}}\mathcal{S'}_0$ and $\widetilde{\mathcal{S}}_1=(t_0)^{\frac{1}{2}} t_1\mathcal{S}_1.$

The normalization condition on $E_r$ is that the highest zig-zag degree term $z^r$ has coefficient $1;$ we  refer to this as $E_r$ being zig-zag monic.  We introduce the following $\widetilde{\zeta}_{i,r}$ notation for the explicit rescaling factors ($r$ of any sign):
$$
\widetilde{\mathcal{S'}}_0E_r= \widetilde{\zeta'}_{0,-(r+1)}   E_{-(r+1)} \hspace{10mm} \widetilde{\mathcal{S}}_1E_{-(r+1)}= \widetilde{\zeta}_{1,r+1}  E_{r+1} .
$$
For $n \ge 0$ we will determine $\widetilde{\zeta'}_{0,-(n+1)}$ and $\widetilde{\zeta}_{1,n}$ below in Proposition \ref{SCALE}. 

In those terms, Theorem 4.1 of \cite{ASKVOL} for the rank 1 case translates to:

\begin{eqnarray}
\label{SRCN0}
 E_{-(n+1)}=\lb  \widetilde{\zeta'}_{0,-(n+1)}  \rb^{-1} \widetilde{\mathcal{S'}}_0E_n& =t_0\lb \widetilde{\zeta'}_{0,-(n+1)}  \rb^{-1} \left[\left(\ds  \frac{ \widetilde{a}_{-(n+1)} } {t_0} \right) \widetilde{U}_0+\widetilde{b}_{-(n+1)} \right]E_{n} \quad&\\   
\label{SRCN1}
E_n= \lb \widetilde{\zeta}_{1,n}  \rb^{-1} \widetilde{\mathcal{S}}_1E_{-n}&=t_1\lb \widetilde{\zeta}_{1,n}  \rb^{-1} \left[\left(\ds  \frac{ \widetilde{c}_{n} } {t_1} \right) \widetilde{T}_1+\widetilde{d}_{n} \right]E_{-n} \quad (n \neq 0) &
\end{eqnarray} 
where
\begin{equation*}
\begin{aligned}
\ds\frac{\widetilde{a}_n}{t_0}= & \begin{cases}\ds \frac{(abcdq^{2n}| q)_1}{cdq^n} & \text{\ if\ } n\geq 0 \\
                                          \ds  - \frac{1-abcdq^{-2(n+1)}}{cdq^{-(n+1)}}&\text{\ if\ } n <0 \\ \end{cases} \\
\ds \widetilde{b}_n=& \begin{cases}  abcdq^{n}\left(\ds \frac{1}{c}+\frac{1}{d}\right)-(a+b)& \text{\ if\ } n\geq 0 \\
                                          q^{- (n  + 1)}\left(\ds \frac{1}{c}+\frac{1}{d}\right)-(a+b)&\text{\ if\ } n <0 \\ \end{cases} \\
\frac{\widetilde{c}_n}{t_1}=&  \begin{cases}\ds  0 & \text{\ if\ } n=0\\
                                         \ds  -  \frac{(abcdq^{2n-1}| q)_1}{abq^{n}} & \text{\ if\ } n>0 \\                                         
                                         \ds  \frac{1-abcdq^{-2n-1}}{abq^{-n}} &\text{\ if\ } n <0 \\ \end{cases} \\
\ds \widetilde{d}_n= &  \begin{cases} \ds \frac{q-abcd}{q} & \text{\ if\ } n=0\\
                                           \ds - \frac{ (abq^n| q)_1 + \left(abcdq^{n-1}| q\right)_1}{abq^{n}} & \text{\ if\ } n>0 \\                                           
                                           \ds  \frac{-cdq^{-n}(ab+1)+(cd+q)}{q}&\text{\ if\ } n <0 \\ \end{cases}  \\                                 
\end{aligned}
\end{equation*}
Theorem 1.2 in \cite{ASKVOL} stated this, but as mentioned in \cite{SIGMA}, has some typos. 
The formulas in equations (\ref{SRCN0}) and (\ref{SRCN1}) hold for any sign of $n$ above, but we emphasize, in this paper, the $n \ge 0$ cases because they,
together with $E_0=1,$ give a straightforward way to determine the $E_r$ inductively, for $r$ increasing in the zig-zag order sense. And it is easier to
prove, as we do later in Proposition \ref{SCALE}, the (also simpler) formulas for $\widetilde{\zeta'}_{0,-(n+1)}$ and $\widetilde{\zeta}_{1,-n}$ in those $n \ge 0$ cases.

\begin{description} 
\item[Comment on Some Awkward Looking Factors:] $\ds \frac{\widetilde{a}_n}{t_0}$ and $\ds \frac{\widetilde{c}_n}{t_1}$ as well as the $t_i$ factors in relating $\mathcal{S'}_0$ (respectively $\mathcal{S}_1$)  to $\ds \frac{\widetilde{a}_n}{t_0} + \widetilde{b}_n$ (respectively $\ds \frac{\widetilde{c}_n}{t_1} + \widetilde{d}_n$) arise   to make $\widetilde{a}_n$ and $\widetilde{c}_n$
differences of eigenvalues of $\widetilde{Y}$ just as $a_n$ and $c_n$ are differences of eigenvalues of $Y;$ e.g. $\widetilde{c}_n = \widetilde{\mu}_{-n} - \widetilde{\mu}_{n}$ in analogy to $c_n= q^{\overline{-n}} - q^{\overline{n}}$ (correcting a typo on page 397 of \cite{ASKVOL}.)
\item[More Detail on the Translation from \cite{ASKVOL}:] 
The starting point is Theorem 4.2 of \cite{ASKVOL} using the value $n=1$ ($1$ variable case) there. We now review the elements in the proof and
application of this theorem.

Besides DAHA identity manipulation, the proof is based on three things:
\begin{enumerate}
\item The relations among eigenvalues of $Y=T_1T_0$ that correspond to the intertwining identities
\begin{equation}
\label{INTERTWINE} Y\mathcal{S}_1=\mathcal{S}_1Y^{-1} \hspace{15mm} Y\mathcal{S}_0= q^{-1}\mathcal{S}_0 Y^{-1}  \text{ where } \mathcal{S}_0= q^{\frac{1}{2}} Y \mathcal{S'}_0.
\end{equation}
\item The creation operator definition, up to normalization, in \cite{ANNALS} of the $E_r.$
\item The fact that the creation operators $\mathcal{S}_0$ and $\mathcal{S}_1$ have squares which, when restricted to their natural invariant $2$-dimensional subspaces,  are multiples of the identity.
\end{enumerate}

The eigenvalue relations may either be viewed:
\begin{enumerate}
\item  As consequences of the formula (\ref{QBARDEF}) (above) for $q^{\overline{n}},$ originally given in \cite{ASKVOL}.
\item Or as applications of Theorem 5.1 in \cite{ANNALS}. In particular, in Theorem 5.1, one can choose 
$\tilde{\nu}= 1+0\delta \in \mathbb{Z} \times  \mathbb{Z} \delta$
satisfying
\begin{eqnarray*}
s_0(1) = -1 -\delta & & s_1(1)=-1.
\end{eqnarray*}
\end{enumerate}
The statements about $\mathcal{S}_0^2$ and $\mathcal{S}_1^2$ are Corollary 5.2 for $n=1$ in \cite{ANNALS}. Here one has to keep in
mind that the $E_r$ are eigenvectors of $Y=T_1T_0.$
The following table summarizes many relationships between the original objects in \cite{ANNALS}, \cite{ASKVOL} and our notation here.

\bc
\begin{tabular}{|c|c|c|c|c|} \hline
$\widetilde{T}_0=t_0^{\frac{1}{2}} T_0$ & $\widetilde{T}_1=t_1^{\frac{1}{2}} T_1$ & $\widetilde{Y}=(t_0t_1)^{\frac{1}{2}} Y$\\ \hline
$U_0=q^{-\frac{1}{2}} T_0^{-1}X$ & $\widetilde{U}_0=t_0 \widetilde{T}_0^{-1} X$ & $\widetilde{U}_0 = (qt_0)^{\frac{1}{2}} U_0$\\ \hline
$U_1=(T_1X)^{-1}$ & $\widetilde{U}_1=t_1(\widetilde{T}_1X)^{-1}$ & $\widetilde{U}_1 = t_1^{\frac{1}{2}} U_1$\\ \hline
$T_1T_0U_0U_1=q^{-\frac{1}{2}}$ & $\widetilde{T}_1\widetilde{T}_0\widetilde{U}_0\widetilde{U}_1=t_0t_1$ & \\ \hline
$\mathcal{S}_1=[T_1,Y] $ & $\widetilde{\mathcal{S}}_1= [\widetilde{T}_1,\widetilde{Y}]$ &$\widetilde{\mathcal{S}}_1= t_0^{\frac{1}{2}}t_1\mathcal{S}_1$ \\ \hline
$\mathcal{S'}_0=[Y,U_0] $ & $\widetilde{\mathcal{S}}_0'= [\widetilde{Y},\widetilde{U}_0]$ &$\widetilde{\mathcal{S}}_0'= t_0(qt_1)^{\frac{1}{2}}\mathcal{S'}_0$ \\ \hline
& & $\mathcal{S'}_0 E_{-(n+1)}= \lp a_{n}U_0 + b_{n} \rp E_{-(n+1)}$ \\ \hline
& & $\mathcal{S}_1 E_{-n}= \lp c_{n}T_1 + d_{n} \rp E_{-n}$  \\ \hline
$\mathcal{S}_0=[Y,U_1^{-1}] $ & $\widetilde{\mathcal{S}}_0= [\widetilde{Y},\widetilde{U}_1^{-1}]$ &$\widetilde{\mathcal{S}}_0= (t_0)^{\frac{1}{2}}\mathcal{S}_0$ \\ \hline
$\ds \frac{\widetilde{a}_n}{t_0}=\lp\frac{t_1}{t_0}\rp^{\frac{1}{2}}a_n$& $\widetilde{b}_n=(qt_1)^{\frac{1}{2}}b_n$ & $\ds \frac{\widetilde{a}_n}{t_0}\widetilde{U}_0+\widetilde{b}_n =(qt_1)^{\frac{1}{2}}\lp a_nU_0+b_n \rp$ \\ \hline
$\ds \frac{\widetilde{c}_n}{t_1}=\lp\frac{t_0}{t_1}\rp^{\frac{1}{2}}c_n$ &$\widetilde{d}_n=t_0^{\frac{1}{2}}d_n$ & $\ds \frac{\widetilde{c}_n}{t_1}\widetilde{T}_1+\widetilde{d}_n =t_0^{\frac{1}{2}}\lp c_nT_1+d_n \rp$ \\ \hline
\end{tabular}
\ec 

To get the $a_n$ and $b_n$ (as in the above table)  from Theorem 4.2 in \cite{ASKVOL} (and then the asserted $\widetilde{a}_n \text{ and } \widetilde{b}_n),$  use $\lambda=\lambda_1=n, \mu=\mu_1=-(n+1).$ Here the $b_n$ comes from the statement about $c_0$ in Theorem 4.2.

To get the $c_n$ and $d_n$ (as in the above table)  from Theorem 4.2 in \cite{ASKVOL} (and then the asserted $\widetilde{c}_n \text{ and } \widetilde{d}_n),$  use $\lambda=\lambda_1=n,\mu = \mu_1=- n$  Here the $d_n$ comes from the statement about (a different)  $c_n$  for $n=1$ in Theorem 4.2.

Expressing everything in terms of $a,b,c,d,q,\text{ and } n$ (of any sign) via equations (\ref{SCALARS}) gives the asserted formulas for 
$\widetilde{a}_n,\widetilde{b}_n,\widetilde{c}_n, \text{ and } \widetilde{d}_n.$ 

\end{description}
\section{Almost Symmetric Bases and their Applications\label{ALMSYM}}

For any integer $n,$ let $\mathcal{R}_n$ denote the elements of the Laurent polynomials $\mathcal{R}$ of zig-zag degree at most $n$ in the zig-zag order.
This gives rise to a filtration
$$
\mathcal{R}_0 \subset  \mathcal{R}_{-1} \subset \mathcal{R}_1 \subset \mathcal{R}_{-2} \subset \mathcal{R}_2 \subset \ldots 
$$
with both $E_n$ and $z^n$ projecting to the same nonzero generator of the $1$-dimensionsal $\mathcal{R}_n/\mathcal{R}_{n-}.$ Here $n-$ denotes
the  predecessor of $n$ in the zig-zag order.

In thinking about the operators $\widetilde{T}_i,$ it is natural to relate them to skew-symmetrization. If we denote the involutions by $s_i:\mathcal{R} \to \mathcal{R},$
namely 
\begin{eqnarray*}
\lb s_1(f) \rb(z)=f\left(z^{-1} \right) &\hspace{15mm} &\lb s_0(f) \rb(z)=f\left(qz^{-1} \right)
\end{eqnarray*}
and the corresponding skew (respectively q-skew) symmetrizations by
$$
\Lambda_1(f)= \frac{1}{2} (1- s_1)\left(f\right)  \hspace{15mm} \text{(respectively } \Lambda_0(f)= \frac{1}{2} (1- s_0)\left(f\right) )
$$
then $\Lambda_1$ (respectively $\Lambda_0$) has eigenfunctions which we will denote for $n \ge 0$ as follows:
\ba
\text{Eigenvalue } 1: \ \ \ f_n =& (z+z^{-1})^n \\
\hspace{5mm} & \text{\big(respectively } f_{n\_\text{q}}= (z+qz^{-1})^n.\big) \\
\text{Eigenvalue } -1: \ \ \ g_n =&(z-z^{-1}) (z+z^{-1})^{n-1} \\
\hspace{5mm} & \text{\big(respectively } g_{n\_\text{q}}= (z-qz^{-1})(z+qz^{-1})^{n-1}.\big) 
\ea
Unfortunately neither  $f_{n+1}$  nor $g_{n+1}$ belong to $\mathcal{R}_{-(n+1)},$ but their difference does. So for $n \ge 0$ we define
$$
h_{n+1}= z^{-1} (z+z^{-1})^n \hspace{15mm}  \text{(respectively } h_{n+1,\_\text{q}}= qz^{-1} (z+qz^{-1})^n.)
$$
Now 
$ \{ f_0,h_1,f_1,h_2,f_2,\ldots\}$ and $ \{ f_{0\_\text{q}},h_{1\_\text{q}},f_{1\_\text{q}},h_{2\_\text{q}},f_{2\_\text{q}},\ldots\}$ both
form bases compatible with the filtration
$$
\mathcal{R}_0 \subset  \mathcal{R}_{-1} \subset \mathcal{R}_1 \subset\ldots \ .
$$
We refer to these as the {\em almost symmetric} and {\em almost q-symmetric} bases respectively.

The operator $\widetilde{T}_0$ involves division by $z^2-q.$ Easy argument shows that for any Laurent polynomial $f$, its q-skew
symmetrization $\Lambda_0(f)$ is divisible by $z^2-q.$ To understand the operators $\widetilde{T}_0,\widetilde{U}_0\text{ and } \widetilde{T}_1$  more fully,
we define Laurent polynomials $f_{\text{skew}},f_{\text{q\_skew}},f_{\text{skew\_rdcd}}, f_{\text{q\_skew\_rdcd}}, f_{\text{sym}},\text{ and } f_{\text{q\_sym}}, $ for any Laurent polynomial $f$ by:
\begin{eqnarray*}
f_{\text{skew}} =  \Lambda_1(f) &\hspace{15mm}& f_{\text{q\_skew}} =  \Lambda_0(f) \\
 f_{\text{skew}} (z) = (z^2-1) \lb f_{\text{skew\_rdcd}}(z) \rb &\hspace{15mm} &f_{\text{q\_skew}} (z) = (z^2-q) \lb f_{\text{q\_skew\_rdcd}}(z) \rb \\
 f_{\text{sym}} = \frac{1}{2} (1+ s_1)\left(f\right) &  & f_{\text{q\_sym}} = \frac{1}{2} (1+ s_0)\left(f\right).
\end{eqnarray*} 
Thus
$$
f= f_{\text{sym}} +f_{\text{skew}} = f_{\text{q\_sym}} +f_{\text{\_skew}}.
$$
In terms of these we have:
\ifJOLT \begin{Proposition} \else \begin{props} \fi
\label{U_0_T_1}
\

\begin{equation}
\label{U0A}
\begin{split}
\widetilde{U}_0f =& \left(\frac{(c+d) z - cd } {z } \right)f -2q\left(\frac{(z-c)(z-d)} {z}\right) f_{\text{q\_skew\_rdcd}}\\
=& - 2qz  f_{\text{q\_skew\_rdcd}} +\left[ c+d+  \left(\frac{cd } {2q } \right)\left(z - \frac{q} {z }\right)\right. \\
 &  \hspace{15mm}  \left.  - \left(\frac{cd } {2q }\right) \left(z + \frac{q} {z }\right)\right]\left(f + 2q f_{\text{q\_skew\_rdcd}}\right) \\
\end{split}
\end{equation}
\begin{equation}
\label{T1A}
\begin{split}
 \widetilde{T}_1 f =& -ab f + 2\left[(1-az)(1-bz)\right]   f_{\text{skew\_rdcd}}\\
 = &-ab f + \left[ -2(a+b) +\left(ab+1 \right) \left(z+\frac{1}{z} \right) + \left(ab - 1 \right) \left(z -\frac{1}{z} \right) \right] \left( zf_{\text{skew\_rdcd}}\right)
\end{split}
\end{equation}
\ifJOLT \end{Proposition} \else \end{props} \fi
\begin{proof}
To establish (\ref{U0A}), using the definitions (\ref{T0DEF})  and (\ref{U0DEF}), we start with
\begin{multline*}
\Big( \widetilde{U}_0 f \Big)(z) = \Big( \widetilde{T}_0 + (1-t_0) \Big) \Big(  z \lb f(z) \rb \Big) \\
=(1-t_0) z \lb f(z) \rb + t_0 \Big( z \lb f(z) \rb \Big) +\Big(\frac{(z-c)(z-d) } {z^2-q } \Big) \Big( \lp \frac{q}{z} \rp \lb f \lp\frac{q}{z} \rp \rb -z \lb f(z) \rb \Big) \\
\end{multline*}
Since 
$$
\frac{1}{z} - \frac{z}{z^2-q} = -\frac{q}{z(z^2-q)},
$$
this gives
\bmlst
(\widetilde{U}_0f)(z)  =zf(z) - (z-c)(z-d) \left[ \frac{1}{z} f\left( \frac{q}{z}\right)\right. \ -  \frac{z}{z^2-q} \left(\left. f\left( \frac{q}{z}\right)-f(z)\right)\right] \\
= z\left[ f_{\text{q\_sym}}(z) + f_{\text{q\_skew}}(z)\right]- (z-c)(z-d)\left[ \frac{1}{z}\left( f_{\text{q\_sym}}(z) - f_{\text{q\_skew}}(z)\right)\right.\\
\left. -\frac{z}{z^2-q}\left( - 2 f_{\text{q\_skew}}(z)\right) \right]
\end{multline*}
\bmlst
\hphantom{(\widetilde{U}_0f)(z)}= z\left[ f_{\text{q\_sym}}(z) + f_{\text{q\_skew}}(z)\right]- \frac{(z-c)(z-d)} {z}\Big[ [ \left( f_{\text{q\_sym}}(z) - f_{\text{q\_skew}}(z)\right)\\ 
\left. + \frac{2z^2}{z^2-q}\left( f_{\text{q\_skew}}(z)\right) \right]\\
= z\left[ f_{\text{q\_sym}}(z) + f_{\text{q\_skew}}(z)\right]- \frac{(z-c)(z-d)} {z} \Big[  f_{\text{q\_sym}}(z)  +\frac{z^2+q}{z^2 - q} \left(f_{\text{q\_skew}}(z)\right)\Big]\\
\end{multline*}
\bmlst
\hphantom{(\widetilde{U}_0f)(z) } = z\left[ f_{\text{q\_sym}}(z) + f_{\text{q\_skew}}(z)\right]- \frac{(z-c)(z-d)} {z}\Big[  f_{\text{q\_sym}}(z) + f_{\text{q\_skew}}(z) \\ 
\left. + \frac{2q}{z^2-q}\left( f_{\text{q\_skew}}(z)\right) \right]
\end{multline*}
\bmlst
\hphantom{(\widetilde{U}_0f)(z) }=\left(z - \frac{(z-c)(z-d)} {z} \right)f  - \left(2q \frac{(z-c)(z-d)} {z}\right) f_{\text{q\_skew\_rdcd}}(z)\\
=\left( \frac{(c+d)z -cd} {z} \right)f  - \left(2q \frac{(z-c)(z-d)} {z}\right) f_{\text{q\_skew\_rdcd}}(z)
\end{multline*}
The proof of the second form of (\ref{U0A}) comes from further observing
$$
 \frac{(c+d) z - cd } {z }  = c+d +  \left(\frac{cd } {2q } \right)\left(z - \frac{q} {z }\right) - \left(\frac{cd } {2q }\right) \left(z + \frac{q} {z }\right)
$$
and
$$
\frac{(z-c)(z-d)} {z}= z - \left(\frac{(c+d) z - cd } {z } \right).
$$

To establish the first form of  (\ref{T1A}), using the definition (\ref{T1DEF}), we start with
\bmlst
 \lp \widetilde{T}_1 f \rp(z)= t_1\lb  f(z) \rb +\lp \ds \frac{(1-az)(1-bz)}{(1-z^2)} \rp \lp f \big( \frac{1}{z}\big)-f(z) \rp \\
 =-ab \lb  f(z) \rb +\lp \ds \frac{(1-az)(1-bz)}{(1-z^2)} \rp \lb -2 f_{\text{\_skew}}(z) \rb \\
  =-ab \lb  f(z) \rb +\lp \ds \frac{(1-az)(1-bz)}{(1-z^2)} \rp \lp -2(z^2-1) \rp \lb f_{\text{\_skew\_rdcd}}(z) \rb \\
    =-ab \lb  f(z) \rb +2(1-az)(1-bz)  \lb f_{\text{\_skew\_rdcd}}(z) \rb 
 \end{multline*}
 which is the first form of (\ref{T1A}).
 
To get the second form, we use
$$
2\frac{(1-az)(1-bz)}{z}=(ab+1)\Big(z + \frac{1}{z} \Big) + (ab-1) \Big(z - \frac{1}{z} \Big) -2(a+b).
$$
\end{proof}

An immediate corollary is that $\widetilde{U}_0$ and $\widetilde{T}_1$ behave very nicely on our filtered bases:
\ifJOLT \begin{Corollary} \else \begin{cors} \fi
\label{T_U_BS}
$$\begin{aligned}
\widetilde{U}_0(f_{n\_\text{q}})&= &  -\left(\frac{cd}{q}\right) h_{n+1,\_\text{q}}& + (c+d) f_{n\_\text{q}}\\
\widetilde{U}_0(h_{n\_\text{q}})& = &q f_{n-1,\_\text{q}} \\
\widetilde{T}_1(f_n) &=& -ab f_n  \\
\widetilde{T}_1(h_n) &=& -ab f_n - h_n &+(a+b) f_{n -1} \\
\end{aligned}$$
\ifJOLT \end{Corollary} \else \end{cors} \fi
\begin{proof}
\begin{enumerate}
\item For $f=f_{n\text{\_q}}=\ds \left(z + \frac{q} {z }\right)^n,$ we have $f_{\text{q\_skew\_rdcd}}=0.$ So
\bmlst
\widetilde{U}_0(f)=\Big((c+d) +  \left(\frac{cd } {2q } \right)\left(z - \frac{q} {z }\right) - \left(\frac{cd } {2q }\right) \left(z + \frac{q} {z }\right) \Big)f_{n\text{\_q}} \\
=\lp (c+d) - \lp \frac{2q}{z} \rp \lp \frac{cd } {2q }\rp \rp f_{n\text{\_q}}\\
= (c+d) f_{n\text{\_q}} -  \lp \frac{cd } {q }\rp h_{n+1,\text{\_q}}.\\
\end{multline*}
\item For $f=h_{n\text{\_q}}= \ds \lp \frac{q}{z} \rp \left(z + \frac{q} {z }\right)^{n-1},$
$$
f_{\text{q\_skew}} = \ds \lp  \frac{1}{2} \rp  \lp \frac{q}{z}  - z\rp \left(z + \frac{q} {z }\right)^{n-1} = \ds- \lp \frac{1}{2z} \rp (z^2-q) \left(z + \frac{q} {z }\right)^{n-1}
$$
So
$$
f_{\text{q\_skew\_rdcd}} = \ds  - \lp \frac{1}{2z} \rp \left(z + \frac{q} {z }\right)^{n-1}=  \ds  - \lp \frac{1}{2q} \rp h_{n\text{\_q}} = \ds - \frac{f}{2q}.
$$
Thus 
$$
f+2q f_{\text{q\_skew\_rdcd}} =0 \text{ and } z f_{\text{q\_skew\_rdcd}}  = \ds - \lp  \frac{1}{2} \rp  \left(z + \frac{q} {z }\right)^{n-1}=  \ds - \lp  \frac{1}{2} \rp f_{n-1,\text{\_q}}.
$$
This gives us
$$
\widetilde{U}_0(h_{n\text{\_q}})= -2q \lp  -\frac{1}{2} f_{n-1,\text{\_q}}\rp +0=qf_{n-1,\text{\_q}}.
$$
\item For $f=f_n= \ds \lp z + \frac{1}{z} \rp^n,$ we have $f_{\text{skew\_rdcd}} =0.$ So by (\ref{T1A}),
$$
\widetilde{T}_1 f_n=-ab f_n.
$$
\item For $f=h_n= \ds\lp \frac{1}{z} \rp \lp z + \frac{1}{z} \rp^n,$
$$
f_{\text{skew}} =\ds \lp \frac{1}{2} \rp \lp \frac{1}{z}-z \rp \lp z + \frac{1}{z} \rp^{n-1}= \ds -\lp \frac{1}{2z} \rp \big(z^2-1\big) \lp z + \frac{1}{z} \rp^{n-1}.
$$
So 
$$
f_{\text{skew\_rdcd}} = \ds -\lp \frac{1}{2z} \rp \lp z + \frac{1}{z} \rp^{n-1}= \ds - \frac{h_n}{2}  \text{ and } zf_{\text{skew\_rdcd}}=\ds - \frac{f_{n-1}}{2} .
$$
Thus by  (\ref{T1A}),
\bmlst
\widetilde{T}_1 h_n = -(ab ) h_n  + \Big\{ -2(a+b) +2ab \lp z + \frac{1}{z} \rp -2(ab-1) \lp\frac{1}{z}  \rp \Big\} \lp- \frac{f_{n-1}}{2} \rp \\
= -(ab ) h_n  + (a+b) f_{n-1} -(ab) f_n + (ab-1) h_n \\
= -(ab) f_n -h_n +  (a+b) f_{n-1} 
\end{multline*}

\end{enumerate}
\end{proof}

\n {\bf Remark:} Corollary \ref{T_U_BS} shows, for $n \ge 0:$
\begin{enumerate}
\item $\widetilde{U}_0$ maps $\mathcal{R}_n$ to $\mathcal{R}_{-(n+1)}$ and $\mathcal{R}_{-(n+1)}$ to $\mathcal{R}_{-(n+1)}.$ 
\item $\widetilde{T}_1$ maps $\mathcal{R}_{-(n+1)}$ to $\mathcal{R}_{n+1}$ and $\mathcal{R}_{n}$ to $\mathcal{R}_{n}.$
\end{enumerate}

Because of Corollary \ref{T_U_BS}  and the recursion (\ref{SRCN0}) and (\ref{SRCN1}), we will need to convert between the almost symmetric bases in low zig-zag
co-degree. 
\ifJOLT \begin{Proposition} \else \begin{props} \fi
\label{COB}
\begin{eqnarray*}
f_n = f_{n\text{\_q}} + (q^{-n}-1)h_{n\text{\_q}} \text{ mod } \mathcal{R}_{-(n-1)} & &
h_n = q^{-n}h_{n\text{\_q}} \text{ mod } \mathcal{R}_{-(n-1)} \\ 
h_{n+1\text{\_q}} = q^{n+1}h_{n+1} \text{ mod } \mathcal{R}_{n-1} & & 
f_{n\text{\_q}} = f_{n} + (q^{n}-1)h_{n}  \text{ mod } \mathcal{R}_{n-1}
\end{eqnarray*}
\ifJOLT \end{Proposition} \else \end{props} \fi
\begin{proof}
\begin{equation*}
h_n = z^{-1} (z+z^{-1})^{n-1} = z^{-n}= q^{-n}h_{n\text{\_q}}  \text{ mod } \mathcal{R}_{-(n-1)}  
\end{equation*}
\bmlst
f_n= (z+z^{-1})^n = z^n + z^{-n} = z^n + q^nz^{-n} + (1-q^n) z^{-n}  \\ = z^n + q^nz^{-n} +  (q^{-n}-1) q^nz^{-n}\\ = f_{n\text{\_q}} + (q^{-n}-1)h_{n\text{\_q}}  \text{ mod } \mathcal{R}_{-(n-1)}  \\
h_{n+1\text{\_q}}= qz^{-1} (z+qz^{-1})^n=q^{n+1}z^{-(n+1)}= q^{n+1}h_{n+1}  \text{ mod } \mathcal{R}_{n-1} \\
f_{n+1\text{\_q}}=  (z+qz^{-1})^{n+1}=z^{n+1} +q^{n+1} z^{-(n+1)} \\=z^{n+1} + z^{-(n+1)} +(q^{n+1}-1)z^{-(n+1)}  \\ = f_{n+1} + (q^{n+1}-1)h_{n+1}   \text{ mod } \mathcal{R}_{n-1}
\end{multline*}

\end{proof}

The almost symmetric bases allow us to make the rescalings involved in the recursive computation of the $E_n$ explicit. This also makes the rationality of
the formulas transparent.

Since the recursion in equations (\ref{SRCN0}) and (\ref{SRCN1}) is in terms of $\widetilde{U}_0 \text{ and } \widetilde{T}_1,$ Corollary \ref{T_U_BS} makes an expansion
of the $E_r$ in terms of the  almost symmetric bases natural. For $n \ge 0,$ we will use the notation (remembering e.g. $h_{m+1} \in \mathcal{R}_{-(m+1)}$):
\begin{eqnarray*}
E_n & =& \sum_{m=0}^n \lambda_{m,n} f_m +  \sum_{m=0}^{n-1} \mu_{-(m+1), n} h_{m+1} \\
E_n &=& \sum_{m=0}^n \lambda_{m,n\text{\_q}} f_{m\text{\_q}} +  \sum_{m=0}^{n-1} \mu_{-(m+1), n} h_{m+1,\text{\_q}} \\
E_{-(n+1)} &= & \sum_{m=0}^n \lambda_{-(m+1),-(n+1)} h_{m+1} +  \sum_{m=0}^{n} \mu_{m, -(n+1)} f_{m} \\
E_{-(n+1)}& = &\sum_{m=0}^n \lambda_{-(m+1),-(n+1)\text{\_q}} h_{m+1\text{\_q}} +  \sum_{m=0}^{n} \mu_{m, -(n+1)\text{\_q}} f_{m\text{\_q}}.
\end{eqnarray*}

Using Corollary \ref{COB}, we quickly see
\ifJOLT \begin{Proposition} \else \begin{props} \fi
\label{TOFROMQ}
For $n \ge 0:$
\begin{eqnarray*}
\lambda_{-(n+1),-(n+1)}&=&\lambda_{-(n+1),-(n+1)\text{\_q}} q^{n+1} \\
 \mu_{n,-(n+1)} &=& \mu_{n,-(n+1)\text{\_q}}\\
\lambda_{-n,-(n+1)} &=&  \mu_{n,-(n+1)\text{\_q}} (q^{n}-1)  + \lambda_{-n,-(n+1)\text{\_q}} q^{n} \ ( n \ge1) 
\end{eqnarray*}
\begin{eqnarray*}
\lambda_{n,n\text{\_q}}& = &\lambda_{n,n}\\
\mu_{-(n+1),n+1\text{\_q}} &=& \lambda_{n+1,n+1}(q^{-(n+1)}-1)+\mu_{-(n+1),n+1}q^{-(n+1)} \\
{} \lambda_{n,n+1\text{\_q}} &=& \lambda_{n,n+1}. 
\end{eqnarray*}
\ifJOLT \end{Proposition} \else \end{props} \fi
\begin{proof}
Modulo $\mathcal{R}_{n-1}:$
\begin{multline*}
E_{-(n+1)}= \lambda_{-(n+1),-(n+1)\text{\_q}} h_{n+1,\text{\_q}}+ \mu_{n,-(n+1)\text{\_q}} f_{n,\text{\_q}}+ \lambda_{-n,-(n+1)\text{\_q}} h_{n,\text{\_q}} \\
               = \lambda_{-(n+1),-(n+1)\text{\_q}} \lb q^{n+1}h_{n+1}  \rb + \mu_{n,-(n+1)\text{\_q}}\lb f_{n} + (q^{n}-1)h_{n} \rb \\
               + \lambda_{-n,-(n+1)\text{\_q}} \lb q^{n}h_{n} \rb \\
               =\lb\lambda_{-(n+1),-(n+1)\text{\_q}}  q^{n+1} \rb h_{n+1} + \lb \mu_{n,-(n+1)\text{\_q}}  \rb f_n \\
               + \lb   \mu_{n,-(n+1)\text{\_q}} (q^{n}-1)  + \lambda_{-n,-(n+1)\text{\_q}} q^{n} \rb h_n 
\end{multline*}
from which we read off the first three statements.

Modulo $\mathcal{R}_{-(n-1)}:$
\begin{multline*}
E_n = \lambda_{n,n} f_n +   \mu_{-n,n} h_n +  \lambda_{n-1,n} f_{n-1} \\
       = \lambda_{n,n} \lb  f_{n\text{\_q}} + (q^{-n}-1)h_{n\text{\_q}}\rb  +   \mu_{-n,n} \lb q^{-n}h_{n\text{\_q}} \rb +  \lambda_{n-1,n} \lb f_{n-1,\text{\_q}} \rb \\
       = \lb  \lambda_{n,n}  \rb f_{n\text{\_q}} +  \lb \lambda_{n,n}  (q^{-n}-1)+  \mu_{-n,n}  q^{-n} \rb h_{n\text{\_q}} + \lb   \lambda_{n-1,n} \rb f_{n-1,\text{\_q}} 
\end{multline*}
from which we read off the last three statements.

\end{proof}

The zig-zag monic condition on the $E_r$  immediately implies
\ifJOLT \begin{Proposition} \else \begin{props} \fi
\label{LNN}
For any $n \ge 0:$
\begin{eqnarray*}
\lambda_{n,n} =1& \hspace{15mm} & \lambda_{n,n\text{\_q}} =  1 \hspace{10mm}\\
\lambda_{-(n+1),-(n+1)} = 1 &\hspace{15mm}& \lambda_{-(n+1),-(n+1)\text{\_q}}= q^{-(n+1)}.
\end{eqnarray*}
\ifJOLT \end{Proposition} \else \end{props} \fi
For our purpose, we need the explicit scaling factors arising in the zig-zag increasing cases of the creation operators acting on $E_r.$
We use the following notation for the exact coefficients:
\begin{equation}
\label{SRCNX}
\begin{aligned}
 {\rm\  (E_{negative} \ case)\ }\ \ \ & E_{-(n+1)}&= \ \ \ &\left( { \hat{a}_{-(n+1)} } \widetilde{U}_0+\hat{b}_{-(n+1)} \right)E_{n}&\ \\   
{\rm\  (E_{positive} \ case)\ }\ \ \ &E_n&=\ \ \ & \left(  { \hat{c}_{n} }  \widetilde{T}_1+\hat{d}_{n} \right)E_{-n}&  \ 
\end{aligned}
\end{equation}

\ifJOLT \begin{Proposition} \else \begin{props} \fi
\label{SCALE}
For $n \ge 0:$
\begin{eqnarray*}
\hat{a}_{-(n+1)} = -\frac{1}{cdq^n} &
& \hat{c}_n = - \frac{1}{ab} \hspace{15mm} \\
\widetilde{\zeta'}_{0,-(n+1)} =  -cdq^{-1} (abcdq^{2n}| q)_1& 
& \widetilde{\zeta}_{1,n} = -abq^{-n}(abcdq^{2n-1}| q)_1  \\
\hat{b}_{-(n+1)} = \frac{(c+d)- cdq^{n }(a+b)}{cdq^{n } (abcdq^{2n}| q)_1}&  
& \hat{d}_n = - \frac{( abq^{n} | q )_1+ab(cdq^{n-1}| q)_1}{ab(1-abcdq^{2n-1})} 
\end{eqnarray*}
\ifJOLT \end{Proposition} \else \end{props} \fi

\begin{proof}
We use the combination of the creation operator point of view with Corollary \ref{T_U_BS}.
\begin{enumerate} 
\item For the $\hat{c}_n$ and $\widetilde{\zeta}_{1,n}$ determination, we start with
\ba
E_n= & \left(   \hat{c}_{n}  \widetilde{T}_1+\hat{d}_{n} \right)E_{-n} \\
     = &   \hat{c}_{n}   \widetilde{T}_1 h_n \text{ mod } \mathcal{R}_{-n}  \\
     =& -ab \hat{c}_n f_n \text{ mod } \mathcal{R}_{-n} \\
     \ea
     Since $E_n = f_n \text{ mod } \mathcal{R}_{-n},$ this gives the asserted formula for $\hat{c}_n.$ 
     This also means 
$$
 \widetilde{T}_1 E_{-n} = \lb \hat{c}_n \rb^{-1} E_n \text{ mod } E_{-n}  
$$
Then
\begin{multline*}
\widetilde{\zeta}_{1, n} E_n = \widetilde{\mathcal{S}}_1 E_{-n} = [\widetilde{T}_1,\widetilde{Y}] E_{-n}= \left(\widetilde{\mu}_{-n} -\widetilde{\mu}_n \right)\left(\widetilde{T}_1 E_{-n}\right) \text{ mod } E_{-n} \\ =\left(\widetilde{\mu}_{-n} -\widetilde{\mu}_n \right) \lb \hat{c}_n \rb^{-1} E_n 
\end{multline*}
gives the $\widetilde{\zeta}_{1,n}$  determination.
\item For the $\hat{a}_{-(n+1)} $ and $\widetilde{\zeta'}_{0,-(n+1)}$ determination, we start with
\ba
E_{-(n+1)}= & \left(   \hat{a}_{-(n+1)}  \widetilde{U}_0+\hat{b}_{-(n+1)} \right)E_{n} \\
     = &   \hat{a}_{-(n+1)}   \widetilde{U}_0 f_{n\text{\_q}} \text{ mod } E_{n}  \\
          = &   \hat{a}_{-(n+1)} \lp- \frac{cd}{q}\rp   h_{n+1\text{\_q}} \text{ mod } \mathcal{R}_{n}  \\
     =& -cdq^n \hat{a}_{-(n+1)} h_{n+1} \text{ mod } \mathcal{R}_{n}. \\
\ea
Since $E_{-(n+1)} = h_{n+1} \text{ mod } \mathcal{R}_{n},$ this gives the asserted formula for $\hat{a}_{-(n+1)}.$
This also means 
$$
 \widetilde{U}_0 E_{n} = \lb \hat{a}_{-(n+1)} \rb^{-1} E_{-(n+1)} \text{ mod } E_{n} . 
$$
Then
\begin{multline*}
\widetilde{\zeta'}_{0,-(n+1)} E_{-(n+1)} =\widetilde{\mathcal{S'}}_0 E_{n} = [\widetilde{Y},\widetilde{U}_0] E_{n}= \left(\widetilde{\mu}_{-(n+1)} -\widetilde{\mu}_n \right)\left(\widetilde{U}_0 E_{n}\right) \text{ mod } E_{n} \\=\left(\widetilde{\mu}_{-(n+1)} -\widetilde{\mu}_n \right) \lb \hat{a}_{-(n+1)} \rb^{-1} E_{-(n+1)} 
\end{multline*}
gives the $\widetilde{\zeta'}_{0,-(n+1)}$  determination.
\item The $\hat{b}_{-(n+1)}$ and $\hat{d}_n$ formulas come from equations (\ref{SRCN0}) and (\ref{SRCN1}) together with the previously translated formulas from \cite{ASKVOL} for $\widetilde{b}_{-(n+1)}$ and $\widetilde{d}_n.$
\end{enumerate} 
\end{proof}

In this paper, we only use the zig-zag increasing cases above and so have  included just the proofs of those.  However Appendix B  includes a table
with the zig--zag decreasing cases as well. (Formulas (4.11) and (4.12) of \cite{ZHED} could be appealed to since they are equivalent to the determination of the 
$\hat{c}_r,\hat{d_r}$ for any sign of $r.$) We mention the simplified forms of those others as well because they may be of interest. They may 
be established, e.g. in the $\hat{a}_n,\hat{b_n}$ case, either using Corollary 5.2 of \cite{ANNALS} about $\mathcal{S}_i^2$ or by noting that the relations
\begin{eqnarray*}
E_{-(n+1)} = \lp \hat{a}_{-(n+1)} \widetilde{U_0} + \hat{b}_{-(n+1)} \rp E_n  & \hspace{15mm}  &E_n = \lp \hat{a}_n \widetilde{U_0} + \hat{b}_n \rp E_{-(n+1)} 
\end{eqnarray*}
are inverse to each other. (Corollary \ref{T_U_BS} also clarifies what happens in low zig-zag co-degree.)

If using Corollary 5.2 of \cite{ANNALS}, one might first confirm, by thinking about the DAHA relation, that $t_i^{\frac{1}{2}} \text{ and } u_i^{\frac{1}{2}}$
are mapped by the involution $\epsilon$ of that paper to what one might guess from its action on $T_i$ and $U_i.$ That then makes immediate
what $\epsilon(a)$ and $\epsilon(c)$ are. Then the automorphism property, together with $t_0=-cdq^{-1} \text{ and } t_1=-ab$ gives the
needed $\epsilon(b)$ and $\epsilon(d)$

An immediate corollary of proposition \ref{SCALE} and recursion relations (\ref{SRCNX})
is the $\boldsymbol{\mathcal{T}}(a,e)$ part of the first step of (\ref{PLANA}), Proof Plan A :
\ifJOLT \begin{Corollary} \else \begin{cors} \fi
\label{PLANA1}
The entries of both $T(a,e)$ and  $\boldsymbol{\mathcal{T}}(a,e) $ are rational functions of $e$ with coefficients in the filed $\mathbb{Q}(a,b,c,d,q).$
\ifJOLT \end{Corollary} \else \end{cors} \fi
(Rationality of entries of $T(a,e)$ is immediate from the formulas written down in theorem \ref{ETC}, since $T(a,e)$ just refers to  the
matrix given by the formulas (\ref{ETCFRM}), perse.)


We now work out the details of the low zig-zag co-degree expansion of the nonsymmetric polynomials in terms of the almost symmetric bases.

The recursive relations (\ref{SRCNX})
with exact scaling factors implies the following for the low zig-zag co-degree almost symmetric basis coefficients:
\ifJOLT \begin{Proposition} \else \begin{props} \fi
\label{RCNL}
For $n \ge 0:$
\begin{equation*}
\begin{aligned}
\lambda_{n+1,n+1} &\ =&\lambda_{-(n+1),-(n+1)} \lb -ab\hat{c}_{n+1} \rb \\
\lambda_{-(n+1),-(n+1)\text{\_q}}&\  = &\lambda_{nn\text{\_q}} \lb \left( -\frac{cd}{q}\right) \hat{a}_{-(n+1)} \rb  \\
\mu_{-(n+1),n+1} &\ =&\lambda_{-(n+1),-(n+1)} \lb -\hat{c}_{n+1}+\hat{d}_{n+1} \rb  \\ 
\mu_{n,-(n+1)\text{\_q}} &\ =&\lambda_{nn\text{\_q}} \lb \left( c+d\right) \hat{a}_{-(n+1)}+ \hat{b}_{-(n+1)} \rb  \\
\end{aligned}
\end{equation*} 
\bmlst
\lambda_{n,n+1} =\lambda_{-(n+1),-(n+1)} \lb (a+b)\hat{c}_{n+1} \rb +\mu_{n,-(n+1)} \lb    (-ab)\hat{c}_{n+1}+ \hat{d}_{n+1}\rb  \\
+\lambda_{-n,-(n+1)} \lb -ab \hat{c}_{n+1}\rb \\
\lambda_{-n,-(n+1)\text{\_q}} =  \mu_{-n,n\text{\_q}}    \hat{b}_{-(n+1)}+\lambda_{n-1,n\text{\_q}} \lb \left( -\frac{cd}{q}\right) \hat{a}_{-(n+1)} \rb
\end{multline*}  
\ifJOLT \end{Proposition} \else \end{props} \fi 
\begin{proof}
For the first $2$ assertions on the left and the next to last one, recall that our standard notation is:
\begin{eqnarray*}
E_{n+1}& =& \lambda_{n+1,n+1} f_{n+1}+\mu_{-(n+1),n+1} h_{n+1}+ \lambda_{n,n+1} f_n \\
E_{-(n+1)} &= &\lambda_{-(n+1),-(n+1)} h_{n+1}+\mu_{n,-(n+1)} f_{n}+ \lambda_{-n,-(n+1)} h_n. 
\end{eqnarray*}
So
\begin{multline*}
\left( \hat{c}_{n+1} \widetilde{T}_1 + \hat{d}_{n+1} \right) E_{-(n+1)}  \text{ mod } \mathcal{R}_{-n}  = \\
\lambda_{-(n+1),-(n+1)} \lb  \hat{c}_{n+1} \lb -abf_{n+1} -h_{n+1} +(a+b)f_n \rb+  \hat{d}_{n+1}h_{n+1}\rb +\\ 
\mu_{n,-(n+1)} \lb  \hat{c}_{n+1} \lb -abf_n \rb+  \hat{d}_{n+1}f_n\rb +  \lambda_{-n,-(n+1)} \lb\hat{c}_{n+1} \lb -abf_n \rb \rb\\
=f_{n+1} \lb \lambda_{-(n+1),-(n+1)} \lb -ab \hat{c}_{n+1}  \rb \rb \\+ h_{n+1} \lb \lambda_{-(n+1),-(n+1)} \lb -\hat{c}_{n+1} +\hat{d}_{n+1} \rb \rb +\\ 
f_{n} \lb \lambda_{-(n+1),-(n+1)} \lb (a+b) \hat{c}_{n+1}   \rb  \right. \\ +  \left. \mu_{n,-(n+1)} \lb -ab\hat{c}_{n+1}+\hat{d}_{n+1} \rb + \lambda_{-n,-(n+1)} \lb -ab \hat{c}_{n+1} \rb \rb
\end{multline*}
from which we can read off the three results.

For the other $3$ assertions, start with our standard notation of:
\begin{eqnarray*}
E_{-(n+1)} &=& \lambda_{-(n+1),-(n+1)\text{\_q}} h_{n+1\text{\_q}}+\mu_{n,-(n+1)\text{\_q}} f_{n\text{\_q}}+ \lambda_{-n,-(n+1)\text{\_q}} h_{n\text{\_q}} \\
E_{n} &=& \lambda_{n,n\text{\_q}} f_{n\text{\_q}}+\mu_{-n,n\text{\_q}} h_{n\text{\_q}}+ \lambda_{n-1,n\text{\_q}} f_{n-1,\text{\_q}}. 
\end{eqnarray*}
So
\begin{multline*}
\left( \hat{a}_{-(n+1)} \widetilde{U}_0 + \hat{b}_{-(n+1)} \right) E_{n}  \text{ mod } \mathcal{R}_{n-1}  = \\
 \lambda_{n,n\text{\_q}} \lb  \hat{a}_{-(n+1)}\lb \left( -\frac{cd}{q}\right) h_{n+1\text{\_q}} +(c+d) f_{n\text{\_q}} \rb  + \hat{b}_{-(n+1)} f_{n\text{\_q}} \rb +\\ \mu_{-n,n\text{\_q}} \lb \hat{a}_{-(n+1)}\lb qf_{n-1\text{\_q}} \rb  + \hat{b}_{-(n+1)}h_{n\text{\_q}} \rb + \lambda_{n-1,n\text{\_q}} \lb \hat{a}_{-(n+1)}\lb\left( -\frac{cd}{q}\right) h_{n\text{\_q}}  \rb  \rb  \\
 =h_{n+1\text{\_q}} \lb \lambda_{n,n\text{\_q}} \lb  \left( -\frac{cd}{q}\right) \hat{a}_{-(n+1)} \rb \rb + f_{n\text{\_q}} \lb \lambda_{n,n\text{\_q}} \lb(c+d)\hat{a}_{-(n+1)}+ \hat{b}_{-(n+1)}  \rb  \rb\\ +h_{n\text{\_q}} \lb \mu_{-n,n\text{\_q}}   \lb \hat{b}_{-(n+1)} \rb +\lambda_{n-1,n\text{\_q}} \lb  \left( -\frac{cd}{q}\right) \hat{a}_{-(n+1)}\rb \rb  
\end{multline*}
from which we can read off the other three results.
\end{proof}
\ifJOLT \begin{Proposition} \else \begin{props} \fi
\label{ZZCD1}
For $n \ge 0:$ 
\begin{eqnarray}
\label{MUQN} \mu_{n,-(n+1)} \   = &\ds \frac{ abq^{n}(c+d) - (a+b)  }{ (abcdq^{2n}| q)_1  } =  \mu_{n,-(n+1)\_\text{q}} \\
\label{MUN} \mu_{-(n+1),n+1}\   = &\ds  - \frac{ (q^{n+1}, cdq^{n}| q)_1 } {  (abcdq^{2n+1}| q)_1  } \hspace{25mm}\\
 \label{MUQP} \mu_{-(n+1),n+1\_\text{q}} \ = & \ds   \frac{ cd(q^{n+1}, ab q^{n+1 }| q)_1  } {  q(abcdq^{2n+1}| q)_1 }.  \hspace{25mm} 
\end{eqnarray}
\ifJOLT \end{Proposition} \else \end{props} \fi

\begin{proof}

By proposition \ref{RCNL}
\begin{multline}
\mu_{n,-(n+1)\text{\_q}} =\lambda_{n,n\text{\_q}} \lb \left( c+d\right) \hat{a}_{-(n+1)}+ \hat{b}_{-(n+1)} \rb= (c+d) \lp -\frac{1}{cdq^n}  \rp +\\
\frac{(c+d)- cdq^{n }(a+b)}{cdq^{n } (abcdq^{2n}| q)_1} =\frac{ (c+d) abcdq^{2n}- cdq^{n }(a+b)}{cdq^{n } (abcdq^{2n}| q)_1} = \text{ as in }  (\ref{MUQN}) \text{ above.} 
\end{multline}
And by proposition \ref{TOFROMQ} this is also the value of $\mu_{n,-(n+1)} .$

Similarly 
\begin{multline}
\mu_{-(n+1),n+1} =\lambda_{-(n+1),-(n+1)} \lb -\hat{c}_{n+1}+\hat{d}_{n+1} \rb \\ = - \lp - \frac{1}{ab}  \rp+ \lp - \frac{( abq^{n+1} | q )_1+ab(cdq^{n}| q)_1}{ab(1-abcdq^{2n+1})} \rp \\
= \frac{- abcdq^{2n+1} +abq^{n+1} -ab + abcdq^{n}  }{ab(1-abcdq^{2n+1}) } =-\frac{(1-q^{n+1})(1-cdq^{n} )}{(1-abcdq^{2n+1})} 
\end{multline}
as in (\ref{MUN}) above. 

And by proposition \ref{TOFROMQ}
\begin{multline*}
 \mu_{-(n+1),n+1\text{\_q}} = \lambda_{n+1,n+1}(q^{-(n+1)}-1)+\mu_{-(n+1),n+1}q^{-(n+1)} \\
 = (q^{-(n+1)}-1) + q^{-(n+1)} \lp-\frac{(1-q^{n+1})(1-cdq^{n} )}{(1-abcdq^{2n+1})}  \rp \\
  =\lp \frac{ 1-q^{n+1} }{q^{n+1} (1-abcdq^{2n+1})} \rp\lp -abcdq^{2n+1}+cdq^{n}\rp = \text{ as in }  (\ref{MUQP}) \text{ above.} 
\end{multline*}
\end{proof}

\ifJOLT \begin{Theorem} \else \begin{thms} \fi
\label{CDG2COEFFS}
The zig-zag co-degree 2  coefficients  of the $E_r(a)$ are given, for $n \geq 1$  by:
\begin{eqnarray*}
\lambda_{-n, -(n+1)\_\text{q}}= -\frac{ (c+d)(q^n,abq^{n}| q)_1+ (a+b)(q^n,cdq^{n}| q)_1}{q^n(q,abcdq^{2n}| q)_{1} } \\ 
\lambda_{-n,-(n+1)} =  - \frac{(c+d)(q^n,abq^{n+1}| q)_1 + q(a+b) (q^n,cdq^{n-1}| q)_1 }{(q,abcdq^{2n}| q)_1} \\ 
 \lambda_{n-1,n} = \lambda_{n-1,n\_\text{q}}
 =  - \frac{ (c+d)(q^n,abq^{n}| q)_1 + q(a+b) (q^n,cdq^{n-1}| q)_1}{(q,abcdq^{2n-1}| q)_1}.\\ 
\end{eqnarray*}
\ifJOLT \end{Theorem} \else \end{thms} \fi 

\begin{proof}
We prove these $\lambda_{s,r} \text{ and } \lambda_{s,r\text{\_q}}$  formulas by induction on zig-zag degree $r.$

Note in the case $r=-1,$  $\lambda_{-0,-1}$ is not even defined (there is a constant term $\mu_{0,-1}$ in $E_{-1}$), so
we can treat it as $0$ by convention; we are not assuming anything from this proposition in proving the $r=1$ case and 
so view this case as the start of the zig-zag induction.
 
Assuming correctness for $r=-(n+1),$ we establish the $r=n+1 \ge 1$ case: (The first case to prove is $s=0,r=1.)$

\begin{multline}
\label{BOTLE2}
\lambda_{n,n+1} =\lambda_{-(n+1),-(n+1)} \lb (a+b)\hat{c}_{n+1} \rb +\mu_{n,-(n+1)} \lb  (-ab)\hat{c}_{n+1}+ \hat{d}_{n+1}  \rb\\+\lambda_{-n,-(n+1)} \lb -ab \hat{c}_{n+1}\rb \\
= \lp 1 \rp \lb (a+b) \lb - \frac{1}{ab} \rb  \rb + \lp\frac{ abq^{n}(c+d) - (a+b)  }{ (abcdq^{2n}| q)_1  }  \rp \cdot \\
\Big\{\Big[  (-ab) \lb - \frac{1}{ab}\rb \Big] - \frac{( abq^{n+1} | q )_1+ab(cdq^{n}| q)_1}{ab(1-abcdq^{2n+1})} \Big\} \\
+ \lp - \frac{(c+d)(q^n,abq^{n+1}| q)_1 + q(a+b) (q^n,cdq^{n-1}| q)_1 }{(q,abcdq^{2n}| q)_1} \rp \lb  (-ab) \lb - \frac{1}{ab}\rb \rb 
\end{multline}
The top line of the right hand side of the last equals sign combines to
\begin{multline}
\label{TOPL2}
\lp - \frac{ 1} {ab (abcdq^{2n},abcdq^{2n+1}| q)_1} \rp \Big\{ (a+b) (1-abcdq^{2n})(1-abcdq^{2n+1}) + \\
\lb abq^{n}(c+d) - (a+b)\rb \big\{ -ab(1-abcdq^{2n+1}) + (1-abq^{n+1}) +ab(1-cdq^n)\big\} \Big\}\\
=\lp - \frac{ 1} {ab (abcdq^{2n},abcdq^{2n+1}| q)_1} \rp \Big\{ (a+b) (1-abcdq^{2n})(1-abcdq^{2n+1}) + \\
\lb abq^{n}(c+d) - (a+b)\rb \big\{ (1-abq^{n+1}) (1-abcdq^{n})\big\} \Big\}
\end{multline}
\begin{multline}
\label{TOPL2B}
=\lp - \frac{ 1} {ab (abcdq^{2n},abcdq^{2n+1}| q)_1} \rp \Big\{ (a+b) \big\{  (1-abcdq^{2n})(1-abcdq^{2n+1}) - \\
(1-abq^{n+1}) (1-abcdq^{n})\big\} +(c+d) \big\{abq^n (1-abq^{n+1}) (1-abcdq^{n}) \big\}.
\end{multline}
For $n=0$ (the start of the induction), this is all there is for $\lambda_{0,1},$ and we note  this simplifies to the asserted formula.
Continuing for $n \ge1,$
combining equation (\ref{TOPL2B}) with the bottom line of (\ref{BOTLE2}) gives:
\begin{multline}
\label{CHUNK3}
\lp - \frac{ 1} {ab (q,abcdq^{2n},abcdq^{2n+1}| q)_1} \rp \Big\{ (a+b) \big\{ (1-q) (1-abcdq^{2n})(1-abcdq^{2n+1}) - \\
(1-q)(1-abq^{n+1}) (1-abcdq^{n})+abq(1-abcdq^{2n+1})(q^n,cdq^{n-1}| q)_1) \big\} \\
+(c+d)\big \{(1-q)abq^n (1-abq^{n+1}) (1-abcdq^{n}) +ab(1-abcdq^{2n+1})  (q^n,abq^{n+1}| q)_1)\big\}\Big\}.
\end{multline}
We view (\ref{CHUNK3}) as 
$$\lp - \frac{ 1} {ab (q,abcdq^{2n},abcdq^{2n+1}| q)_1} \rp \Big\{ (a+b)p_1 + (c+d)p_2 \Big\} $$
where $p_1$ and $p_2$ are polynomials. We will finish this step in the inductive proof of the $ \lambda_{n,n+1} $ formula in the theorem (with $n$ replaced by $n+1$) by showing
$$
p_1=abq(1-abcdq^{2n})\lb (q^{n+1},cdq^{n}| q)_1\rb   \hspace{5mm} p_2=ab(1-abcdq^{2n})\lb(q^{n+1},abq^{n+1}| q)_1 \rb.
$$
To verify the $p_2$ claim:
\begin{multline*}
p_2=(1-q)abq^n (1-abq^{n+1}) (1-abcdq^{n}) +ab(1-abcdq^{2n+1})  (q^n,abq^{n+1}| q)_1) \\
  = (abq^{n+1}| q)_1\big\{ (1-q)abq^n (1-abcdq^{n}) + ab(1-abcdq^{2n+1})  (1-q^n)\big\} \\
  \text{ (because the } abq^n \text{ and } a^2b^2cdq^{2n+1} \text{ terms cancel)} \\
  = (abq^{n+1}| q)_1\big\{ -abq^{n+1}-a^2b^2cdq^{2n} +ab + a^2b^2cdq^{3n+1}  \big\} \\
  =ab(abq^{n+1}| q)_1\big\{ (1-q^{n+1})(1-abcd q^{2n}) \big\} 
\end{multline*}
as asserted.
The $p_1$ simplification is a little more involved:
\begin{multline*}
p_1=(1-q) (1-abcdq^{2n})(1-abcdq^{2n+1}) -(1-q)(1-abq^{n+1}) (1-abcdq^{n})\\
+abq(1-abcdq^{2n+1})(q^n,cdq^{n-1}| q)_1) \\
=(1-q) (1-abcdq^{2n})(1-abcdq^{2n+1}) -(1-q)(1-abq^{n+1}) (1-abcdq^{2n}+abcdq^{2n}-abcdq^{n}) \\
+abq(1-abcdq^{2n}+abcdq^{2n}-abcdq^{2n+1})(q^n,cdq^{n-1}| q)_1) 
\end{multline*}
\begin{multline*}
=(1-abcdq^{2n})\big\{(1-q) (1-abcdq^{2n+1}) -(1-q) (1-abq^{n+1})+abq(1-q^n)(1-cdq^{n-1}) \big\} \\
-(1-q) (1-abq^{n+1}) \lb abcdq^n(q^n-1)\rb+ \lb a^2b^2cdq^{2n+1}(1-q) \rb (1-q^n) (1-cdq^{n-1})
\end{multline*}
\begin{multline*}
=(1-abcdq^{2n})\big\{(1-q) abq^{n+1}(1-cdq^n)+abq(1-q^n)(1-cdq^{n-1}) \big\} \\
+(1-q)  abcdq^n(1-q^n)\big\{(1-abq^{n+1}) +  abq^{n+1}(1-cdq^{n-1}) \big\}\\
=(1-abcdq^{2n})\big\{ (1-q) abq^{n+1}(1-cdq^n) \\+abq(1-q^n)(1-cdq^{n-1})+abcdq^n(1-q)(1-q^n) \big\} 
\end{multline*}
\begin{multline*}
=(1-abcdq^{2n})\big\{ (1-q) abq^{n+1}(1-cdq^n)+abq(1-q^n)\big[ (1-cdq^{n-1})+ cdq^{n-1}(1-q)\big] \big\} \\
=(1-abcdq^{2n})(1-cdq^n)\big\{(1-q) abq^{n+1}+abq(1-q^n)\big\} \\
=abq(1-abcdq^{2n})(1-cdq^n)\big\{1-q^{n+1}\big\} 
\end{multline*}
as claimed. This finishes showing that the the zig-zag degree $-(n+1)$ case implies the asserted $\lambda_{n,n+1}$ formula.

We have $\lambda_{n,n+1} =\lambda_{n,n+1\text{\_q}} $ by proposition \ref{TOFROMQ} finishing induction starting as well as going from $r=-(n+1)$ to  $r=n+1.$ 

Assuming correctness for $r=n \ge 1,$ we now establish the $r=-(n+1)$ case:

By proposition \ref{RCNL}
\begin{multline*}
\lambda_{-n,-(n+1)\text{\_q}} = \mu_{-n,n\text{\_q}}    \hat{b}_{-(n+1)}+\lambda_{n-1,n\text{\_q}} \lb \left( -\frac{cd}{q}\right) \hat{a}_{-(n+1)} \rb  \\
=   \lp  \frac{ cd(q^{n}, ab q^{n }| q)_1  } {  q(abcdq^{2n-1}| q)_1 } \rp \lb \frac{(c+d)- cdq^{n }(a+b)}{cdq^{n } (abcdq^{2n}| q)_1}  \rb \\ + \lp  -\frac{ (c+d)(q^n,abq^{n}| q)_1 + q(a+b) (q^n,cdq^{n-1}| q)_1}{(q,abcdq^{2n-1}| q)_1}\rp \lb \left( -\frac{cd}{q}\right) \lp -\frac{1}{cdq^n} \rp\rb \\
=  \lp  \frac{ (q^{n}| q)_1  } {  q^{n+1}(q,abcdq^{2n-1},abcdq^{2n}| q)_1 } \rp \Big\{ (1-q)(abq^n| q)_1\lb(c+d) -cdq^n(a+b) \rb  \\
-(1-abcdq^{2n}) \lb (c+d)(abq^{n}| q)_1 + q(a+b) (cdq^{n-1}| q)_1 \rb \Big\}
\end{multline*}
\begin{multline*}
=  \lp  \frac{ (q^{n}| q)_1  } {  q^{n+1}(q,abcdq^{2n-1},abcdq^{2n}| q)_1 } \rp \Big\{(c+d) (abq^n| q)_1\lb abcdq^{2n}-q  \rb \\
+(a+b)\lb -(1-q)(1-abq^n)(cdq^n)-q(1-abcdq^{2n})(1-cdq^{n-1})\rb \Big\} \\
=  \lp  \frac{ (q^{n}| q)_1  } {  q^{n+1}(q,abcdq^{2n-1},abcdq^{2n}| q)_1 } \rp \Big\{(c+d) (abq^n| q)_1\lb-q(abcdq^{2n-1}| q)_1  \rb \\
+(a+b)\lb abcdq^{2n} +cdq^{n+1} -q -abc^2d^2q^{3n} \rb \Big\} \\
=  \lp  \frac{ (q^{n}| q)_1  } {  q^{n+1}(q,abcdq^{2n-1},abcdq^{2n}| q)_1 } \rp \Big\{(c+d) (abq^n| q)_1\lb-q(abcdq^{2n-1}| q)_1  \rb \\
+(a+b)\lb abcdq^{2n}(cdq^{n}| q)_1-q(cdq^{n}| q)_1 \rb \Big\} \\
=  - \lp  \frac{ (q^{n}| q)_1  } {  q^{n}(q,abcdq^{2n}| q)_1 } \rp \Big\{ (c+d) (abq^n| q)_1 +(a+b)(cdq^{n}| q)_1 \Big\}  
\end{multline*}
as asserted for the $r=-(n+1)$ case. 

Now that we have $\lambda_{-n,-(n+1)\text{\_q}},$
\begin{multline*}
\lambda_{-n,-(n+1)}=\mu_{n,-(n+1)\text{\_q}}(q^n-1)+\lambda_{-n,-(n+1)\text{\_q}} q^{n} =(q^n-1)\lp \frac{ abq^{n}(c+d) - (a+b)  }{ (abcdq^{2n}| q)_1  }  \rp \\
+ q^n \lp -\frac{ (c+d)(q^n,abq^{n}| q)_1+ (a+b)(q^n,cdq^{n}| q)_1}{q^n(q,abcdq^{2n}| q)_{1} }  \rp \\
=-\lp \frac{1-q^n}{(q,abcdq^{2n}| q)_{1} } \rp\big\{ abq^{n}(c+d)(1-q) - (a+b)(1-q) +(c+d)(1-abq^n)  \\  +(a+b)(1-cdq^n) \big\} 
=-\lp \frac{1-q^n}{(q,abcdq^{2n}| q)_{1} } \rp\big\{  (c+d)(1-abq^{n+1})+(a+b)(q-cdq^n) \big\}. 
\end{multline*}
which agrees with the asserted $\lambda_{-n,-(n+1)}  .$
\end{proof}

The formulas (\ref{ETCFRM}) for $T(a,aq)$ involve factors $(ae^{-1}| q)_r.$ When $e=aq,$, these vanish for $r \ge 2.$ Consequently the decomposition
$$
T(a,aq)=\begin{bmatrix}
T ^{00}(a,aq) & T^{01}(a,aq) \\
T^{10}(a,aq) & T^{11}(a,aq)\\
\end{bmatrix}
$$
specializes to
$$
T^{00}(a,aq)=\begin{bmatrix}
\tau_{0,0}(a,aq)  & \tau_{0,1}(a,aq)  &   0 &0\\
0  &  \tau_{1,1}(a,aq)  &\tau_{1,2}(a,aq)     &\ddots \\
 0 &\ddots   & \ddots  &\ddots \\
\end{bmatrix}
$$
$$
\ \ T^{01}(a,aq)=\begin{bmatrix}
\sigma_{0, -1} (a,aq)  & 0 & 0   \\
0  & \sigma_{1,-2}(a,aq)  & 0  &   \\
0  &\ddots & \ddots  \\
\end{bmatrix}
$$
$$
T^{10}(a,aq)=
\begin{bmatrix}
0  & \sigma_{-1,1} (a,aq) & 0  &  0 \\
0  & 0  & \sigma_{-2,2} (a,aq)  &   \ddots \\
0  & 0 &0  &  \ddots \\
\end{bmatrix}
$$
$$
T^{11}(a,aq)=\begin{bmatrix}
\tau_{-1, -1} (a,aq)  & \tau_{-1,-2} (a,aq) & 0  &  0  \\
0  & \tau_{-2,-2}(a,aq)  & \tau_{-2,-3} (a,aq)  &  \ddots \\
0  & 0 & \tau_{-3,-3} (a,aq)  &  \ddots \\
\end{bmatrix}.
$$

This has strong implications for the form of the discrete co-cycle condition (\ref{COCYCAQ}) as well, which we come back to in section \ref{COCYCPRF}. 

\section{Combining the $d_{r,s}c_{|r|,|s|}$ Products into the $\tau \text{ and } \sigma$ \label{TAUSIGMA}}
The most combinatorially involved part of our proof of Theorem \ref{ETC}  involves co-cycle condition verification. 

To facilitate its formulation and the clarity of correctness of our arguments, it is
useful to combine the $d_{rs}c_{|r|,|s|}$ expressions into more distinctive and mnemonic expressions. As mentioned 
earlier, equations (\ref{TAUDEF}) and (\ref{SIGMADEF})  give the main notation we use, namely 
\begin{align*}
\tau_{r,s}&=&d_{r,s}c_{|r|,|s|} & \hspace{10mm} & \text{if } (r \ge 0 \text{ and } s \ge 0)  \text{ or } (r <  0 \text{ and } s < 0)\\
\sigma_{r,s}&=& d_{r,s}c_{|r|,|s|} & \hspace{10mm} & \text{if } (r \ge 0 \text{ and } s < 0)  \text{ or } (r <  0 \text{ and } s \ge 0). 
\end{align*}
In this way, the sign portion of the $r \preceq s$ relation leads to the natural $2 \times 2$ block matrix 
form of the transition matrices. (This is where the $4$ cases of the $d_{rs}$ come from.) Explicitly
this means we are using the notation:

For  $n \geq 0$ 
\begin{multline}
E_n(z; a,b,c,d|q) = \sum_{m=0}^n \lb \tau_{m,n}(a, e ; b, c, d | q) \rb E_m(z; e,b,c,d|q)\\  + \sum_{m=0}^{n-1} \lb \sigma_{-(m+1),n}(a, e ; b, c, d | q) \rb E_{-(m+1)}(z; e,b,c,d|q) \\
E_{-(n+1)}(z; a,b,c,d|q) = \sum_{m=0}^n \lb \tau_{-(m+1),-(n+1)}(a, e ; b, c, d | q) \rb E_{-(m+1)}(z; e,b,c,d|q) \\+ \sum_{m=0}^n \lb \sigma_{m,-(n+1)}(a, e ; b, c, d | q) \rb E_{m}(z; e,b,c,d|q) 
\end{multline}
where for $k,n \geq 0$
\begin{multline}
\label{ETCFRM}
\tau_{k,n}(a, e ; b, c, d | q) =  \frac{ (q^{n-k+1}| q)_k (eq)^{n-k}(bc q^k, bd q^k, cd q^k, ae^{-1} | q )_{n-k}   }  { (q| q)_k (abcdq^{n+k}, bcdeq^{2k} | q)_{n-k} } \\
\sigma_{k,-(n+1)}(a, e ; b, c, d | q) =\frac{ (q^{n-k+1}| q)_k e^{n+1-k}(bc q^k, bd  q^k, ae^{-1} | q )_{n+1-k}   (cd  q^k| q )_{n-k}   }  { (q| q)_k (abcdq^{n+k} , bcdeq^{2k} | q)_{n+1-k} } \\
\sigma_{-(k+1),n}(a, e ; b, c, d | q) =   \frac{(q^{n-k}| q)_{k+1} (bc q^{k + 1}, bd q^{k + 1} | q )_{n-k - 1}   (cd q^{k }, ae^{-1}| q )_{n-k}  }  { (q| q)_k (abcdq^{n+k} , bcdeq^{2k+1} | q)_{n-k} } \\
\times bcd e^{n-k} q^{n + k} \\
\tau_{-(k+1),-(n+1)}(a, e ; b, c, d | q) = \frac{(q^{n-k+1}| q)_k  e^{n-k}(bc q^{k+1}, bd q^{k+1}, cdq^k, ae^{-1}| q )_{n-k}  }  { (q| q)_k  (abcdq^{n+k+1}, bcdeq^{2k+1} | q)_{n-k} }.
\end{multline}

\ifEXTRAPROOFS

\bc \bf (Verification below is of $\tau,\sigma$ given previous definitions of $c_{rs},d_{rs}.$) \ec
\begin{proof} For $k,n \ge 0,$
\begin{enumerate}
\item 
\bmlst
\tau_{k,n}(a, e ; b, c, d | q) = d_{k,n} c_{k,n} = \left\{ \frac{q^{n-k} (abcdq^{n+k-1}| q)_1}{(abcdq^{2n-1}| q)_1}  \right\} \cdot \\
\left\{ \frac{e^{n-k} (q^{n-k+1} | q )_k(bc q^k, bd q^k, cd q^k, ae^{-1}  | q )_{n-k}   }  { (q  | q )_k (abcdq^{n+k-1}, bcdeq^{2k}  | q )_{n-k} }  \right\} \\
=  \frac{ (q^{n-k+1}| q)_k (eq)^{n-k}(bc q^k, bd q^k, cd q^k, ae^{-1} | q )_{n-k}   }  { (q| q)_k (abcdq^{n+k}, bcdeq^{2k} | q)_{n-k} } 
\end{multline*}
\item 
\bmlst
\sigma_{k,-(n+1)}(a, e ; b, c, d | q) = d_{k,-(n+1)} c_{k,n+1} = \left\{  \frac{ (q^{n-k+1 }|q)_{ 1 } } {( q^{n+1}, cdq^{n }|q)_{ 1 }}  \right\} \cdot \\
 \left\{\frac{e^{n-k+1} (q^{n-k+2} | q )_k(bc q^k, bd q^k, cd q^k, ae^{-1}  | q )_{n-k+1}   }  { (q  | q )_k (abcdq^{n+k}, bcdeq^{2k}  | q )_{n-k+1} }  \right\} \\
 = \frac{ (q^{n-k+1}| q)_k e^{n+1-k}(bc q^k, bd  q^k, ae^{-1} | q )_{n+1-k}   (cd  q^k| q )_{n-k}   }  { (q| q)_k (abcdq^{n+k} , bcdeq^{2k} | q)_{n+1-k} } 
\end{multline*}
 \item 
 \bmlst
\sigma_{-(k+1), n}(a, e ; b, c, d | q) = d_{-(k+1), n} c_{k+1,n} = \left\{ \frac{bcdeq^{n+k} (q^{k+1}, cdq^k, ae^{-1}q^{n-k-1 }|q)_{1 } } { (abcdq^{2n-1}, bcdeq^{2k+1 }|q)_{ 1}} \right\} \\
 \left\{ \frac{e^{n-k-1} (q^{n-k} | q )_{k+1}(bc q^{k+1}, bd q^{k+1}, cd q^{k+1}, ae^{-1}  | q )_{n-k-1}   }  { (q  | q )_{k+1} (abcdq^{n+k}, bcdeq^{2(k+1)}  | q )_{n-k-1} } \right\} \\
= \frac{(q^{n-k}| q)_{k+1} bcd e^{n-k} q^{n + k}(bc q^{k + 1}, bd q^{k + 1} | q )_{n-k - 1}   (cd q^{k }, ae^{-1}| q )_{n-k}  }  { (q| q)_k (abcdq^{n+k} , bcdeq^{2k+1} | q)_{n-k} } \\
\end{multline*}
 \item 
 \bmlst
\tau_{-(k+1),-(n+1)}(a, e ; b, c, d | q) = d_{-(k+1),-(n+1)} c_{k+1,n+1} = \\ \left\{  \frac{ (q^{k+1}, cdq^{k}, bcdeq^{n+k+1 }|q)_{1 } } { (q^{n+1}, cdq^n, bcdeq^{2k+1 }|q)_{ 1} }  \right\} \\
\left\{ \frac{e^{n-k} (q^{n-k+1} | q )_{k+1}(bc q^{k+1}, bd q^{k+1}, cd q^{k+1}, ae^{-1}  | q )_{n-k}   }  { (q  | q )_{k+1} (abcdq^{n+k+1}, bcdeq^{2(k+1)}  | q )_{n-k} }  \right\} \\
=  \frac{(q^{n-k+1}| q)_k  e^{n-k}(bc q^{k+1}, bd q^{k+1}, cdq^k, ae^{-1}| q )_{n-k}  }  { (q| q)_k  (abcdq^{n+k+1}, bcdeq^{2k+1} | q)_{n-k} } 
\end{multline*}
\end{enumerate}
\end{proof}
\fi 

\section{The True $\boldsymbol{\mathcal{T}}(a,aq)$ Matches the $\boldsymbol{T}(a,aq) $ of Theorem \ref{ETC} \label{TAQ}}

Here we use the notation 
$$
\Delta_{ae} f = f(e) -f(a)
$$
for the change in values of a function  $f$ as one moves from argument $a$ to argument $e.$

Since the $\{\lambda_{r,s},\mu_{r,s}\}$ are also entries in a change of basis relationship, the true  transition functions $\boldsymbol{\mathcal{T}}(a,aq) $ may be expressed
in terms of them.  However, in the low zig-zag co-degree cases, the transition matrix entries may also be conveniently determined by a 
successive substitution argument.

\ifJOLT \begin{Proposition} \else \begin{props} \fi
\label{LCDGTR}
\begin{eqnarray}
\boldsymbol{\mathcal{T}}^{00}_{n-1, n}(a,e) &=&\tau_{n-1,n}= -\Delta_{ae} \lambda_{n-1,n} +\lb \mu_{n-1,-n}(e) \rb \Delta_{ae} \mu_{-n,n} \label{DELTA00} \\
 \boldsymbol{\mathcal{T}}^{01}_{n, n}(a,e) &=&\sigma_{n,-(n+1)}=-\Delta_{ae}\mu_{n,-(n+1)}  \label{DELTA01}  \\
\boldsymbol{\mathcal{T}}^{10}_{n, n+1}(a,e) &=&\sigma_{-(n+1),n+1}= -\Delta_{ae}\mu_{-(n+1),n+1} \label{DELTA10} \\
\boldsymbol{\mathcal{T}}^{11}_{n-1,n}(a,e)&=&\tau_{-n,-(n+1)} =-\Delta_{ae} \lambda_{-n, -(n+1)} +\lb \mu_{-n,n}(e) \rb \Delta_{ae} \mu_{n,-(n+1)} 
\end{eqnarray}
\ifJOLT \end{Proposition} \else \end{props} \fi
\begin{proof}
We write down the proof of (\ref{DELTA00}) and  (\ref{DELTA10}), the other two being similar.

Recall $\lambda_{r,r}=1$ for all parameters and any sign of $r.$ Then  modulo $\mathcal{R}_{-(n-1)}:$
\begin{eqnarray*}
E_{n-1}(e) =\lb  \lambda_{n-1,n-1}(e)\rb f_{n-1} & \Rightarrow & f_{n-1}=E_{n-1}(e)\\
E_{-n}(e)= \lb \lambda_{-n,-n}(e) \rb h_n + \lb \mu_{n-1,-n}(e) f_{n-1}(e) \rb &  \Rightarrow & \\
 h_{n}=E_{-n}(e) - \lb \mu_{n-1,-n}(e) \rb E_{n-1}(e) \\
f_n = E_n(e) -  \lb \mu_{-n,n}(e) \rb h_n - \lb \lambda_{n-1,n}(e) \rb f_{n-1} & 
\end{eqnarray*}
$$\hphantom{f_n}    =E_n(e)   -  \lb \mu_{-n,n}(e) \rb \big( E_{-n}(e) - \lb \mu_{n-1,-n}(e) \rb E_{n-1}(e)  \big) - \lb  \lambda_{n-1,n}(e)\rb E_{n-1}(e) $$
So
\begin{multline*}
E_n(a) = f_n+ \lb \mu_{-n,n}(a) \rb h_n + \lb  \lambda_{n-1,n}(a)\rb f_{n-1} \\
= E_n(e)   -  \lb \mu_{-n,n}(e) \rb \big( E_{-n}(e) - \lb \mu_{n-1,-n}(e) \rb E_{n-1}(e)  \big) - \lb  \lambda_{n-1,n}(e)\rb E_{n-1}(e) \\
+ \lb \mu_{-n,n}(a) \rb \big( E_{-n}(e) - \lb \mu_{n-1,-n}(e) \rb E_{n-1}(e) \big) +  \lb  \lambda_{n-1,n}(a)\rb E_{n-1}(e) 
\end{multline*}
Combining terms and comparing coefficients of $E_{-n}(e)$ (for  (\ref{DELTA10}))  and $E_{n-1}(e)$ (for  (\ref{DELTA00}))with the definition
$$
E_n(a)= \tau_{n,n}E_n(e) + \sigma_{-n,n} E_{-n}(e) + \tau_{n-1,n}E_{n-1}(e) \ \ \  \text{ mod }\mathcal{R}_{-(n-1)}
$$
gives the asserted formulas.
\end{proof}

The zig-zag co-degree $1$ formulas (\ref{DELTA01}) and (\ref{DELTA10}) above have an obvious linearity based on 
$$
 \Delta_{uw} f= \Delta_{uv} f + \Delta_{vw} f \text{ since } \lp f(w) - f(u) \rp = \lp f(v) - f(u) \rp + \lp f(w) - f(v) \rp.
$$
This immediately implies that the zig-zag co-degree $1$ matrix entries of $T$ satisfy what is required by the discrete co-cycle condition. We shall need
these special cases in the next section, so we record them in the corollary below.
\ifJOLT \begin{Corollary} \else \begin{cors} \fi
\label{CDG1LIN}
For any $a,e,\text{ and } f,$
\begin{equation*}
\sigma_{n,-(n+1)}(a,e) = \sigma_{n,-(n+1)}(f,e) +  \sigma_{n,-(n+1)}(a,f).
\end{equation*}
\begin{equation*}
   \sigma_{-(n+1), n+1}(a,e) =  \sigma_{-(n+1), n+1}(f,e) +  \sigma_{-(n+1), n+1}(a,f). 
\end{equation*}
\ifJOLT \end{Corollary} \else \end{cors} \fi
\ifJOLT \begin{Proposition} \else \begin{props} \fi
\label{CDG1AEQ}

In the special case of $e=aq,$ the zig-zag co-degree $1$ transition functions satisfy:  
\begin{enumerate}
\item $\ds \boldsymbol{\mathcal{T}}^{01}_{n,n}(a,aq) =\sigma_{n, -(n+1)}(a,aq) = \frac{  aq (bc q^n, bd  q^n, q^{-1} | q )_{1}   }  {  (abcdq^{2n} | q)_{1} (abcdq^{2n+1} | q)_{1}   }$
\item $\ds \boldsymbol{\mathcal{T}}^{10}_{n, n+1}(a,aq) = \sigma_{-(n+1), n+1}(a,aq) =
 \frac{abcd    q^{2(n +1)}   (q^{n+1}, cd q^{n }, q^{-1}| q )_{1}  }  {  (abcdq^{2n+1} | q)_{1} (abcdq^{2(n+1)} | q)_{1}   } $
\end{enumerate}
\ifJOLT \end{Proposition} \else \end{props} \fi

\begin{proof}

\begin{enumerate}
\item 
By proposition \ref{LCDGTR}
$$ \boldsymbol{\mathcal{T}}^{01}_{n, n}(a,aq) =\sigma_{n, -(n+1)}(a,aq) =-\Delta_{ae}\mu_{n,-(n+1)} $$
So 
\begin{multline*}
\boldsymbol{\mathcal{T}}^{01}_{nn}(a,aq) = - \Big\{ \Big\{ - \frac{ \lb ( \lb aq \rb bq^{n} | q)_1 -1 \rb (c+d) + ( \lb aq \rb +b)  }{  ( \lb aq \rb bcdq^{2n} | q)_1 } \Big\} \\
- \Big\{  - \frac{ \lb (abq^{n} | q)_1 -1 \rb (c+d) + (a+b)  }{  (abcdq^{2n} | q)_1 }\Big\} \Big\} \\
=\Big\{ \frac{1}{(abcdq^{2n} | q)_1 (abcdq^{2n+1} | q)_1  } \Big\}  \cdot \\
\big\{ (abcdq^{2n} | q)_1\lb  - abq^{n+1} (c+d)  +(aq+b)\rb  \\- (abcdq^{2n+1} | q)_1\lb  - abq^{n} (c+d) +(a+b)\rb \big\}
\end{multline*}
Comparing with the statement of the proposition and keeping in mind that $$ q(q^{-1}| q)_1= -(q| q)_1,$$ 
we see that we need to verify equality of the numerators, both of which are polynomials; i.e. the vanishing of
\begin{multline*}
  (abcdq^{2n} | q)_1\lb  - abq^{n+1} (c+d) +(aq+b)\rb\\ - (abcdq^{2n+1} | q)_1\lb  - abq^{n}  (c+d) +(a+b)\rb \\
   +a(bc q^n| q )_{1}  (bd  q^n| q )_{1}    (q | q )_{1}  = 0
\end{multline*}
Introducing an additional variable $u$ (intuitively, replacing $q^n$), it
suffices to show that the {\em linear} in $q$  polynomial
\begin{multline*}
  (abcdu^{2} | q)_1\lb  - abuq (c+d) +(aq+b)\rb - (abcdu^2q | q)_1\lb  - abu  (c+d) +(a+b)\rb \\
   +a(bc u| q )_{1}  (bd  u| q )_{1}    (q | q )_{1}  = 0  
\end{multline*}
vanishes for all $a,b,c,d,u.$ 
\begin{enumerate}
\item For $q=0,$ we have 
$$ b(abcdu^{2} | q)_1- \lb  - abu  (c+d) +(a+b)\rb + a(bc u| q )_{1}  (bd  u| q )_{1}=0$$
\item For $q=1,$ we have 
\bmlst
(abcdu^{2} | q)_1\lb  - abu (c+d) +(a+b)\rb \\- (abcdu^2 | q)_1\lb  - abu  (c+d) +(a+b)\rb =0
\end{multline*}
\end{enumerate}
So, since the $1$ variable polynomial  in $q$ with coefficients in $\mathbb{Q}(a,b,c,d,u)$ is at most degree $1$  and vanishes at two distinct values of $q,$ it must be identically
zero and the required identity has been established.

\item By proposition \ref{LCDGTR}
 $$ \boldsymbol{\mathcal{T}}^{10}_{n, n+1}(a,aq) = \sigma_{-(n+1), n+1}(a,aq)= -\Delta_{ae}\mu_{-(n+1),n+1} $$
So
\bmlst
\boldsymbol{\mathcal{T}}^{10}_{n,n+1}(a,aq) = - \Big\{ \Big\{   - \frac{ (cdq^{n } | q)_1 ( q^{n+1}| q  )_1 } {  ( \lb aq \rb bcdq^{2n+1} | q)_1 } \Big\} 
 - \Big\{ - \frac{ (cdq^{n} | q)_1 ( q^{n+1}| q  )_1 } {  (abcdq^{2n+1} | q)_1 }  \Big\} \Big\}\\
= \frac{ (cdq^{n} | q)_1 ( q^{n+1}| q  )_1  }{ (abcdq^{2n+1} | q)_1 (abcdq^{2(n+1)} | q)_1 }\Big\{(abcdq^{2n+1} | q)_1 -(abcdq^{2(n+1)} | q)_1\Big\}\\
= \frac{ (cdq^{n} | q)_1 ( q^{n+1}| q  )_1 \big\{ - abcdq^{2n+1}(q| q)_1 \big\} }{ (abcdq^{2n+1} | q)_1 (abcdq^{2(n+1)} | q)_1 }\\
= \frac{ (cdq^{n} | q)_1 ( q^{n+1}| q  )_1 \big\{  abcdq^{2(n+1)}(q^{-1}| q)_1 \big\} }{ (abcdq^{2n+1} | q)_1 (abcdq^{2(n+1)} | q)_1 }
\end{multline*}

\end{enumerate}
\end{proof}

\ifJOLT \begin{Proposition} \else \begin{props} \fi
\label{CDG2AEQ}

In the special case of $e=aq,$ the zig-zag co-degree $2$ transition functions satisfy:
\ 

\begin{enumerate}
\item $\ds \boldsymbol{\mathcal{T}}^{00}_{n-1, n}(a,aq) =\tau_{n-1,n}= - \frac{   aq  (q^n,bc q^{n-1}, bd q^{n-1}, cd q^{n-1}| q )_{1}    }  {\lb (abcdq^{2n-1} | q)_{1}  \rb^2  }$
\item $\ds \boldsymbol{\mathcal{T}}^{11}_{n-1,n}(a,aq) =\tau_{-n,-(n+1)}=  -\frac{  a(q^n, bc q^{n}, bd q^{n},cdq^{n-1}| q )_{1}  }  {\lb  (abcdq^{2n} | q)_{1} \rb^2 }  $
\end{enumerate}

\ifJOLT \end{Proposition} \else \end{props} \fi

\begin{proof}
\begin{enumerate}
\item
By proposition \ref{LCDGTR}
$$\boldsymbol{\mathcal{T}}^{00}_{n-1, n}(a,aq) =-\Delta_{ae} \lambda_{n-1,n} +\lb \mu_{n-1,-n}(e) \rb \Delta_{ae} \mu_{-n,n} .$$
So 
\bmlst
\boldsymbol{\mathcal{T}}^{00}_{n-1,n}(a,aq)  = \\ - \Big\{\Big\{-{\ds \frac{(q^n| q)_1}{(q| q)_1  (\lb aq \rb bcdq^{2n-1}| q)_1}} \left[(c+d)(\lb aq \rb bq^{n}| q)_1 + q( \lb aq \rb +b) (cdq^{n-1}| q)_1 \right] \Big\}\\
   - \Big\{ -{\ds \frac{(q^n| q)_1}{(q| q)_1  (abcdq^{2n-1}| q)_1}} \left[(c+d)(abq^{n}| q)_1 + q(a+b) (cdq^{n-1}| q)_1 \right]\Big\}\Big\} \\
   + \Big\{  - \frac{ \lb (\lb aq \rb bq^{n-1} | q)_1 -1 \rb (c+d) + (  \lb aq \rb +b)  }{  ( \lb aq \rb bcdq^{2(n-1)} | q)_1 } \Big\} \cdot \Big\{ \\
   \Big\{- \frac{ (cdq^{n - 1} | q)_1 ( q^n| q  )_1 } {  ( \lb aq \rb bcdq^{2n-1} | q)_1 }  \Big\} 
-  \Big\{ - \frac{ (cdq^{n - 1} | q)_1 ( q^n| q  )_1 } {  (abcdq^{2n-1} | q)_1 }  \Big\}\Big\}  
\end{multline*}

Using the notation $\big[a,b,c,\dots\big]=(a| q)_1(b| q)_1(c| q)_1\ldots,$ 
the abbreviation variables $u=q^{n-1},y=abcd$
and multiplying by 
$$
 \frac{ (q,q)_1 \lb ( yq^{2n-1} | q)_1 \rb^2 ( yq^{2n} | q)_1 } {(q^n| q)_1}=\frac{ \lb q,yu^2q , yu^2q , yu^2q^{2}\rb } {(uq| q)_1}
$$ 
it suffices to show the vanishing of 
\bmlst 
p_1=\Big\{aq \big[q,  yu^2q^{2}, bcu, bdu, cdu \big]    \Big\} \\
+ \Big\{ \lb ( yu^2q| q)_1\rb^2  \left[(c+d)(a buq^{2}| q)_1 + q(  aq  +b) (cd u| q)_1 \right] \Big\}\\
- \Big\{   ( yu^2q | q)_1  ( yu^2q^{2} | q)_1 \left[(c+d)(abuq| q)_1 + q(a+b) (cdu| q)_1 \right] \Big\}  \\ 
+\Big\{   (q,q)_1 ( yu^2q | q)_1 \Big[ \lb (abuq | q)_1 -1 \rb (c+d) + (   aq  +b) \Big]  (cdu | q)_1\Big\} \\
 -  \Big\{   (q,q)_1( yu^2q^{2} | q)_1 \Big[ \lb ( abuq | q)_1 -1 \rb (c+d) + (   aq  +b) \Big] (cdu | q)_1 \Big\} \\
=\Big\{aq \big[q,  yu^2q^{2}, bcu, bdu, cdu \big]    \Big\} \\
+ ( yu^2q| q)_1 \Big\{(c+d) \big\{ ( yu^2q| q)_1 (abuq^2|q)_1 - ( yu^2q^2| q)_1 (abuq|q)_1 \big\}\\
                                 +(cdu|q)_1 \Big\{ aq \big\{ ( yu^2q | q)_1  (q) -  ( yu^2q^{2} | q)_1 (1)  \big\} + bq  \big\{ ( yu^2q | q)_1 -  ( yu^2q^{2} | q)_1 \big\} \Big\} \Big\} \\                                 
 +  \lb q,cdu \rb \big\{  \lp - abuq  \rp (c+d) + (   aq  +b)    \big\} \big\{ ( yu^2q | q)_1 -  ( yu^2q^{2} | q)_1 \big\}                                
\end{multline*}

Note (keeping in mind, e.g., $(yu^2q)(q) - (yu^2q^2)(1) =0$) that
\begin{eqnarray*}
( yu^2q | q)_1  (q) -  ( yu^2q^{2} | q)_1 (1)  = q-1  = -  (q| q)_1 \\
( yu^2q | q)_1 -  ( yu^2q^{2} | q)_1=- yu^2q (q| q)_1 
\end{eqnarray*}
\bml
( yu^2q | q)_1 (ab uq^2 | q)_1  - ( yu^2q^2 | q)_1  (ab u q | q)_1 \\
= -yu^2q -abuq^2 +yu^2q^2 +ab u q \\
=uq(q|q)_1 (ab -yu) \\
= uq(q|q)_1(ab -abcdu) \\ 
= ab u q (q|q)_1 (cdu | q)_1 
\end{multline}

So $p_1$ is also equal to
\bmlst
p_2=\Big\{aq \big[q,  yu^2q^{2}, bcu, bdu, cdu \big]    \Big\} \\
+ ( yu^2q| q)_1 \Big\{(c+d) \big\{ ab u q (q|q)_1 (cdu | q)_1  \big\}\\
                                 +(cdu|q)_1 \Big\{ aq \big\{  -  (q| q)_1 \big\} + bq  \big\{ - yu^2q (q| q)_1 \big\} \Big\} \Big\} \\
 +  \lb q,cdu \rb \big\{  \lp - abuq  \rp (c+d) + (   aq  +b)    \big\} \big\{ - yu^2q (q| q)_1 \big\} \\
\end{multline*}
Multiplying $p_2$ by $\lb q (q|q)_1\  (cdu | q)_1 \rb^{-1}$
we are reduced to showing the vanishing of
\bmlst
p_3=\Big\{a \big[  yu^2q^{2}, bcu, bdu \big]    \Big\} 
+ ( yu^2q| q)_1 \Big\{(c+d) \big\{ ab u     \big\}\\
                                 +\Big\{ a \big\{  -1   \big\} + b  \big\{ - yu^2q \big\} \Big\} \Big\} 
 +  \big\{  \lp - abuq  \rp (c+d) + (   aq  +b)    \big\} \big\{ - yu^2 (q| q)_1 \big\} \\
 = a \big[  abcdu^2q^{2}, bcu, bdu \big]    
+ a( abcdu^2q| q)_1 \Big\{b u(c+d)
                                 -\lp 1  +    b^2cdu^2q \rp  \Big\} \\
  - abcdu^2 (q| q)_1  \big\{  \lp - abuq  \rp (c+d) + (   aq  +b)    \big\}  \\
 = a \big[  abcdu^2q^{2}, bcu, bdu \big]    
+ a( abcdu^2q| q)_1 \Big\{-(1-bcu)(1-bdu)
                                +b^2cdu^2(1-q)   \Big\} \\
 - abcdu^2 (q| q)_1  \big\{ aq(1-bcu)(1-bdu)  +b(1-abcdu^2q)    \big\}  \\                                
 = a \big[ bcu, bdu \big]  \big\{( abcdu^2q^2| q)_1 - ( abcdu^2q| q)_1 - abcdu^2q (q| q)_1\big\}  \\
 + ab^2cdu^2 \lb q, abcdu^2q \rb \big\{1  - 1 \big\}\\
= a \big[ bcu, bdu \big]  \big\{(  abcdu^2q  (q| q)_1  - abcdu^2q (q| q)_1\big\} =0,
\end{multline*} 
 thus proving the first formula.
\item By proposition \ref{LCDGTR}
$$
\boldsymbol{\mathcal{T}}^{11}_{n-1,n}(a,aq) =\\ -\Delta_{ae} \lambda_{-n, -(n+1)} +\lb \mu_{-n,n}(e) \rb \Delta_{ae} \mu_{n,-(n+1)} .
$$
So  
\bmlst
\boldsymbol{\mathcal{T}}^{11}_{n-1,n}(a,aq) =\\- \Big\{ \Big\{   -{\ds \frac{(q^n| q)_1}{(q| q)_1 ( \lb aq \rb bcdq^{2n}| q)_1}} \left[(c+d)( \lb aq \rb bq^{n+1}| q)_1 + q( \lb aq \rb +b) (cdq^{n-1}| q)_1 \right] \Big\}  \\
 - \Big\{  -{\ds \frac{(q^n| q)_1}{(q| q)_1 (abcdq^{2n}| q)_1}} \left[(c+d)(abq^{n+1}| q)_1 + q(a+b) (cdq^{n-1}| q)_1 \right]   \Big\} \Big\}   
\end{multline*}
 \bmlst
\hphantom{\boldsymbol{\mathcal{T}}^{11}_{n-1,n}(a,aq) =} 
+\Big\{  - \frac{ (cdq^{n - 1} | q)_1 ( q^n| q  )_1 } {  ( \lb aq \rb bcdq^{2n-1} | q)_1 }\Big\} \cdot  \Big\{   \\
 \Big\{   - \frac{ \lb ( \lb aq \rb bq^{n} | q)_1 -1 \rb (c+d) + ( \lb aq \rb +b)  }{  ( \lb aq \rb bcdq^{2n} | q)_1 }\Big\}  \\
 -   \Big\{ - \frac{ \lb (abq^{n} | q)_1 -1 \rb (c+d) + (a+b)  }{  (abcdq^{2n} | q)_1 } \Big\}\Big\}  
\end{multline*}

Using the notation $\big[a,b,c,\dots\big]=(a| q)_1(b| q)_1(c| q)_1\ldots,$ 
the abbreviation variables $u=q^{n},y=abcd$
and multiplying by 
$$
 \frac{ (q,q)_1  \lb (yu^2| q)_1 \rb^2 (yu^2q| q)_1 } {(u| q)_1  }=\frac{\lb q,yu^2, yu^2, yu^2q \rb } { (u| q)_1},
$$

it suffices to show the vanishing of
\bmlst 
p_1=\Big\{ a\lb q, yu^2q ,bc u ,bd u \rb (cd q^{n-1}| q )_{1}    \Big\} \\
+ \Big\{ \lb (yu^{2},yu^2 \rb   \left\{ (c+d) (a buq^{2}| q)_1 + q(aq +b) (cd q^{n-1}| q )_{1} \right\} \Big\} \\
- \Big\{   \lb yu^2,yu^2q \rb \left\{(c+d)(abuq| q)_1 + q(a+b)(cd q^{n-1}| q )_{1}  \right\}\Big\}  \\
+ \Big\{  \lb q, yu^2\rb(cd q^{n-1}| q )_{1}\Big\{    -a buq  (c+d) +  aq +b)\Big\} \\
 -  \lb q, yu^2q \rb (cd q^{n-1}| q )_{1}\Big\{   - abu  (c+d) + (a+b)  \Big\}  \Big\} \\
 =\Big\{ a\lb q, yu^2q ,bc u ,bd u \rb (cd q^{n-1}| q )_{1}    \Big\} \\
+ (yu^2 | q)_1 \Big\{ (c+d) \big\{ (yu^2 | q)_1 (a buq^{2}| q)_1 - (yu^2 q | q)_1 (a buq | q)_1 \big\} \\
 + (cd q^{n-1} | q)_1  \Big\{ aq \big\{ (yu^2 | q)_1 (q) -(yu^2 q | q)_1 (1)\big\} \\ 
 +bq  \big\{ (yu^2 | q)_1  -(yu^2 q | q)_1 \big\}  \Big\} \Big\}\\
  +\lb q, cd q^{n-1} \rb \Big\{ a(1-bu(c+d)) \big\{ (yu^2 | q)_1 (q) -(yu^2 q | q)_1 (1)  \big\} \\ 
  +b  \big\{ (yu^2 | q)_1  -(yu^2 q | q)_1 \big\}   \Big\}
\end{multline*}

Note (keeping in mind, e.g., $(yu^2)(q) - (yu^2q)(1) =0$) that
\begin{eqnarray*}
( yu^2 | q)_1  (q) -  ( yu^2q | q)_1 (1)  = q-1  = -  (q| q)_1 \\
( yu^2 | q)_1 -  ( yu^2q | q)_1=- yu^2 (q| q)_1 
\end{eqnarray*}
\bml
( yu^2 | q)_1 (ab uq^2 | q)_1  - ( yu^2q | q)_1  (ab u q| q)_1 \\
= -yu^2 -abuq^2 +yu^2q +abuq \\
=u(q|q)_1 (abq -yu) 
=abu (q|q)_1 (q-cdu) \\
=abuq (q|q)_1 (1-cdq^{n-1}) 
=abuq (q|q)_1 (cdq^{n-1} | q)_1 
\end{multline}

So $p_1$ is also equal to
\bmlst
p_2=\Big\{ a\lb q, yu^2q ,bc u ,bd u \rb (cd q^{n-1}| q )_{1}    \Big\} \\
+ (yu^2 | q)_1 \Big\{ (c+d) \big\{abu q (q|q)_1 (cduq^{n-1}|q)_1 \big\} \\
 + (cd q^{n-1} | q)_1  \Big\{ aq \big\{-  (q| q)_1   \big\} 
 +bq  \big\{ - yu^2 (q| q)_1  \big\}  \Big\} \Big\}\\
  +\lb q, cd q^{n-1} \rb \Big\{ a(1-bu(c+d)) \big\{ -  (q| q)_1   \big\}  
 +b  \big\{ - yu^2 (q| q)_1  \big\}   \Big\} 
\end{multline*}

So, upon multiplying by $\lb (q| q)_1(cdq^{n-1}| q)_1 \rb^{-1}$ 
we are reduced to showing the vanishing of
\bmlst
p_3=\Big\{ a\lb yu^2q ,bc u ,bd u \rb    \Big\} 
+ (yu^2 | q)_1 \Big\{ (c+d) \big\{abu q \big\} \\
 + \Big\{ aq \big\{- 1   \big\}  
 +bq  \big\{ - yu^2 \big\}  \Big\} \Big\}
  +(q| q)_1\Big\{ a(1-bu(c+d)) \big\{ -  1   \big\}  
 +b  \big\{ - yu^2   \big\}   \Big\} \\
 =  a\lb abcdu^2q ,bc u ,bd u \rb     
+ q (abcdu^2 | q)_1 \Big\{ abu (c+d) 
  -a (1+b^2cdu^2)  \Big\} \\
  +a(q| q)_1\Big\{- (1-bu(c+d)) 
   - b^2cdu^2     \Big\} \\
=  a\lb abcdu^2q ,bc u ,bd u \rb     
- aq \lb abcdu^2, bcu, bdu \rb \\
  - a(q| q)_1 \lb bcu, bdu \rb \\  
=  a \lb bc u ,bd u \rb  \big\{   (abcdu^2q |q)_1 - q (abcdu^2 |q)_1 - (q|q)_1 \big\} \\
=  a \lb bc u ,bd u \rb  \big\{ (q | q )_1 -(q|q)_1 \big\} =0 .
\end{multline*}
Thus the proposition proof is complete.
\end{enumerate}

\end{proof}

Propositions \ref{CDG1AEQ} and \ref{CDG2AEQ} complete the proof of the second step of (\ref{PLANA}), Proof Plan A.

\section{Proof of the Discrete Co-cycle Identity for $\boldsymbol{T}$  \label{COCYCPRF}}

The block form of the equation 
$$T(aq^p,aq^{p+1}) T(a,aq^p) = T(a,aq^{p+1})$$

is
\bmlst
\begin{bmatrix}
T ^{00} (aq^p,aq^{p+1}& T^{01}(aq^p,aq^{p+1} \\
T^{10}(aq^p,aq^{p+1} & T^{11}(aq^p,aq^{p+1} \\
\end{bmatrix}
\begin{bmatrix}
T ^{00}(a,aq^p) & T^{01}(a,aq^p) \\
T^{10}(a,aq^p) & T^{11}(a,aq^p) \\
\end{bmatrix} \\
=
\begin{bmatrix}
T ^{00}(a,aq^{p+1}) & T^{01}(a,aq^{p+1}) \\
T^{10}(a,aq^{p+1}) & T^{11}(a,aq^{p+1}) \\
\end{bmatrix}
\end{multline*}
So we have the four matrix equations:
\begin{eqnarray*}
T ^{00}(a,aq^{p+1})  =  \Big[ T ^{00} (aq^p,aq^{p+1}) \Big] \Big[T^{00}(a,aq^p) \Big] &  +  & \Big[T ^{01} (aq^p,aq^{p+1}) \Big] \Big[T ^{10}(a,aq^p) \Big] \\
T ^{01}(a,aq^{p+1})  =  \Big[T ^{00} (aq^p,aq^{p+1}) \Big] \Big[T^{01}(a,aq^p) \Big] & +  & \Big[T ^{01} (aq^p,aq^{p+1})  \Big]\Big[T ^{11}(a,aq^p) \Big] \\
T ^{10}(a,aq^{p+1})  =  \Big[T ^{10} (aq^p,aq^{p+1}) \Big] \Big[T^{00}(a,aq^p) \Big] &+  & \Big[T ^{11} (aq^p,aq^{p+1}) \Big] \Big[T ^{10}(a,aq^p) \Big]\\
T ^{11}(a,aq^{p+1})  =  \Big[T ^{10} (aq^p,aq^{p+1}) \Big]\Big[T^{01}(a,aq^p) \Big] & +&  \Big[T ^{11} (aq^p,aq^{p+1} ) \Big] \Big[T ^{11}(a,aq^p) \Big]
\end{eqnarray*}
Keep in mind, as described 
at then end of section \ref{ALMSYM},
that the entries of each $T^{ij}(aq^p,aq^{p+1})$ are zero except (possibly) on the diagonal and the superdiagonal.

In terms of $\tau \text{ and } \sigma,$ these equations are:
\bml
\label{COCYC00DSP}
\tau_{k,n}(a,aq^{p+1}) = \big[ \tau_{k,k}(aq^p,aq^{p+1})\big] \big[ \tau_{k,n}(a,aq^p) \big]+ \big[ \tau_{k,k+1}(aq^p,aq^{p+1})\big] \cdot \big\{ \\ \big[ \tau_{k+1, n}(a,aq^p)\big] \big\}
   +\big[ \sigma_{k,-(k+1)}(aq^p,aq^{p+1})\big] \big[ \sigma_{-(k+1),n}(a,aq^p) \big] 
   \end{multline}  
   \bml
   \label{COCYC01DSP}
\sigma_{k,-(n+1)}(a,aq^{p+1}) = \big[ \tau_{k,k}(aq^p,aq^{p+1})\big] \big[ \sigma_{k,-(n+1)}(a,aq^p) \big]+ \big[ \tau_{k,k+1}(aq^p,aq^{p+1})\big] \cdot \big\{ \\\big[ \sigma_{k+1, -(n+1))}(a,aq^p)\big] \big\}
    +\big[ \sigma_{k,-(k+1)}(aq^p,aq^{p+1})\big] \big[ \tau_{-(k+1),-(n+1)}(a,aq^p) \big]
       \end{multline}  
   \bml
   \label{COCYC10DSP}
  \sigma_{-(k+1), n}(a,aq^{p+1}) = 
  \big[ \sigma_{-(k+1),k+1}(aq^p,aq^{p+1})\big] \big[ \tau_{k+1, n}(a,aq^p)\big] 
  \\  +\big[ \tau_{-(k+1),-(k+1)}(aq^p,aq^{p+1})\big]\cdot \big\{ 
     \big[ \sigma_{-(k+1),n}(a,aq^p) \big] \big\} \\
    + \big[ \tau_{-(k+1),-(k+2)}(aq^p,aq^{p+1})\big] \big[ \sigma_{-(k+2), n}(a,aq^p)\big]  
       \end{multline}  
   \bml
   \label{COCYC11DSP}
 \tau_{-(k+1),-(n+1)}(a,aq^{p+1}) = 
    \big[ \sigma_{-(k+1),k+1}(aq^p,aq^{p+1})\big] \big[ \sigma_{k+1, -(n+1)}(a,aq^p)\big] \\
  + \big[ \tau_{-(k+1),-(k+1)}(aq^p,aq^{p+1})\big] 
  \big[ \tau_{-(k+1),-(n+1)}(a,aq^p) \big] \\
  + \big[ \tau_{-(k+1),-(k+2)}(aq^p,aq^{p+1})\big] \big[ \tau_{-(k+2), -(n+1)}(a,aq^p)\big] 
  \end{multline}
In the proofs of these identites, we will often reduce them to the vanishing of a polynomial. To further that, the following easily proven identities will often be used:
  
  \ifJOLT \begin{Lemma} \else \begin{lemmas} \fi
\label{POCHNEG}

\ 

\begin{eqnarray*}
q^d(q^{e}| q)_1 &=& (q^{d+e}| q)_1-(q^d| q)_1 \\
q^d(q^{f-g}| q)_1 &= &(q^{d+f-g}| q)_1-(q^d| q)_1 \\
(q^{f-g}| q)_1 &=& q^{-g}\left\{ (q^{f}| q)_1-(q^g| q)_1\right\} \\
(q^{-p}| q)_1 &=& - q^{-p}(q^p| q)_1 \\
\end{eqnarray*}
\ifJOLT \end{Lemma} \else \end{lemmas} \fi


%
%
%

It is feasible to directly check the identities (\ref{COCYC00DSP}),(\ref{COCYC01DSP}),(\ref{COCYC10DSP}) and (\ref{COCYC11DSP})
since there is a great deal of cancellation. However the readability of the quantities involved is enhanced by formulating some simplification lemmas
for ratios which may be interpreted as appearing in the identities.

It may be helpful motivationally to note that in both Lemmas   \ref{URT} and  \ref{VRT}, the first index of the transition coefficient in the denominator is always the
zig-zag successor of that of the numerator.  And in Lemma \ref{URT}, there is a further difference in the powers of $q$ between numerator and denominator; $aq^{p+1}$ in the numerator vs. $aq^{p}$ in the denominator.

\ifJOLT \begin{Lemma} \else \begin{lemmas} \fi For $k,n \ge 0$
\label{URT}
\ 

\begin{enumerate}
\item $\ds \frac{ \tau_{k,n}(a,aq^{p+1}) } { \sigma_{-(k+1),n}(a,aq^{p}) }=\frac{q^{n-k}(bcq^{k}, bdq^{k} , q^{-(p+1)}|q)_1 } {bcdq^{2k} (q^{n-k}, q^{n-k-p-1}|q)_1  }. $

\item $\ds  \frac{ \sigma_{k,-(n+1)}(a,aq^{p+1}) } { \tau_{-(k+1),-(n+1)}(a,aq^{p}) } = \frac{aq^{n-k+p+1} (bcq^k, bdq^k, q^{-(p+1)}|q)_1} { 
(abcdq^{n+k}, abcdq^{n+k+p+1} |q)_1}.  $

\item $\ds  \frac{ \sigma_{-(k+1),n}(a,aq^{p+1}) } { \tau_{k+1,n}(a,aq^{p}) } = \frac{abcd q^{n+k+p+1} (q^{k+1}, cdq^{k}, q^{-(p+1)}|q)_1} { (abcdq^{n+k}, abcdq^{n+k+p+1} |q)_1}.$

\item $\ds \frac{ \tau_{-(k+1),-(n+1)}(a,aq^{p+1}) } { \sigma_{k+1,-(n+1)}(a,aq^{p}) } = \frac{q^{n-k} (q^{k+1}, cdq^k, q^{-(p+1)}|q)_1} {(q^{n-k},q^{n-k-p-1}|q)_1 }. $ 
\end{enumerate}
\ifJOLT \end{Lemma} \else \end{lemmas} \fi

\ifEXTRAPROOFS
\begin{proof}
\begin{enumerate}
\item 
\bmlst
\Big( \left\{ (q^{n-k+1}| q)_k [aq^{p+1}]^{n-k}q^{n-k}(bc q^k | q )_{n-k} \right\}  \cdot \\
\left. \left\{\frac{  (bd q^k| q )_{n-k}   (cd q^k| q )_{n-k}   {\ds  (q^{-(p+1)}  | q )_{n-k}  }  }  { (q| q)_k (abcdq^{n+k} | q)_{n-k}  (bcd[aq^{p+1}]q^{2k} | q)_{n-k} }\right\} \rp \cdot \\
\Big( \Big\{
\frac{1}
{ (q^{n-k}| q)_{k+1} bcd [aq^p]^{n-k} q^{n + k}(bc q^{k + 1}| q )_{n-k - 1}} \Big\}  \cdot \\
\Big\{ \left. \frac
{ (q| q)_k (abcdq^{n+k} | q)_{n-k}  (bcd[aq^p]q^{2k+1} | q)_{n-k} }
{  (bd q^{k + 1} | q )_{n-k - 1}   (cd q^{k }| q )_{n-k}  {\ds   (q^{-p} | q )_{n-k}  } }  
\Big\} \rp \\
= \frac{q^{n-k}(bcq^{k}, bdq^{k} , q^{-(p+1)}|q)_1 } {bcdq^{2k} (q^{n-k}, q^{n-k-p-1}|q)_1  }.  \\
\end{multline*}
\item 
\bmlst
\Big( 
\left\{(q^{n-k+1}| q)_k  \lb aq^{p+1} \rb^{n+1-k}(bc q^k| q )_{n+1-k}  \right\} \cdot \\
\Big\{ \frac{  (bd  q^k| q )_{n+1-k}   (cd  q^k| q )_{n-k}   {(q^{-(p+1)} | q )_{n +1- k} }  }  { (q| q)_k (abcdq^{n+k} | q)_{n+1-k}  (bcd \lb aq^{p+1} \rb q^{2k} | q)_{n+1-k} } \Big\}
\Big) \\
\Big(  
\left\{
\frac{1}{(q^{n-k+1}| q)_k  \lb aq^{p} \rb ^{n-k}} 
 \right\} \cdot \\
\Big\{ \frac
{ (q| q)_k (abcdq^{n+k+1} | q)_{n-k}  (bcd \lb aq^{p} \rb q^{2k+1} | q)_{n-k} }  
{ (bc q^{k+1} | q )_{n-k}  (bd q^{k+1}| q )_{n-k}   (cd q^k| q )_{n-k}   {\ds  (q^{-p} | q )_{n-k}  }  }  
\Big\}
\Big) \\
= \frac{aq^{n-k+p+1} (bcq^k, bdq^k, q^{-(p+1)}|q)_1} { 
(abcdq^{n+k}, abcdq^{n+k+p+1} |q)_1}. \\
\end{multline*}
\item 
\bmlst
\Big( 
\left\{ (q^{n-k}| q)_{k+1} bcd \lb aq^{p+1} \rb^{n-k} q^{n + k}(bc q^{k + 1}| q )_{n-k - 1}  \right\} \cdot \\
\Big\{ \frac
{ (bd q^{k + 1} | q )_{n-k - 1}   (cd q^{k }| q )_{n-k}  {\ds   (q^{-(p+1)} | q )_{n-k}  } }  { (q| q)_k (abcdq^{n+k} | q)_{n-k}  (bcd\lb aq^{p+1} \rb q^{2k+1} | q)_{n-k} }  
\Big\} \Big) \\
\Big( 
\left\{ \frac{1}{ (q^{n-k}| q)_{k+1}  \lb aq^{p} \rb ^{n-k-1}q^{n-k-1} } \right\} \cdot \\
\Big\{ \frac
{ (q| q)_{k+1} (abcdq^{n+k+1} | q)_{n-k-1}  (bcd \lb aq^{p} \rb q^{2(k+1)} | q)_{n-k-1} }
{(bc q^{k+1} | q )_{n-k-1}  (bd q^{k+1}| q )_{n-k-1}    (cd q^{k+1}| q )_{n-k-1}   {\ds  (q^{-p}| q )_{n-k-1}  }  }  \Big\} 
\Big) \\
= \frac{abcd q^{n+k+p+1} (q^{k+1}, cdq^{k}, q^{-(p+1)}|q)_1} { (abcdq^{n+k}, abcdq^{n+k+p+1} |q)_1}. \\
\end{multline*}
\item 
\bmlst
\Big(  
\left\{ (q^{n-k+1}| q)_k  \lb aq^{p+1} \rb^{n-k} \right\} \cdot \\
\Big\{ \frac{(bc q^{k+1}| q )_{n-k}  (bd q^{k+1}| q )_{n-k}   (cdq^k | q )_{n-k}  {\ds   (q^{-(p+1)} | q )_{n-k}  } }  { (q| q)_k  (abcdq^{n+k+1} | q)_{n-k}  (bcd\lb aq^{p+1} \rb q^{2k+1} | q)_{n-k} } 
\Big\} \Big) \\
\Big( 
\Big\{
\frac{1}{ (q^{n-k}| q)_{k+1}  \lb aq^{p} \rb ^{n-k}}
\Big\} \cdot \\
\Big\{ \frac
  { (q| q)_{k+1} (abcdq^{n+k+1} | q)_{n-k}  (bcd \lb aq^{p} \rb q^{2(k+1)} | q)_{n-k} }
{ (bc q^{k+1}| q )_{n-k}  (bd  q^{k+1}| q )_{n-k}   (cd  q^{k+1}| q )_{n-k-1}   
  {\ds  (q^{-p} | q )_{n - k} }  }  
  \Big\}
  \Big) \\
= \frac{q^{n-k} (q^{k+1}, cdq^k, q^{-(p+1)}|q)_1} {(q^{n-k},q^{n-k-p-1}|q)_1 }.\\
\end{multline*}
\end{enumerate}

\end{proof}
\fi 


\ifJOLT \begin{Lemma} \else \begin{lemmas} \fi  For $k,n \ge 0:$
\label{VRT}
\ 

\begin{enumerate}
\item $\ds \frac{ \tau_{k,n}(a,aq^{p}) } { \sigma_{-(k+1),n}(a,aq^{p}) } = \frac{(bcq^k,bdq^k, abcdq^{n+k+p} |q)_1} {bcd q^{2k} (q^{n-k}, abcdq^{2k+p} |q)_1}.$

\item $\ds \frac{ \sigma_{k,-(n+1)}(a,aq^{p}) } { \tau_{-(k+1),-(n+1)}(a,aq^{p}) } =\frac{aq^p(bcq^k,bdq^k, q^{n-k-p}|q)_1 } {(abcdq^{n+k},abcdq^{2k+p }|q)_1 }.
$

\item  $\ds \frac{ \sigma_{-(k+1),n}(a,aq^{p}) } { \tau_{k+1,n}(a,aq^{p}) } = \frac{abcdq^{2k+p+1}(q^{k+1},q^{n-k-p-1},cdq^k|q)_1 } {(abcdq^{n+k},abcdq^{2k+p+1 }|q)_1 }.$

\item $\ds \frac{ \tau_{-(k+1),-(n+1)}(a,aq^{p}) } { \sigma_{k+1,-(n+1)}(a,aq^{p}) } = \frac{ (q^{k+1}, cdq^k, abcdq^{n+k+p+1}|q)_1} {(q^{n-k}, abcdq^{2k+p+1 }|q)_1 }. $

\end{enumerate}
\ifJOLT \end{Lemma} \else \end{lemmas} \fi

\ifEXTRAPROOFS

\begin{proof}

\begin{enumerate}
\item 
\bmlst
\Big( 
\left\{ (q^{n-k+1}| q)_k [aq^p]^{n-k}q^{n-k} \right\} \cdot \\
\Big\{ \frac{ (bc q^k | q )_{n-k}  (bd q^k| q )_{n-k}   (cd q^k| q )_{n-k}   {\ds  (q^{-p}  | q )_{n-k}  }  }  { (q| q)_k (abcdq^{n+k} | q)_{n-k}  (bcd[aq^p]q^{2k} | q)_{n-k} } \Big\}
\Big) \\
\Big( 
\Big\{ \frac{1}{ (q^{n-k}| q)_{k+1} bcd [aq^p]^{n-k} q^{n + k}} \Big\} \cdot \\ 
\Big\{ \frac
{ (q| q)_k (abcdq^{n+k} | q)_{n-k}  (bcd[aq^p]q^{2k+1} | q)_{n-k} } 
{(bc q^{k + 1}| q )_{n-k - 1}  (bd q^{k + 1} | q )_{n-k - 1}   (cd q^{k }| q )_{n-k}  {\ds   (q^{-p} | q )_{n-k}  } }  \Big\}
\Big) \\
= \frac{(bcq^k,bdq^k, abcdq^{n+k+p} |q)_1} {bcd q^{2k} (q^{n-k}, abcdq^{2k+p} |q)_1}.\\
\end{multline*}
\item 
\bmlst
\Big( 
\left\{ (q^{n-k+1}| q)_k  \lb aq^{p} \rb ^{n+1-k} \right\} \cdot \\
\frac{ (bc q^k| q )_{n+1-k}  (bd  q^k| q )_{n+1-k}   (cd  q^k| q )_{n-k}   {(q^{-p} | q )_{n +1- k} }  }  { (q| q)_k (abcdq^{n+k} | q)_{n+1-k}  (bcd \lb aq^{p} \rb q^{2k} | q)_{n+1-k} }   \Big) \\
\Big( 
\Big\{ \frac{1}{ (q^{n-k+1}| q)_k  \lb aq^{p} \rb ^{n-k}(bc q^{k+1} | q )_{n-k}} \Big\} \cdot \\ 
\Big\{  \frac
{ (q| q)_k (abcdq^{n+k+1} | q)_{n-k}  (bcd \lb aq^{p} \rb q^{2k+1} | q)_{n-k} }  
{   (bd q^{k+1}| q )_{n-k}   (cd q^k| q )_{n-k}   {\ds  (q^{-p} | q )_{n-k}  }  }  
\Big\}
\Big) \\
= \frac{aq^p(bcq^k,bdq^k, q^{n-k-p}|q)_1 } {(abcdq^{n+k},abcdq^{2k+p }|q)_1 }. \\
\end{multline*}
\item 
\bmlst
\Big( 
\left\{ (q^{n-k}| q)_{k+1} bcd  \lb aq^{p} \rb ^{n-k} q^{n + k}\right\} \cdot \\
\frac{(bc q^{k + 1}| q )_{n-k - 1}  (bd q^{k + 1} | q )_{n-k - 1}   (cd q^{k }| q )_{n-k}  {\ds   (q^{-p} | q )_{n-k}  } }  { (q| q)_k (abcdq^{n+k} | q)_{n-k}  (bcd \lb aq^{p} \rb q^{2k+1} | q)_{n-k} }  \Big) \\
\Big(
\Big\{ \frac{1}{(q^{n-k}| q)_{k+1}  \lb aq^{p} \rb ^{n-k-1}q^{n-k-1} } \Big\} \cdot \\  
\Big\{  \frac
{ (q| q)_{k+1} (abcdq^{n+k+1} | q)_{n-k-1}  (bcd \lb aq^{p} \rb q^{2(k+1)} | q)_{n-k-1} }
{ (bc q^{k+1} | q )_{n-k-1}  (bd q^{k+1}| q )_{n-k-1}    (cd q^{k+1}| q )_{n-k-1}   {\ds  (q^{-p}| q )_{n-k-1}  }  }  
\Big\}
\Big) \\
= \frac{abcdq^{2k+p+1}(q^{k+1},q^{n-k-p-1},cdq^k|q)_1 } {(abcdq^{n+k},abcdq^{2k+p+1 }|q)_1 }. \\
\end{multline*}
\item 
\bmlst
\Big(  
\left\{ (q^{n-k+1}| q)_k   \lb aq^{p} \rb ^{n-k} \right\} \cdot \\
\Big\{ \frac{(bc q^{k+1}| q )_{n-k}  (bd q^{k+1}| q )_{n-k}   (cdq^k | q )_{n-k}  
  {\ds   (q^{-p} | q )_{n-k}  } 
  }  { (q| q)_k  (abcdq^{n+k+1} | q)_{n-k}  (bcd \lb aq^{p} \rb q^{2k+1} | q)_{n-k} } 
  \Big\} \Big) \\
\Big( 
\Big\{ \frac{1}{ (q^{n-k}| q)_{k+1}  \lb aq^{p} \rb ^{n-k} } \Big\} \cdot \\ 
\frac
  { (q| q)_{k+1} (abcdq^{n+k+1} | q)_{n-k}  (bcd \lb aq^{p} \rb q^{2(k+1)} | q)_{n-k} }  
{ (bc q^{k+1}| q )_{n-k}  (bd  q^{k+1}| q )_{n-k}   (cd  q^{k+1}| q )_{n-k-1}   
  {\ds  (q^{-p} | q )_{n - k} }   } 
  \Big) \\
= \frac{ (q^{k+1}, cdq^k, abcdq^{n+k+p+1}|q)_1} {(q^{n-k}, abcdq^{2k+p+1 }|q)_1 }.  \\
\end{multline*}
\end{enumerate}

\end{proof}
\fi 

\bc \bf The $T^{00}$ Identity\label{T00} \ec
\ifJOLT \begin{Proposition} \else \begin{props} \fi
\label{T00PROP}
When $0 \le k \le n-1$
\bmlst
\tau_{k,n}(a,aq^{p+1}) = \big[ \tau_{k,k}(aq^p,aq^{p+1})\big] \big[ \tau_{k,n}(a,aq^p) \big] + \big[ \tau_{k,k+1}(aq^p,aq^{p+1})\big] \big[ \tau_{k+1, n}(a,aq^p)\big] \\
  +\big[ \sigma_{k,-(k+1)}(aq^p,aq^{p+1})\big] \big[ \sigma_{-(k+1),n}(a,aq^p) \big].
\end{multline*}
 And 
$$  \tau_{n,n}(a,aq^{p+1}) = \big[ \tau_{n,n}(aq^p,aq^{p+1})\big] \big[ \tau_{n,n}(a,aq^p) \big].$$
\ifJOLT \end{Proposition} \else \end{props} \fi

\begin{proof}
The second identity just says $1=1 \cdot 1.$


For the first, upon dividing by $\sigma_{-(k+1),n}(a,aq^p),$ we see it is sufficient to prove
\bmlst
\frac{\tau_{k,n}(a,aq^{p+1})}{\sigma_{-(k+1),n}(a,aq^p)} = \left\{ \tau_{k,k}(aq^p,aq^{p+1})\right\} \left\{ \frac{\tau_{k,n}(a,aq^p)}{\sigma_{-(k+1),n}(a,aq^p)} \right\} \\ + \left\{ \tau_{k,k+1}(aq^p,aq^{p+1})\right\} \left\{ \frac{\tau_{k+1, n}(a,aq^p)}{\sigma_{-(k+1),n}(a,aq^p)} \right\} 
  +\left\{ \sigma_{k,-(k+1)}(aq^p,aq^{p+1})\right\} .
\end{multline*}
Using Lemmas \ref{URT} and \ref{VRT}, this means we need to show
\bml
\label{EQN1a}
\left\{ 
\frac{q^{n-k}(bcq^{k}, bdq^{k} , q^{-(p+1)}|q)_1 } {bcdq^{2k} (q^{n-k}, q^{n-k-p-1}|q)_1  }
\right\} 
  =   \left\{ 1 \right\} 
\frac{(bcq^k,bdq^k, abcdq^{n+k+p} |q)_1} {bcd q^{2k} (q^{n-k}, abcdq^{2k+p} |q)_1}
 \\
+ \left\{ \frac{ (q^{2}| q)_k [aq^{p+1}]q(bc q^k | q )_{1}  (bd q^k| q )_{1}   (cd q^k| q )_{1}   (q^{-1}| q )_{1}    }  { (q| q)_k ([aq^p]bcdq^{2k+1} | q)_{1}  (bcd[aq^{p+1}]q^{2k} | q)_{1} } \right\}   \cdot  \\ 
\left\{ \frac
{(abcdq^{n+k},abcdq^{2k+p+1 }|q)_1 } 
{abcdq^{2k+p+1}(q^{k+1},q^{n-k-p-1},cdq^k|q)_1 } 
\right\} \\
  +  \left\{\frac{ (q| q)_k [aq^{p+1}](bc q^k| q )_{1}  (bd  q^k| q )_{1}     {\ds  (q^{-1} | q )_{1} }  }  { (q| q)_k ([aq^p]bcdq^{2k} | q)_{1}  (bcd[ aq^{p+1}]q^{2k} | q)_{1} } \right\} \cdot 1 \\
\end{multline}

Multiplying (\ref{EQN1a}) by
$$ \frac{ bcdq^{3k+p+1} ( q, q^{n-k}, abcdq^{2k+p}, abcdq^{2k+p+1}, q^{n-k-p-1} | q )_1 } { (bcq^k, bdq^k | q)_1 }$$
we see it is sufficient to show the vanishing of the  polynomial $p_1$ below.


$\big($To arrive at the final form of $p_1,$ we use Lemma \ref{POCHNEG} above  to simplify the following expressions:
$$ 
(q^{-1}| q)_1,\ 
(q^{-(p+1)}| q)_1,\  
(q^{n-k}| q)_1, \ \text{and }
(q^{n-k-p-1}| q)_1. \ 
\big)$$

We will eventually reduce this identity to the vanishing of a 1-variable polynomial in $q$ with coefficients in the field $\mathbb{Q}(a,b,c,d,y,u,v,w)$ 
with the property that when 
$$
y=abcd \hspace{5mm} u= q^n \hspace{5mm} v= q^k\hspace{5mm} w= q^p
$$
we obtain a unit in the coefficient field times the difference between the two sides of equation (\ref{EQN1a}) above. So, effectively, we can use the variables $y,u,v,w$
as abbreviations for these expressions.

\begin{multline}
p_1=-\left\{ (yq^{2k+p} | q)_1 (yq^{2k+p+1} | q)_1  (q| q)_1   q^{n+p+1}     (q^{-(p+1)}| q )_{1}   \right\}   \\
 +   \left\{ 1 \right\}
  \left\{  (yq^{n+k+p} | q)_1  (yq^{2k+p+1} | q)_1 q^{k+p+1}    (q| q)_1   (q^{n-k-p-1} | q )_{1}    \right\} \\
    + \left\{  (yq^{n+k} | q)_1  (yq^{2k+p} | q)_1(q^{n-k}| q)_1    q^{k+p+2}       (q^{-1} | q )_{1}    
       \right\}  \\
    +  \left\{ (q^{n-k}| q)_1y (q| q)_1 q^{3k+2p+2}   (q^{-1} | q )_{1}  
     (q^{n-k-p-1} | q )_{1}    \right\} 
\end{multline}
\begin{multline*}
\hphantom{p_1}=-\left\{ (yq^{2k+p} | q)_1 (yq^{2k+p+1} | q)_1  (q| q)_1  q^{n+p+1}   \lb - q^{-(p+1)}(q^{p+1}| q)_1\rb   \right\} \\
 +   \left\{ 1 \right\}
  \left\{  (yq^{n+k+p} | q)_1  (yq^{2k+p+1} | q)_1  (q| q)_1  q^{k+p+1}   \lb q^{-(k+p+1)} \lb (q^{n}| q)_1-(q^{k+p+1}| q)_1 \rb \rb    \right\} \\
  + \left\{  (yq^{n+k} | q)_1  (yq^{2k+p} | q)_1  q^{k+p+2} \lb q^{-k} \lb (q^{n}| q)_1-(q^k| q)_1 \rb\rb       \lb - q^{-1}(q| q)_1 \rb   
       \right\} \\
 +  \left\{  \lb q^{-k} \lb (q^{n}| q)_1-(q^k| q)_1 \rb \rb (q| q)_1  yq^{3k+2p+2}  \lb- q^{-1}(q| q)_1 \rb \right. \cdot \\ 
\left. \lb q^{-(k+p+1)} \lb (q^{n}| q)_1-(q^{k+p+1}| q)_1 \rb \rb       \right\} 
\end{multline*}
\begin{multline*}
\hphantom{p_1}=\left\{ (yq^{2k+p} | q)_1 (yq^{2k+p+1} | q)_1  (q| q)_1   q^{n}    (q^{p+1}| q)_1   \right\}   \\
+   \left\{ 1 \right\}%
  \left\{  (yq^{n+k+p} | q)_1  (yq^{2k+p+1} | q)_1  (q| q)_1   \lb   (q^{n}| q)_1-(q^{k+p+1}| q)_1 \rb \    \right\} \\
    - \left\{  (yq^{n+k} | q)_1  (yq^{2k+p} | q)_1  q^{p+1}   \lb (q^{n}| q)_1-(q^k| q)_1 \rb       (q| q)_1    
       \right\} \hspace{40mm} \\
  -  \left\{  \lb (q^{n}| q)_1-(q^k| q)_1  \rb (q| q)_1  yq^{k+p}  \lb (q| q)_1 \rb   
\lb (q^{n}| q)_1-(q^{k+p+1}| q)_1 \rb       \right\} 
\end{multline*}

Using the notation $\big[a,b,c,\dots\big]=(a| q)_1(b| q)_1(c| q)_1\ldots,$ 
and multiplying by $\lb  (q| q)_1 \rb^{-1}$ 
it suffices to show the vanishing of 
\bmlst
p_2 = u[yv^2w, yv^2wq, wq] 
        + (vwq - u)[yuvw, yv^2wq] \\
       - wq(v-u)[yuv, yv^2w]
       -yvw(v-u)(vwq-u) (q| q)_1
\end{multline*}

This expression may be interpreted as a one variable polynomial of degree at most $2$ in $q$ with coefficients in the field $\mathbb{Q}(y,u,v,w).$

The coefficient of $q^2$ is
$$
u (1-yv^2w) (-yv^2w) (-w) + (vw) (1-yuvw) (-yv^2w) -yvw (v-u) (vw) (-1) =0.
$$
So the polynomial $p_2$ is of degree at most $1$ in $q.$
Evaluating 
\begin{description}
\item[at $q=0$ ] $p_2(0)=u(1-yv^2w)-u(1-yuvw) -yvw(v-u)(-u)=0$
\item[at $ q=u(vw)^{-1}$ ] Note at this point $yv^2wq=yuv$ and $wq=uv^{-1}.$ So
$$ p_2\big(u(vw)^{-1}\big) = u[yv^2w, yuv](1-uv^{-1}) - uv^{-1}(v-u)[yv^2w,yuv]=0. $$
\end{description}
Thus the polynomial $p_2$ is $0$ and the proposition is proven.

\end{proof}

\bc \bf The $T^{01}$ Identity \ec

\ifJOLT \begin{Proposition} \else \begin{props} \fi
\label{T01PROP}
When $0 \le k \le n-1$
\bml
\sigma_{k,-(n+1)}(a,aq^{p+1}) = \big[ \tau_{k,k}(aq^p,aq^{p+1})\big] \big[ \sigma_{k,-(n+1)}(a,aq^p) \big] \\
+ \big[ \tau_{k,k+1}(aq^p,aq^{p+1})\big] \big[ \sigma_{k+1, -(n+1))}(a,aq^p)\big] \\
  +\big[ \sigma_{k,-(k+1)}(aq^p,aq^{p+1})\big] \big[ \tau_{-(k+1),-(n+1)}(a,aq^p) \big].
\end{multline}
And
\bml
\sigma_{n,-(n+1)}(a,aq^{p+1}) = \big[ \tau_{n,n}(aq^p,aq^{p+1})\big] \big[ \sigma_{n,-(n+1)}(a,aq^p) \big] \\
  +\big[ \sigma_{n,-(n+1)}(aq^p,aq^{p+1})\big] \big[ \tau_{-(n+1),-(n+1)}(a,aq^p) \big].
\end{multline}

\ifJOLT \end{Proposition} \else \end{props} \fi

\begin{proof} The second is immediate from corollary \ref{CDG1LIN} together with the observation that $\tau_{r,r}=1$ for any sign of $r.$

For the first, upon dividing by $\tau_{-(k+1),-(n+1)}(a,aq^p),$ we see it is sufficient to prove
\bmlst
\frac{\sigma_{k,-(n+1)}(a,aq^{p+1}) } {\tau_{-(k+1),-(n+1)}(a,aq^p)} 
= \left\{ \tau_{kk}(aq^p,aq^{p+1})\right\} \left\{ \frac{\sigma_{k,-(n+1)}(a,aq^p)} {\tau_{-(k+1),-(n+1)}(a,aq^p)} \right\} \\
+ \left\{ \tau_{k,k+1}(aq^p,aq^{p+1})\right\} \left\{ \frac{\sigma_{k+1, -(n+1))}(a,aq^p) } {\tau_{-(k+1),-(n+1)}(a,aq^p)} \right\}
  +\left\{ \sigma_{k,-(k+1)}(aq^p,aq^{p+1})\right\} \cdot 1.
\end{multline*}
That means we need to show
\bml
\label{EQN1b} 
\frac{aq^{n-k+p+1} (bcq^k, bdq^k, q^{-(p+1)}|q)_1} { (abcdq^{n+k}, abcdq^{n+k+p+1} |q)_1}
  =   \left\{1  \right\} 
 \left\{  \frac{aq^p( bcq^k,bdq^k, q^{n-k-p} |q)_1 } {(abcdq^{n+k},abcdq^{2k+p }|q)_1 } \right\} \\
+ \left\{\frac{ (q^{2}| q)_k \lb aq^{p+1} \rb q(bc q^k | q )_{1}  (bd q^k| q )_{1}   (cd q^k| q )_{1}   {\ds  (q^{-1} | q )_{1}  }  }  { (q| q)_k ( \lb aq^{p} \rb bcdq^{2k+1} | q)_{1}  (bcd\lb aq^{p+1} \rb q^{2k} | q)_{1} }  \right\}  \cdot \\
 \left\{ \frac
  {(q^{n-k}, abcdq^{2k+p+1 }|q)_1 }
 { (q^{k+1}, cdq^k, abcdq^{n+k+p+1}|q)_1} 
 \right\} \\
  +  \left\{ \frac{ (q| q)_k  \lb aq^{p+1} \rb  (bc q^k| q )_{1}  (bd  q^k| q )_{1}     {\ds  (q^{-1}| q )_{1} }  }  { (q| q)_k ( \lb aq^{p} \rb bcdq^{2k} | q)_{1}  (bcd \lb aq^{p+1} \rb  q^{2k} | q)_{1} } \right\}  \cdot   
  1
\end{multline}


We will eventually reduce this identity to the vanishing of a 1-variable polynomial in $y$ with coefficients in the field $\mathbb{Q}(a,b,c,d,u,v,w,q)$ 
with the property that when 
$$
y=abcd \hspace{5mm} u= q^n \hspace{5mm} v= q^k\hspace{5mm} w= q^p
$$
we obtain a unit in the coefficient field times the difference between the two sides of equation (\ref{EQN1b})  above. So, effectively, we can use the variables $y,u,v,w$
as abbreviations for the above expressions.

Multiplying (\ref{EQN1b})  by  
\begin{equation}
\frac{q^k ( abcdq^{n+k}, abcdq^{n+k+p+1}, abcdq^{2k+p}, abcdq^{2k+p+1} |q)_1 } { a ( bcq^k,bdq^k | q)_1 }
\end{equation}
we see it is sufficient to show the vanishing of the  polynomial $p_1$ below. We use $y$ as an abbreviation for $abcd.$

$\big($To arrive at the final form of $p_1,$ we use Lemma \ref{POCHNEG} above  to simplify the following expressions:
$$
(q^{-1}| q)_1,\  
(q^{-p}| q)_1,\ 
(q^{-(p+1)}| q)_1,\ 
(q^{n-k}| q)_1, \ \text{and }
(q^{n-k-p}| q)_1. \big)
$$
\bmlst
p_1 =  -\left\{  (yq^{2k +p} | q)_{1}   (yq^{2k +p+1} | q)_{1}  
   q^{n+p+1}  (q^{-(p+1)}| q)_1  \right\}\\
  +   \left\{1  \right\}   
  \left\{    (yq^{n+k +p+1} | q)_{1}  (yq^{2k +p+1} | q)_{1}  q^{k+p}  
    (q^{n-k-p}| q)_1      \right\} \\
 + \left\{\frac{ (q^{n-k}| q)_1(yq^{n+k} | q)_{1}   (y q^{2k+p} | q)_{1}   q^{k+p+2} 
   (q^{-1}| q)_1 
  }  { (q| q)_1 
   }  \right\} \\ 
 +  \left\{ (y q^{n+k} | q)_{1} (y q^{n+k+p+1} | q)_{1}    q^{k+p+1}      (q^{-1}| q)_1   \right\} 
   \end{multline*}      
\bmlst
\hphantom{p_1 }=  -\left\{  (yq^{2k +p} | q)_{1}   (yq^{2k +p+1} | q)_{1}  q^{n+p+1} \Big\{   - q^{-(p+1)}(q^{p+1}| q)_1\Big\}   \right\}\\
  +   \left\{1  \right\}   
  \left\{    (yq^{n+k +p+1} | q)_{1}  (yq^{2k +p+1} | q)_{1}   q^{k+p} 
 \big\{ q^{-(k+p)} \lb (q^n| q)_1 -(q^{k+p}| q)_1 \rb  \big\}     \right\}   \\
 +\frac{ \left\{ \big\{q^{-k} \lb(q^n| q)_1 -(q^k| q)_1 \rb  \big\} (yq^{n+k} | q)_{1}   (y q^{2k+p} | q)_{1}  q^{k+p+2} 
  \lb -q^{-1}(q| q)_1 \rb 
  \right\}}  { (q| q)_1 } \\ 
 +  \left\{ (y q^{n+k} | q)_{1} (y q^{n+k+p} | q)_{1}    q^{k+p+1}     \lb -q^{-1}(q| q)_1 \rb  \right\} 
  \end{multline*}
 \bmlst
\hphantom{p_1 }=  \left\{  (yq^{2k +p} | q)_{1}   (yq^{2k +p+1} | q)_{1}   q^{n} \Big\{  (q^{p+1}| q)_1\Big\}   \right\}\\
  +   \left\{1  \right\}   
  \left\{    (yq^{n+k +p+1} | q)_{1}  (yq^{2k +p+1} | q)_{1}   
   \big\{ \lb (q^{n}| q)_1-(q^{k+p}| q)_1\rb    \big\}
      \right\} \\
 - \left\{(yq^{n+k} | q)_{1}   (y q^{2k+p} | q)_{1}   q^{p+1}
    \right\} 
 \big\{\lb (q^{n}| q)_1-(q^{k}| q)_1 \rb \big\} \\
 -  \left\{ (y q^{n+k} | q)_{1} (y q^{n+k+p+1} | q)_{1}    q^{k+p}     \lb (q| q)_1 \rb  \right\} 
  \end{multline*} 
 
  Using the notation $\big[a,b,c,\dots\big]=(a| q)_1(b| q)_1(c| q)_1\ldots,$ 
setting $y=abcd, u=q^n,v=q^k, \text{ and } w=q^p,$ we see it is sufficient to show the vanishing of

\bmlst
p_2=u\lb yv^2w, yv^2wq, wq\rb 
+(vw-u)\lb yuvwq, yv^2wq \rb\\
-wq(v-u) \lb yuv, yv^2w\rb
-vw \lb yuv, yuvwq,q\rb 
\end{multline*}
 
 This expression may be interpreted as a one variable polynomial of degree at most $2$ in $y$ with coefficients in the field $\mathbb{Q}(u,v,w,q).$
 Evaluating
 \begin{description}
\item [at $y=0$] $p_2(0)=u(1-wq)+(vw-u) -wq(v-u)-vw(1-q) =0.$
\item [at $y=(v^2wq)^{-1}$] When $y=(v^2wq)^{-1},$ note $yv^2w =q^{-1},$ $yuv= u(vwq)^{-1},$ and $yuvwq=uv^{-1}.$ So
\bmlst
p_2\big( (v^2wq)^{-1} \big) = -wq(v-u) \lb u(vwq)^{-1}, q^{-1} \rb -vw \lb u(vwq)^{-1}, uv^{-1}, q \rb\\
= (1-u(vwq)^{-1})\big\{  -wq(v-u)  \lp -q^{-1} (1-q) \rp -vw(1-uv^{-1}) (1-q) \big\} \\
= w(1-u(vwq)^{-1})(1-q)\big\{ v-u -v(1-uv^{-1}) \big\} =0
\end{multline*}
\item [at $y=(uv)^{-1}$] When $y=(uv)^{-1},$ note $yv^2w= u^{-1}vw,$ $yv^2wq=u^{-1}vwq ,$ and $yuvwq=wq .$ So
\bmlst
p_2\big( (uv)^{-1} \big) = (1-u^{-1}vwq)(1-wq) \big\{u(1 -  u^{-1}vw) +vw-u\big\}=0 
\end{multline*}
\end{description}
Thus the polynomial $p_2$ is $0$ and the proposition is proven.

%
\end{proof}

\bc \bf The $T^{10}$ Identity \ec

\ifJOLT \begin{Proposition} \else \begin{props} \fi
\label{T10PROP}
When $0 \le k \le n-2$
\bml
    \sigma_{-(k+1), n}(a,aq^{p+1}) =  \big[ \sigma_{-(k+1),k+1}(aq^p,aq^{p+1})\big] \big[ \tau_{k+1, n}(a,aq^p)\big] \\
+\big[ \tau_{-(k+1),-(k+1)}(aq^p,aq^{p+1})\big] \big[ \sigma_{-(k+1),n}(a,aq^p) \big] \\
  + \big[ \tau_{-(k+1),-(k+2)}(aq^p,aq^{p+1})\big] \big[ \sigma_{-(k+2), n}(a,aq^p)\big] 
\end{multline}
And
\bmlst
   \sigma_{-n, n}(a,aq^{p+1}) =  \big[ \sigma_{-n,n}(aq^p,aq^{p+1})\big] \big[ \tau_{n,n}(a,aq^p)\big] \\
+\big[ \tau_{-n,-n}(aq^p,aq^{p+1})\big] \big[ \sigma_{-n,n}(a,aq^p) \big]. 
\end{multline*}
\ifJOLT \end{Proposition} \else \end{props} \fi

\begin{proof} The second is immediate from corollary \ref{CDG1LIN} together with the observation that $\tau_{r,r}=1$ for any sign of $r.$

For the first, upon dividing by $\tau_{k+1,n}(a,aq^p),$ we see it is sufficient to prove
\bmlst
    \frac{ \sigma_{-(k+1), n}(a,aq^{p+1}) } { \tau_{k+1,n}(a,aq^p)} =  \left\{ \sigma_{-(k+1),k+1}(aq^p,aq^{p+1}) \right\} \cdot 1 \\
+\left\{  \tau_{-(k+1),-(k+1)}(aq^p,aq^{p+1})  \right\} \left\{ \frac{ \sigma_{-(k+1),n}(a,aq^p) } { \tau_{k+1,n}(a,aq^p) } \right\} \\
  + \left\{  \tau_{-(k+1),-(k+2)}(aq^p,aq^{p+1}) \right\} \left\{ \frac { \sigma_{-(k+2), n}(a,aq^p) } { \tau_{k+1,n}(a,aq^p) } \right\} .
\end{multline*}
That means we need to show
\bml
\label{EQN1c}\ 
\frac{abcd q^{n+k+p+1} (q^{k+1}, cdq^{k}, q^{-(p+1)}|q)_1} { (abcdq^{n+k}, abcdq^{n+k+p+1} |q)_1} \\
= 
   \left\{ \frac{(q| q)_{k+1} bcd  \lb aq^{p+1} \rb  q^{2k+1}  (cd q^{k }| q )_{1}  {\ds   (q^{-1} | q )_{1}  } }  { (q| q)_k ( \lb aq^{p} \rb bcdq^{2k+1} | q)_{1}  (bcd \lb aq^{p+1} \rb  q^{2k+1} | q)_{1} } \right\}  \cdot 
   1 \\
+  \left\{1 \right\}  \cdot 
\left\{ \frac{abcdq^{2k+p+1}(q^{k+1},q^{n-k-p-1},cdq^k|q)_1 } {(abcdq^{n+k},abcdq^{2k+p+1 }|q)_1 } \right\} \\
+  \left\{ \frac{(q^{2}| q)_k   \lb aq^{p+1} \rb  (bc q^{k+1}| q )_{1}  (bd q^{k+1}| q )_{1}   (cdq^k | q )_{1}  {\ds   (q^{-1}| q )_{1}  } }  { (q| q)_k  ( \lb aq^{p} \rb bcdq^{2(k+1)} | q)_{1}  (bcd \lb aq^{p+1} \rb  q^{2k+1} | q)_{1} } \right\} \cdot \\   
\left\{ \frac
{bcd q^{2(k+1)} (q^{n-k-1}, abcdq^{2k+p+2} |q)_1} 
{(bcq^{k+1},bdq^{k+1}, abcdq^{n+k+p+1} |q)_1} \right\} .
\end{multline}

We will eventually reduce this identity to the vanishing of a 1-variable polynomial in $y$ with coefficients in the field $\mathbb{Q}(a,b,c,d,u,v,w,q)$ 
with the property that when 
$$
y=abcd \hspace{5mm} u= q^n \hspace{5mm} v= q^k\hspace{5mm} w= q^p
$$
we obtain a unit in the coefficient field times the difference between the two sides of equation (\ref{EQN1c})  above. So, effectively, we can use the variables $y,u,v,w$
as abbreviations for the above expressions.

Multiplying (\ref{EQN1c})  by
$$
\frac{ ( q, abcdq^{n+k}, abcdq^{n+k+p+1}, abcdq^{ 2k+p+1}, abcdq^{2k+p+2} | q )_1 } {abcd (q^{k+1}, cdq^k | q)_1 } \\
$$
we see it is sufficient to show the vanishing of the  polynomial $p_1$ below.

$\big($To arrive at the final form of  $p_1,$ we use Lemma \ref{POCHNEG} above  to simplify the following expressions:
$$ 
(q^{-1}| q)_1,\  
  (q^{-(p+1)} | q )_{1}, 
\text{ and }
(q^{n-k-p-1}| q)_1. 
\big) 
$$


\bmlst
p_1  =  
-\left\{ (yq^{2k+p+1}| q)_1(yq^{2k+p+2}| q)_1  (q| q)_1 
 q^{n + k+p+1} (q^{-(p+1)} | q )_{1}    \right\}  \\
+  
  \left\{ (yq^{n+k}| q)_1 (yq^{n+k+p+1}| q)_1  (q| q)_1  q^{2k+p+2} (q^{-1}| q)_1    \right\}  \\  
   +  \left\{1 \right\}  \cdot    
\left\{ (yq^{n+k+p+1}| q)_1 (yq^{2k+p+2}| q)_1  (q| q)_1 
 q^{2k+p+1}  (q^{n-k-p-1}| q)_1    \right\} \\
+  \left\{ (yq^{n+k}| q)_1 (yq^{2k+p+1}| q)_1(q^{n-k-1}| q)_1 q^{2k+p+3}    
(q^{-1}| q)_1    \right\} \cdot 
\end{multline*}
\bmlst 
\hphantom{p_1}  =  
-\left\{ (yq^{2k+p+1}| q)_1(yq^{2k+p+2}| q)_1  (q| q)_1 
 q^{n + k+p+1} \lb - q^{-(p+1)}(q^{p+1}| q)_1  \rb   \right\}  \\
+  
     \left\{ (yq^{n+k}| q)_1 (yq^{n+k+p+1}| q)_1  (q| q)_1  q^{2k+p+2} \lb  - q^{-1}(q| q)_1 \rb   \right\}   \\  
  +  \left\{1 \right\}  \cdot   
\left\{ (yq^{n+k+p+1}| q)_1 (yq^{2k+p+2}| q)_1(q| q)_1 
 q^{2k+p+1}  \lb q^{-(k+p+1)}   \lb (q^{n}| q)_1-(q^{k+p+1}| q)_1 \rb \rb  \right\} \\
+  \left\{ (yq^{n+k}| q)_1 (yq^{2k+p+1}| q)_1 \lb  q^{-(k+1)} \lb (q^{n}| q)_1-(q^{k+1}| q)_1 \rb\ \rb  q^{2k+p+3} 
 \lb - q^{-1}(q| q)_1 \rb   \right\}   
\end{multline*}
\bmlst
\hphantom{p_1} =  
\left\{ (yq^{2k+p+1}| q)_1(yq^{2k+p+2}| q)_1 (q| q)_1 
 q^{n + k}   (q^{p+1}| q)_1     \right\}\\
 -   
     \left\{ (yq^{n+k}| q)_1 (yq^{n+k+p+1}| q)_1   \lb (q| q)_1 \rb^2  q^{2k+p+1}   \right\}   \hspace{50mm} \\  
  +  \left\{1 \right\}  \cdot   %
\left\{ (yq^{n+k+p+1}| q)_1 (yq^{2k+p+2}| q)_1 (q| q)_1 
 q^{k}  \lb  \lb (q^{n}| q)_1-(q^{k+p+1}| q)_1 \rb \rb  \right\} \\
- \left\{ (yq^{n+k}| q)_1 (yq^{2k+p+1}| q)_1 \lb  \lb (q^{n}| q)_1-(q^{k+1}| q)_1 \rb\ \rb q^{k+p+1}
 (q| q)_1   \right\} 
\end{multline*}
Using the notation $\big[a,b,c,\dots\big]=(a| q)_1(b| q)_1(c| q)_1\ldots,$ 
recalling our abbreviation variables $u,v,\text{ and } w,$
and multiplying by $v^{-1} (q| q)^{-1}$ 
it suffices to show the vanishing of 
\bmlst
p_2 =  u[ yv^2wq,yv^2wq^2,wq]
-vwq[yuv, yuvwq,q] \\
 +(vwq-u ) [yuvwq, yv^2wq^2 ] 
-wq(vq-u)[yuv, yv^2wq] 
\end{multline*}

This expression may be interpreted as a one variable polynomial of degree at most $2$ in $y$ with coefficients in the field $\mathbb{Q}(u,v,w,q).$
Evaluating
\begin{description}
\item[at $y=0$ ] 
$$
p_2(0)= u(1-wq) -vwq(1-q) +(vwq-u) -wq(vq-u)=0.
$$
\item[at $y=(v^2wq)^{-1}$ ] Note at this point $yv^2wq^2=q, \ yuvwq=uv^{-1}, \text{ and }yuv= u(vwq)^{-1}.\ $ So
\bmlst
p_2\big((v^2wq)^{-1} \big) =(1-uv^{-1}) \big\{  -vwq(1- u(vwq)^{-1}) (1-q) +(vwq-u)(1-q)\big\} \\
                                         =(1-uv^{-1})(1-q)(vwq-u)  \big\{ -1+1 \big\} =0.
\end{multline*}
\item[at $y=(uvwq)^{-1}$ ] Note at this point $yv^2wq=u^{-1}v \ yv^2wq^2= u^{-1}v q\  \text{ and }yuv=(wq)^{-1} $ So

\bmlst
p_2\big((uvwq)^{-1} \big)= (1-u^{-1}v)\big\{ u(1-u^{-1}vq)(1-wq)-wq(vq-u)(1-(wq)^{-1} \big\} \\
                                       = (1-u^{-1}v)(vq-u)(1-wq)\big\{-1 +1\big\} =0.
\end{multline*}
\end{description}

Thus $p_2$ being of degree at most  $2$ and vanishing at $3$ points implies $p_2$ is identically $0$ and the proposition proof is complete.

\end{proof}

\pagebreak
\bc \bf The $T^{11}$ Identity \ec
\ifJOLT \begin{Proposition} \else \begin{props} \fi
\label{T11PROP}
When $0 \le k \le n-1$
\bml
  \tau_{-(k+1),-(n+1)}(a,aq^{p+1}) = 
   \big[ \sigma_{-(k+1),k+1}(aq^p,aq^{p+1})\big] \big[ \sigma_{k+1, -(n+1)}(a,aq^p)\big] \\
+ \big[ \tau_{-(k+1),-(k+1)}(aq^p,aq^{p+1})\big] \big[ \tau_{-(k+1),-(n+1)}(a,aq^p) \big] \\
+ \big[ \tau_{-(k+1),-(k+2)}(aq^p,aq^{p+1})\big] \big[ \tau_{-(k+2), -(n+1)}(a,aq^p)\big].
\end{multline}
And 
$$  \tau_{-(n+1),-(n+1)}(a,aq^{p+1}) = \big[ \tau_{-(n+1),-(n+1)}(aq^p,aq^{p+1})\big] \big[ \tau_{-(n+1),-(n+1)}(a,aq^p) \big].$$
\ifJOLT \end{Proposition} \else \end{props} \fi

\begin{proof}
The second identity just says $1=1 \cdot 1.$

For the first, upon dividing by $ \sigma_{k+1,-(n+1)}(a,aq^p) ,$ we see it is sufficient to prove
\bmlst
  \frac{ \tau_{-(k+1),-(n+1)}(a,aq^{p+1}) } { \sigma_{k+1,-(n+1)}(a,aq^p) }  = 
   \left\{ \sigma_{-(k+1),k+1}(aq^p,aq^{p+1}) \right\} \cdot 1  \\
+ \left\{ \tau_{-(k+1),-(k+1)}(aq^p,aq^{p+1})\right\} \left\{  \frac{ \tau_{-(k+1),-(n+1)}(a,aq^p) } { \sigma_{k+1,-(n+1)}(a,aq^p) } \right\} \\
+ \left\{ \tau_{-(k+1),-(k+2)}(aq^p,aq^{p+1})\right\} \left\{  \frac{ \tau_{-(k+2), -(n+1)}(a,aq^p) }{ \sigma_{k+1,-(n+1)}(a,aq^p) } \right\}.
\end{multline*}
That means we need to show
\bml
\label{EQN1d}
\frac{q^{n-k} (q^{k+1}, cdq^k, q^{-(p+1)}|q)_1} {(q^{n-k},q^{n-k-p-1}|q)_1 } \\
= 
   \left\{ \frac{(q| q)_{k+1} bcd  \lb aq^{p+1} \rb   q^{2k+1}  (cd q^{k }| q )_{1}  {\ds   (q^{-1} | q )_{1}  } }  { (q| q)_k ( \lb aq^{p} \rb bcdq^{2k+1} | q)_{1}  (bcd \lb aq^{p+1} \rb  q^{2k+1} | q)_{1} }  \right\}  \cdot 1 \\
+  \left\{1 \right\} 
\left\{ \frac{ (q^{k+1}, cdq^k, abcdq^{n+k+p+1}|q)_1} {(q^{n-k}, abcdq^{2k+p+1 }|q)_1 } \right\} \\
+  \left\{ \frac{(q^{2}| q)_k   \lb aq^{p+1} \rb  (bc q^{k+1}| q )_{1}  (bd q^{k+1}| q )_{1}   (cdq^k | q )_{1}  {\ds   (q^{-1} | q )_{1}  } }  { (q| q)_k  ( \lb aq^{p} \rb bcdq^{2(k+1)} | q)_{1}  (bcd \lb aq^{p+1} \rb  q^{2k+1} | q)_{1} } \right\} \cdot      \\
   \left\{ \frac
  {(abcdq^{n+k+1},abcdq^{2k+p+2 }|q)_1 } 
  {aq^p(bcq^{k+1},bdq^{k+1}, q^{n-k-p-1}|q)_1 } \right\}
\end{multline}

Multiplying (\ref{EQN1d}) by
$$
 \frac{q^{k+p+1} ( q, q^{n-k}, q^{n-k-p-1}, abcdq^{2k+p+1},  abcdq^{2k+p+2} | q)_1} { (q^{k+1}, cdq^k | q)_1} \\
$$
we see it is sufficient to show the vanishing of the  polynomial $p_1$ below.

$\big($To arrive at the final form of $p_1,$ we use Lemma \ref{POCHNEG} above  to simplify the following expressions:
$$ 
(q^{-1}| q)_1,\ 
(q^{-p}| q)_1, \ 
(q^{-(p+1)}| q)_1,\  
 \text{and }
q^{n-k-p-1}| q)_1 \ 
\big) 
$$

We will eventually reduce this identity to the vanishing of a 1-variable polynomial in $y$ with coefficients in the field $\mathbb{Q}(a,b,c,d,u,v,w,q)$ 
with the property that when 
$$
y=abcd \hspace{5mm} u= q^n \hspace{5mm} v= q^k\hspace{5mm} w= q^p
$$
we obtain a unit in the coefficient field times the difference between the two sides of equation (\ref{EQN1d}) above. So, effectively, we can use the variables $y,u,v,w$
as abbreviations for the above expressions.


\bmlst
p_1=
-\left\{(yq^{2k+p+1} | q)_1 (yq^{2k+p+2} | q)_1  (q| q)_1   q^{n+p+1}     (q^{-(p+1)} | q )_{1}    \right\} \\
 +
   \left\{    (q^{n-k}| q)_{1}  (q| q)_1 y q^{3k+2p+3}     (q^{-1} | q )_{1}     
   (q^{n-k-p-1} | q )_{1}  
   \right\} \\
  +  \left\{1 \right\} 
   \left\{  (yq^{n+k+p+1} | q)_1  (yq^{2k+p+2} | q)_1(q| q)_1  q^{k+p+1}  (q^{n-k-p-1} | q )_{1}   
   \right\} \\
 +  \left\{  (yq^{n+k+1} | q)_1  (yq^{2k+p+1} | q)_1  (q^{n-k}| q)_1   (q^{-1} | q )_{1}        
 q^{k+p+2} 
  \right\}  
 \end{multline*} 
\bmlst
\hphantom{p_1}=
-\left\{(yq^{2k+p+1} | q)_1 (yq^{2k+p+2} | q)_1  (q| q)_1   q^{n+p+1}    \lb  - q^{-(p+1)}(q^{p+1}| q)_1 \rb    \right.   \\
 +
   \left\{    \lb  q^{-k} \lb (q^{n}| q)_1-(q^k| q)_1 \rb\rb   (q| q)_1   y  q^{3k+2p+3}     \lb  - q^{-1}(q| q)_1  \rb  \right.  \cdot \\
   \left. \lb q^{-k-p-1} \lb (q^n| q)_1 - (q^{k+p+1}| q)_1 \rb  \rb   \right\} \\
  +  \left\{1 \right\} 
   \left\{  (yq^{n+k+p+1} | q)_1  (yq^{2k+p+2} | q)_1(q| q)_1     q^{k+p+1}  \lb q^{-k-p-1} \lb (q^n| q)_1 - (q^{k+p+1}| q)_1 \rb  \rb         \right\} \\
+  \left\{  (yq^{n+k+1} | q)_1  (yq^{2k+p+1} | q)_1   \lb q^{-k} \lb (q^{n}| q)_1-(q^k| q)_1 \rb \rb      \lb  - q^{-1}(q| q)_1  \rb     
   q^{k+p+2} 
  \right\}  
 \end{multline*}
\bmlst
\hphantom{p_1}=
\left\{(yq^{2k+p+1} | q)_1 (yq^{2k+p+2} | q)_1  (q| q)_1   q^{n}      (q^{p+1}| q)_1      \right\} \\
 - \left\{     \lb (q^{n}| q)_1-(q^k| q)_1 \rb      y  q^{k+p+1}    \lb (q| q)_1  \rb^2
   \lb (q^n| q)_1 - (q^{k+p+1}| q)_1 \rb        \right\} \\
+ \left\{1 \right\} 
   \left\{  (yq^{n+k+p+1} | q)_1  (yq^{2k+p+2} | q)_1(q| q)_1       \lb (q^n| q)_1 - (q^{k+p+1}| q)_1 \rb          \right\} \\
 -  \left\{  (yq^{n+k+1} | q)_1  (yq^{2k+p+1} | q)_1  \lb (q^{n}| q)_1-(q^k| q)_1 \rb    (q| q)_1      
   q^{p+1} 
  \right\}  
 \end{multline*}
Using the notation $\big[a,b,c,\dots\big]=(a| q)_1(b| q)_1(c| q)_1\ldots,$ 
recalling our abbreviation variables $u,v,\text{ and } w,$
and multiplying by $\lb (q| q)_1 \rb^{-1}$ 
it suffices to show the vanishing of 
\bmlst
p_2 = u[yv^2wq, yv^2wq^2, wq ] 
- yvwq(v-u)(vwq-u)(q| q)_1\\
+ (vwq-u) [yuvwq,yv^2wq^2] 
 -wq(v-u) [yuvq,yv^2wq] 
\end{multline*}
This expression may be interpreted as a one variable polynomial of degree at most $2$ in $y$ with coefficients in the field $\mathbb{Q}(u,v,w,q).$
Evaluating 
\begin{description}
\item[at $y=0$ ] 
$$
p_2(0)= u(1-wq) + (vwq-u)  -wq(v-u)=0.
$$
\item[at $y=(v^2wq)^{-1}$ ] Note at this point $yv^2wq^2=q, \ yvwq=v^{-1}, \text{ and } yuvwq=uv^{-1}.$ So
\bmlst
p_2\big((v^2wq)^{-1} \big)=(vwq-u) \big\{ -v^{-1}(v-u)(1-q) + (1-uv^{-1})(1-q)  \big\} \\
                                         =(vwq-u) (1-uv^{-1})(1-q)  \big\{ -1+1 \big\} =0.
\end{multline*}
\item[at $y=(v^2wq^2)^{-1}$ ] Note at this point $yv^2wq=q^{-1}, \ yvwq=(vq)^{-1}, \text{ and } yuvq=u(vwq)^{-1}.$ So
\bmlst
p_2\big((v^2wq^2)^{-1} \big)=(v-u)  \big\{-(vq)^{-1}(vwq-u)(1-q) - wq(1- u(vwq)^{-1})(1-q^{-1}) \big\} \\
                                       =(v-u) (1-q)\big\{ -( w -u (vq)^{-1}  ) + (w- u(vq)^{-1})\big\}=0. 
\end{multline*}
\end{description}
Thus $p_2$ being of degree $2$ and vanishing at $3$ points implies $p_2$ is identically $0$ and the proposition proof is complete.

\end{proof}

The combination of propositions \ref{PLANA1}, \ref{CDG1AEQ}, \ref{CDG2AEQ}, \ref{T00PROP}, \ref{T01PROP}, \ref{T10PROP}, and \ref{T11PROP}   finishes the proof of all $3$ steps of (\ref{PLANA}), Proof Plan A, and so completes the proof of theorem \ref{ETC}. 

\section{Reproof of the Askey-Wilson Result Using the Nonsymmetric Version \label{REPROOF}}

Having given a direct proof of the shift-a connection coefficient formula Theorem \ref{ETC}, we now show that it can be used to give another
proof of the original Askey-Wilson result Theorem \ref{PTCFRM}. 

For this purpose, consistent with $x=\cos{\theta}, z=e^{i\theta}$, so $2x=\lp z +z^{-1}\rp,$ we view the $n$'th Askey-Wilson polynomial $P_n$ 
as a zig-zag monic Laurent polynomial in $z.$ So as an ordinary polynomial in $x,$ the leading coefficient would be $2^{n}.$

As earlier for the $E_r,$ we will shorten $P_n(z; a,b,c,d|q)$ to $P_n(a).$

First note that the usual DAHA relation on $T_1$ 
$$
T_1-T_1^{-1} = t_1^{\frac{1}{2}} - t_1^{-\frac{1}{2}}
$$
translates to the quadratic relation 
\begin{equation}
\lp T_1+ t_1^{-\frac{1}{2}} \rp \lp T_1- t_1^{\frac{1}{2}} \rp = 0.
\end{equation}
showing that the possible eigenvalues of $T_1$ are $ -t_1^{-\frac{1}{2}} \text{ and } t_1^{\frac{1}{2}}.$

Define
\begin{eqnarray}
L_1 = T_1+t_1^{-\frac{1}{2}}& & L_2= -T_1 +t_1^{\frac{1}{2}}.
\end{eqnarray}
Definition 1.2 in \cite{ASKVOL} defined the Askey-Wilson $P_n$ (respectively $Q_n$) as normalized multiples of $L_1E_n$ (respectively
$L_2 E_n$) for $n \ge 0.$ Theorem 1.3 there (referring to \cite{ANNALS}) pointed out that this definition of $P_n$ agrees with the usual $P_n.$

It was pointed out in Theorem 1.3 of \cite{ASKVOL} that, up to a normalizing factor, $P_{|n|}$ is $(T_1+t_1^{-\frac{1}{2}})E_n$ for $n \ge 0.$

A similar fact holds when $E_n$ is replaced by $E_{-n}$ for $n > 0.$ The reason for this is that by the recursion (\ref{SRCN0}) and 
(\ref{SRCN1}), (originally
proved in \cite{ASKVOL}) for $n>0,$ the $2$-dimensional subspace spanned by $E_{-n} \text{ and } E_n$ is invariant under $T_1.$ So
$L_1+L_2$ is a multiple of the identity on these $2$-dimensional subspaces and easy DAHA calculations show:
\ifJOLT \begin{Lemma} \else \begin{lemmas} \fi 
\label{TWODPROJCD} Up to scalar factors, the operators $L_1$ and $L_2$ are algebraically orthogonal projections onto $1$-dimensional
subspaces of the  $2$-dimensional subspace (for $n \ge0$)  spanned by $E_{-(n+1)} \text{ and } E_{n+1}.$ In fact
\begin{enumerate}
\item $L_1L_2=L_2 L_1=0$
\item $L_1^2=( t_1^{\frac{1}{2}}+ t_1^{-\frac{1}{2}} )L_1$
\item $L_2^2=( t_1^{\frac{1}{2}}+ t_1^{-\frac{1}{2}} )L_2$
\end{enumerate}
\ifJOLT \end{Lemma} \else \end{lemmas} \fi 
{\bf Remark: } By  algebraically orthogonal projections, we are referring to endomorphisms $\pi_1$ and $\pi_2=\text{Identity} -\pi_1$ satisfying
$\pi_1^2=\pi_1.$ This of course also implies $\pi_2^2=\pi_2$ and $\pi_1\pi_2 =\pi_2\pi_1=0.$ These projections can also  be related to the natural
inner product specified in Definition 1.4 of \cite{ASKVOL} but we omit the details here.
\begin{proof}
\begin{enumerate}
\item By the basic DAHA relation
$$ L_1= T_1+t_1^{-\frac{1}{2}} = T_1^{-1} +t_1^{\frac{1}{2}}.$$
Then
\bmlst
L_1L_2= \lp T_1^{-1} +t_1^{\frac{1}{2}} \rp \lp  -T_1 +t_1^{\frac{1}{2}} \rp \\
= -1 -  t_1^{\frac{1}{2}}  \lp T_1-T_1^{-1} \rp +t_1 \\
= -1 -  t_1^{\frac{1}{2}}  \lp  t_1^{\frac{1}{2}} -  t_1^{-\frac{1}{2}}  \rp +t_1 =0.
\end{multline*}
And $L_2L_1=L_1L_2.$
\item $\ds L_1^2=L_1(L_1+L_2)= \lp t_1^{\frac{1}{2}} + t_1^{-\frac{1}{2}} \rp L_1.$
\item $\ds L_2^2=(L_1+L_2)L_2= \lp t_1^{\frac{1}{2}} + t_1^{-\frac{1}{2}} \rp L_2.$
\end{enumerate}
\end{proof}

Consequently, for $n>0,$ $P_n$ is also a multiple of $\lp T_1+t_1^{-\frac{1}{2}} \rp E_{-n}.$ (The operator $L_1,$ up to constant factors, 
is sometimes referred to as {\em Hecke-symmetrization}.) Translating to the operator $\widetilde{T}_1= t_1^{\frac{1}{2}} T_1 $ which we have mostly been using (and which
avoids square roots), we see that $P_n$ for $n>0$ is a multiple of $(\widetilde{T}_1 +1)E_{-n}.$

With the zig-zag monic normalization conventions, Proposition \ref{SCALE} and equation
(\ref{SRCNX}),
give us more explicitly for $n \ge 0$ 

\begin{equation}
\label{T1ECD} 
\widetilde{T}_1 E_{-(n+1)} =\lb  \hat{c}_{n+1} \rb^{-1} \Big\{E_{n+1} -  \hat{d}_{n+1} E_{-(n+1)} \Big\}.
\end{equation}
where 
$$
 \hat{c}_{n+1} = - \frac{1}{ab} =t_1^{-1} 
 \hspace{15mm} 
\hat{d}_{n+1} = \ds - \frac{( abq^{n+1} | q )_1+ab(cdq^{n}| q)_1}{ab(1-abcdq^{2n+1})} 
$$

\ifJOLT \begin{Proposition} \else \begin{props} \fi
\label{HECKESYMCD}
For $n \ge 0,$
$$ P_{n+1}= t_1^{-1} (\widetilde{T}_1 +1)E_{-(n+1)}.$$
(Also $P_0=E_0=1.)$

More explicitly
$$
P_{n+1}(a)  = E_{n+1} (a) + \lb \gamma_{n+1} (a) \rb E_{-(n+1)}(a)
$$
where $\gamma_{n+1} (a)$ (really a function of $a,b,c,d,\text{ and } q$)   is the scalar defined by
\begin{equation}
\label{GAMMACD}
\gamma_{n+1} (a) = \frac{(q^{n+1},cdq^n|q)_1}{(abcdq^{2n+1}|q)_1}.
\end{equation}
\ifJOLT \end{Proposition} \else \end{props} \fi
\begin{proof}
We've already argued  $t_1^{-1} (\widetilde{T}_1 +1)E_{-(n+1)}$ is a multiple of $P_{n+1}.$ Equation (\ref{T1ECD}) implies $t_1^{-1} \lb  \hat{c}_{n+1} \rb^{-1}=1.$ And $P_{n+1} =E_{n+1} = f_{n+1} \text{ mod } \mathcal{R}_{-(n+1)}$ while $E_{-(n+1)} \in \mathcal{R}_{-(n+1)}.$ So $t_1^{-1} (\widetilde{T}_1 +1)E_{-(n+1)}$ gives $P_{n+1}$ exactly.

Remembering $t_1\hat{c}_{n+1}=1$ and  plugging the formulas of (\ref{T1ECD}) into $P_{n+1}= t_1^{-1} (\widetilde{T}_1 +1)E_{-(n+1)}$ gives us
$$
P_{n+1} = E_{n+1} + \lp  \hat{c}_{n+1} -  \hat{d}_{n+1}  \rp E_{-(n+1)}.
$$
(Since our normalization convention is that both $P_{n+1}$ and $E_{n+1}$ are zig-zag monic, the coefficient $1$ above  in front of $E_{n+1}$ was
known in advance.)
Our definition of $\gamma_{n+1} (a)$ is just a simplification of $\hat{c}_{n+1} -  \hat{d}_{n+1}$ as demonstrated in the easy Lemma \ref{CMDSIMPCD} below.
\end{proof}

\ifJOLT \begin{Lemma} \else \begin{lemmas} \fi For $n \ge 0$
\label{CMDSIMPCD}
$$
\hat{c}_{n+1} -  \hat{d}_{n+1}  = \frac{(q^{n+1},cdq^n|q)_1}{(abcdq^{2n+1}|q)_1},
$$
this quantity being the $ \gamma_{n+1} (a)$ (really $ \gamma_{n+1} (a,b,c,d|q)$) defined in Proposition \ref{HECKESYMCD}.
\ifJOLT \end{Lemma} \else \end{lemmas} \fi 
\begin{proof}
\bmlst
\hat{c}_{n+1} -  \hat{d}_{n+1}  =-\frac{1}{ab} + \frac{( abq^{n+1} | q )_1+ab(cdq^{n}| q)_1}{ab(1-abcdq^{2n+1})}  \\
= -\lp\frac{1} {ab(1-abcdq^{2n+1})} \rp  \big\{ 1 - abcdq^{2n+1} + (abq^{n+1} - 1) +ab(cdq^n-1)\big\}\\
=\frac{(1-q^{n+1})(1-cdq^n)} {1-abcdq^{2n+1}}.\\
\end{multline*}
\end{proof}

Now that we know
\begin{eqnarray*}
P_{n+1}(a) &=& E_{n+1}(a) + \lb   \gamma_{n+1}(a) \rb \lb E_{-(n+1)}(a) \rb \\
P_{n+1}(e) &=& E_{n+1}(e) + \lb   \gamma_{n+1}(e) \rb \lb E_{-(n+1)}(e) \rb, 
\end{eqnarray*}
a natural way to obtain the Askey-Wilson connection coefficient relation Theorem \ref{PTCFRM} is:

\begin{equation} 
\label{PLANPCD}
\text{ \em (PROOF PLAN B)  }
\end{equation}
\begin{enumerate}
\item Start with the top line for $P_{n+1}(a).$
\item Apply Theorem \ref{ETC} to express each of $E_{\pm(n+1)}(a)$ in terms of the $E_r(e).$ (Here $r$ can be of any sign.) 
\item Show that the combinations of $E_{\pm(m+1)}(e)$ (for $m\ge 0$) which result are in fact the $c_{m+1,n+1}P_{m+1}(e)$ of the
Askey-Wilson result. (As well as the $E_0(e)$ coefficient matching $c_{0,n+1}.$)
\end{enumerate}

We will now show that 
{\em PROOF PLAN B} 
can be carried out to prove the Askey-Wilson result.

Recall our generic notation for the $E$ connection coefficient relations.
For  $n \geq 0:$ 
\begin{equation}
\label{EGENERICPLUSCD0} E_n(a) = \sum_{m=0}^n \lb \tau_{m,n} \rb E_m(e)  + \sum_{m=0}^{n-1} \lb \sigma_{-(m+1),n} \rb E_{-(m+1)}(e). \end{equation} 
\begin{equation}
\label{EGENERICMINUSCD0} E_{-(n+1)}(a) =  \sum_{m=0}^n \lb \tau_{-(m+1),-(n+1)}\rb E_{-(m+1)}(e) + \sum_{m=0}^n \lb \sigma_{m,-(n+1)}\rb E_{m}(e) . 
\end{equation} 

Our original formulation of Theorem \ref{ETC}, introduced the notation $d_{r,s}$ which is related to the (variously subscripted) $\tau,\sigma$ by:
\begin{equation}
d_{r,s} = \begin{cases}
\tau_{r,s}/c_{|r|,|s|} & \text{if } (r \ge 0 \text{ and } s \ge 0)  \text{ or } (r <  0 \text{ and } s < 0)\\
\sigma_{r,s}/c_{|r|,|s|} &  \text{if } (r \ge 0 \text{ and } s < 0)  \text{ or } (r <  0 \text{ and } s > 0). \\
\end{cases}
\end{equation}
We return to that notation now and write:
\begin{equation}
\label{EGENERICPLUSCD} E_n(a) = \sum_{m=0}^n \lb d_{m,n} c_{m,n}  \rb E_m(e)  + \sum_{m=0}^{n-1} \lb d_{-(m+1),n} c_{m+1,n}\rb E_{-(m+1)}(e).
\end{equation}
\begin{equation}  
\label{EGENERICMINUSCD} E_{-(n+1)}(a) = \sum_{m=0}^n \lb d_{-(m+1),-(n+1)} c_{m+1, n+1} \rb E_{-(m+1)}(e) + \sum_{m=0}^n \lb d_{m,-(n+1)} c_{m,n+1}  \rb E_{m}(e).
\end{equation}


Plugging these two into 
$$ P_{n+1}(a) = E_{n+1}(a) + \lb   \gamma_{n+1}(a) \rb \lb E_{-(n+1)}(a) \rb $$
and considering the coefficients of $E_{\pm(m+1)}(e)$ which result makes clear the relevance of the following two propositions.
\ifJOLT \begin{Proposition} \else \begin{props} \fi For $0 \le m \le n$
\label{EGAMMAEPLUSCD}
\begin{equation}
d_{m+1,n+1} + \gamma_{n+1}(a) d_{m+1,-(n+1)} =1.
\end{equation}
\ifJOLT \end{Proposition} \else \end{props} \fi
\begin{proof}
\bmlst
d_{m+1,n+1} + \gamma_{n+1}(a) d_{m+1,-(n+1)} \\
=\Big\{ \ds  \frac{q^{n-m} (abcdq^{n+m+1}| q)_1}{(abcdq^{2n+1}| q)_1}\Big\} 
+\Big\{ \ds \frac{(q^{n+1},cdq^n|q)_1}{(abcdq^{2n+1}|q)_1} \Big\} \Big\{  \frac{ (q^{n-m }|q)_{ 1 } } {( q^{n+1}, cdq^{n }|q)_{ 1 }}\Big\} \\
\lp \frac{1}{(abcdq^{2n+1}|q)_1 } \rp \big\{q^{n-m} - abcdq^{ 2n+1} +1 - q^{n-m}\big\} =1.\\
\end{multline*}
\end{proof}

\ifJOLT \begin{Proposition} \else \begin{props} \fi
\label{EGAMMAEMINUSCD} For $0 \le m \le n$
\begin{equation}
d_{-(m+1),n+1} + \gamma_{n+1}(a) d_{-(m+1),-(n+1)} =\gamma_{m+1}(e).
\end{equation}
\ifJOLT \end{Proposition} \else \end{props} \fi 
\begin{proof}
\bmlst
d_{-(m+1),n+1} + \gamma_{n+1}(a) d_{-(m+1),-(n+1)} \\
= \Big\{ \ds  \frac{bcdeq^{n+m+1} (q^{m+1}, cdq^{m}, ae^{-1}q^{n-m}|q)_{1 } } { (abcdq^{2n+1}, bcdeq^{2m+1 }|q)_{ 1}} \Big\} \\
+ \Big\{ \ds \frac{(q^{n+1},cdq^n|q)_1}{(abcdq^{2n+1}|q)_1}  \Big\}  \Big\{ \ds  \frac{ (q^{m+1}, cdq^{m}, bcdeq^{n+m+1 }|q)_{1 } } { (q^{n+1}, cdq^{n}, bcdeq^{2m+1}|q)_{ 1} }  \Big\}  \\
=\Big\{ \frac{(q^{m+1}, cdq^{m} |q)_{1 } } { (abcdq^{2n+1}, bcdeq^{2m+1 }|q)_{ 1}} \Big\} \big\{ bcdeq^{n+m+1}(ae^{-1}q^{n-m}|q)_{1 } +(bcdeq^{n+m+1 }|q)_{ 1} \big\}\\
=\Big\{ \frac{(q^{m+1}, cdq^{m} |q)_{1 } } { (abcdq^{2n+1}, bcdeq^{2m+1 }|q)_{ 1}} \Big\} \big\{ bcdeq^{n+m+1} - abcd q^{2n+1} + 1- bcdeq^{n+m+1 } \big\}\\
=  \frac{(q^{m+1}, cdq^{m} |q)_{1 } } { ( bcdeq^{2m+1 }|q)_{ 1}}=\gamma_{m+1}(e) \\
\end{multline*}
as asserted.
\end{proof}

Now we prove the Askey-Wilson connection coefficient result:
\begin{proof} (of Theorem \ref{PTCFRM})

The case of $n=0$ checks since $P_0$ is $1,$ independent of parameters, and $c_{0,0}=1.$

For $n \ge 0,$ using equations (\ref{EGENERICPLUSCD}) and (\ref{EGENERICMINUSCD}) as well as Propositions \ref{EGAMMAEPLUSCD} and \ref{EGAMMAEMINUSCD}
\bmlst
P_{n+1}(a)  = E_{n+1} (a) + \lb \gamma_{n+1} (a) \rb E_{-(n+1)}(a) \\
=  \sum_{m=0}^{n+1} \lb d_{m,n+1} c_{m,n+1} \rb E_m(e)  + \sum_{m=0}^{n} \lb d_{-(m+1),n+1} c_{m+1,n+1} \rb E_{-(m+1)}(e) \\
+ \lb \gamma_{n+1} (a) \rb  \Big\{ \sum_{m=0}^n \lb d_{-(m+1),-(n+1)} c_{m+1,n+1} \rb E_{-(m+1)}(e) + \sum_{m=0}^n \lb d_{m,-(n+1)}  c_{m,n+1}\ \rb E_{m}(e) \Big\}\\
\end{multline*}
\bmlst
= c_{0,n+1} \Big\{  d_{0,n+1}  E_0(e) +  \gamma_{n+1} (a)  d_{0,-(n+1)} \Big\} E_{0}(e) \\
+\sum_{m=1}^{n} c_{m,n+1}  \big\{  d_{m,n+1}  + \gamma_{n+1}(a) d_{m,-(n+1)} \big\} E_{m}(e) + 1 \cdot E_{n+1}(e)\\
+\sum_{m=0}^n c_{m+1,n+1} \big\{ d_{-(m+1),n+1}  + \gamma_{n+1}(a) d_{-(m+1),-(n+1)}  \big\} E_{-(m+1)}(e)  \\
\end{multline*}
\bmlst
= \text{ (by Propositions \ref{EGAMMAEPLUSCD} and \ref{EGAMMAEMINUSCD} }) \\
= c_{0,n+1} \cdot 1 \cdot E_{0}(e) 
+\sum_{m=0}^{n-1} c_{m+1,n+1} \cdot 1 \cdot E_{m+1}(e) + c_{n+1,n+1} E_{n+1}(e)\\
+\sum_{m=0}^{n-1} c_{m+1,n+1} \cdot \gamma_{m+1}(e) \cdot E_{-(m+1)}(e) + c_{n+1,n+1} \cdot \gamma_{n+1}(e) \cdot E_{-(m+1)}(e) \\
= \sum_{m=0}^{n+1} c_{m,n+1} P_m \\
\end{multline*}
completing the re-proof of the Askey-Wilson Theorem \ref{PTCFRM}.
\end{proof}

%
%

\n {\bf Remark:} Instead of proving Lemma \ref{TWODPROJCD} and Proposition \ref{HECKESYMCD}, we could have based those
aspects of our proof on the discussion in Section 3 of \cite{ZHED}. In particular, that reference clearly explains how, for $n>0,$ using the alternate 
choice
$$D=\widetilde{Y}+q^{-1}abcd \widetilde{Y}^{-1} \hspace{5mm} \text{(with } \widetilde{Y}= \widetilde{T}_1\widetilde{T}_0) $$
for second order operator with Askey-Wilson polynomials $P_n$ as eigenfunctions, the corresponding eigenspace of $D$ is $4$ dimensional and spanned by
$P_n,Q_n,E_n,\text{ and } E_{-n}.$ Moreover $D$ commutes with both $\widetilde{T}_1$ and $\widetilde{T}_0,$ with $P_n$ and $Q_n$
being eigenfunctions of $\widetilde{T}_1.$ The respective eigenvalues are $t_1$ and $-1.$ Exact formulas expressing $E_{\pm n}$
as linear combinations of $P_n$ and the eigenvalue  $-1$ eigenfunction $Q_n^\dagger$ of $\widetilde{T}_0$ are also written down there.

\section{Appendix A on Change of Basis Conventions \label{APPA} }

Let $\mathcal{B}=\{e_{\alpha}\}_{\alpha \in \mathcal{I}}$ be an ordered basis of a vector space $\mathcal{V}$ where the index
set $\mathcal{I}=0,1,2,\ldots\ .$ For $x \in \mathcal{V},$ we represent the linear combination  $x=\sum_{\alpha \in \mathcal{I}} v^{\alpha} e_{\alpha}$
by the column vector
$$
[x]_{\mathcal{B}} = \begin{bmatrix} v^0 \\ v^1 \\ \vdots \\ \end{bmatrix}.
$$

If $\overline{\mathcal{B}}=\{\bar{e}_{\beta}\}_{\beta \in \mathcal{I}}$ is another ordered basis, then 
$$
[x]_{\overline{\mathcal{B}}} = \begin{bmatrix} \bar{v}^0 \\ \bar{v}^1 \\ \vdots \\ \end{bmatrix}.
$$
is the column vector corresponding to the linear combination $x= \sum_{\beta \in \mathcal{I}} \bar{v}^{\beta} \bar{e}_{\beta}.$

Then a change of basis relationship
$$
e_{\alpha}=\sum_{\beta \in \mathcal{I}} T_{\beta\alpha}  \bar{e}_{\beta}
$$
corresponds to
$$
 \begin{bmatrix} \bar{v}^0 \\ \bar{v}^1 \\ \vdots \\ \end{bmatrix} = T \begin{bmatrix} v^0 \\ v^1 \\ \vdots \\ \end{bmatrix}
$$
where the row $\beta$ column $\alpha$ entry of $T$ is $T_{\beta\alpha}.$

This is easily confirmed by noting
$$
x=
 \sum_{\alpha} v^{\alpha} e_{\alpha} =
 \sum_{\alpha,\beta} v^{\alpha} T_{\beta\alpha} \bar{e}_{\beta}=
 \sum_{\beta} \Big(\sum_{\alpha} T_{\beta\alpha} v^{\alpha} \Big)  \bar{e}_{\beta} =
 \sum_{\beta} \bar{v}^{\beta} \bar{e}_{\beta}=
 x.
$$

In the case 
\bmlst
n\geq 0: \  E_n(z; a,b,c,d| q) = \sum_{m=0}^n\tau_{m,n} E_m(z; e,b,c,d| q) \\+ \sum_{m=0}^{n-1}\sigma_{-(m+1),n} E_{-(m+1)}(z; e,b,c,d| q) 
\end{multline*}
\bmlst
n\geq0: \  E_{-(n+1)}(z; a,b,c,d| q) = \sum_{m=0}^{n}\tau_{-(m+1),-(n+1)} E_{-(m+1)}(z; e,b,c,d| q) \\+ \sum_{m=0}^{n}\sigma_{m,-(n+1)} E_{m}(z; e,b,c,d| q) 
\end{multline*}
we are viewing $E_0,E_1,\ldots$ as an ordered basis for $\mathcal{R}^0$ and $E_{-1},E_{-2}, \ldots$ as an ordered basis for $\mathcal{R}^1$

Note that $E_{-1}$ is in the initial position (index $0$) of the second list, in conformity with viewing it as $E_{-(n+1)}$ for $n=0.$

Thinking of parameter $e$ as giving rise to the $\bar{e}_{\beta}$ basis and parameter $a$ the $e_{\alpha}$ basis, then our generic
change of basis formula
$$
e_{\alpha}=\sum_{\beta } T_{\beta\alpha}  \bar{e}_{\beta}
$$
 becomes, in block form,
 
 $$
 T=\begin{bmatrix}
T ^{00} & T^{01} \\
T^{10} & T^{11} \\
\end{bmatrix}
$$
 so that {\em for  column vectors of (blocks of length $n$ or $n+1$) components} relative to these bases
 $$
  \begin{bmatrix} \bar{v}^0 \\ \bar{v}^1 \\  \end{bmatrix} = T \begin{bmatrix} v^0 \\ v^1  \\ \end{bmatrix}=\begin{bmatrix}
T ^{00} & T^{01} \\
T^{10} & T^{11} \\
\end{bmatrix}
\begin{bmatrix} v^0 \\ v^1  \\ \end{bmatrix}.
 $$
 
 Thus for $n\ge 0$ 
 \bmlst
 E_n(z; a,b,c,d| q) = \sum_{m=0}^n\tau_{m,n} E_m(z; e,b,c,d | q) \\+ \sum_{m=0}^{n-1}\sigma_{-(m+1),n} E_{-(m+1)}(z; e,b,c,d| q) 
\end{multline*}
corresponds to 
\bmlst
E_n(z; a,b,c,d| q) = \sum_{m=0}^nT^{00}_{m,n} E_m(z; e,b,c,d| q) \\+ \sum_{m=0}^{n-1}T^{10}_{m,n} E_{-(m+1)}(z; e,b,c,d| q) 
\end{multline*}
(Note the second matrix entry above really is $T^{10}_{m,n} $ (rather than $T^{01}_{m,n})$ in line with the general property (of our
conventions) that the first  $n+1$ {\em columns} of $T$ are expressing the decomposition of the first $n+1$ vectors $E_k(z; a,b,c,d| q) $
$(0 \le k \le n)$ as linear combinations of the ordered basis  
\bmlst
E_0(z; e,b,c,d| q), E_1(z; e,b,c,d| q),\ldots, \\ E_n(z; e,b,c,d| q), E_{-1} (z; e,b,c,d| q) ,\ldots,  E_{-n}(z; e,b,c,d| q).\big)
\end{multline*}

Similarly
 \bmlst
E_{-(n+1)}(z; a,b,c,d| q) = \sum_{m=0}^{n}\tau_{-(m+1),-(n+1)} E_{-(m+1)}(z; e,b,c,d| q) \\+ \sum_{m=0}^{n}\sigma_{m,-(n+1)} E_{m}(z; e,b,c,d| q) \\
\end{multline*}
corresponds to 
\bmlst
 E_{-(n+1)}(z; a,b,c,d| q) = \sum_{m=0}^{n}T^{11}_{m,n} E_{-(m+1)}(z; e,b,c,d| q) \\ + \sum_{m=0}^{n}T^{01}_{m,n} E_{m}(z; e,b,c,d| q)
\end{multline*}

These tell us that the row index $m$ and column index $n$ entries of the four matrices are given by
\begin{eqnarray*}
T^{00}_{m,n} = \tau_{m,n} & \hspace{15mm} & T^{01}_{m,n} = \sigma_{m,-(n+1)} \\
T^{10}_{m,n} = \sigma_{-(m+1),n} &  \hspace{15mm} & T^{11}_{m,n} = \tau_{-(m+1),-(n+1)}. \\
\end{eqnarray*}

\section{Appendix B \label{HAT} }

In this paper, we only use the zig-zag increasing cases in the following table and so have  included just the proofs of those. We mention the others, which we
have also proven, because they may be of interest.
\begin{table}[htp] 
\caption{Summary of $\hat{a}_m,\hat{b}_m,\hat{c}_m,\hat{d}_m,$  \ $m =n \text{ or } -(n+1),  \  n \geq 0.$}
\begin{tabular}{|c|c|c|} \hline
 $n \geq 0$ {\bf (Non-Negative Cases)} & $-(n+1) < 0$ {\bf (Negative Cases)} \\ \hline
$ \hat{a}_n ={\ds \frac{q^{n}\lb(abcdq^{2n}| q)_1\rb^2}{(acq^{n}, bcq^{n}, adq^{n}, bdq^{n}| q)_1}}$          
                                     & $\ds \hat{a}_{-(n+1)} =
                                                       - \frac{1}{cdq^{n}}
         $ \\ \hline         
$ \hat{b}_n = \lb cdq^{2n}(abcdq^{2n}| q)_1\rb \cdot$  
& $\hat{b}_{-(n+1)}=
                                         \ds     \frac{(c+d)- cdq^{n }(a+b)}{cdq^{n } (abcdq^{2n}| q)_1} 
      $\\ 
 $\ds \frac{ \lb ab(c+d)q^{n}  - (a +b) \rb }{(acq^{n}, bcq^{n}, adq^{n}, bdq^{n}| q)_1}$    &  \\ \hline    
$\ds  \hat{c}_n = 
               - \frac{1}{ab} $
                           & $\ds \hat{c}_{-(n+1)} =\frac{\lb(abcdq^{2n+1}| q)_1\rb^2}{(q^{n+1}, abq^{n+1}, cdq^{n}, abcdq^{n}| q)_1}$ \\ \hline                              
$\ds \hat{d}_{n} =   
              \frac{( abq^{n} -1 )+ab(cdq^{n-1} -1 )}{ab(1-abcdq^{2n-1})} 
                           $& $ \hat{d}_{-(n+1)}=\lb abq^{n}(abcdq^{2n+1}| q)_1 \rb \cdot$ \\ 
& $\ds \frac{\lb q(abcdq^{n}| q)_1  +cd(q^{n+1}| q)_1\rb} 
{(q^{n+1}, abq^{n+1}, cdq^{n}, abcdq^{n}| q)_1}                                                                                               
       $ \\ \hline
 $ \widetilde{\mu}_{n}=
                          t_0t_1q^{n}=
                          abcdq^{n-1}
                          $ & $\widetilde{\mu}_{-(n+1)}= 
                                                   q^{-(n+1)}
       $ \\ \hline
$ \widetilde{\zeta}_{0n}=\lb 
 \hat{a}_{n}\rb^{-1}\lb \widetilde{\mu}_{n} - \widetilde{\mu}_{-(n+1)}\rb $ 
                                                  & $ \widetilde{\zeta}_{0,-(n+1)}= \lb\hat{a}_{-(n+1)}\rb^{-1}\lb \widetilde{\mu}_{-(n+1)} - \widetilde{\mu}_{n}\rb $  \\
 \ \ \ $={ \ds  -\frac
 {(acq^{n}, bcq^{n}, adq^{n}, bdq^{n}| q)_1} 
 {q^{2n+1} (abcdq^{2n}| q)_1 }  } $     & \ \ \ $= -cdq^{-1} (abcdq^{2n}| q)_1$ \\ \hline      
$ \widetilde{\zeta}_{1n}=  \lb\hat{c}_{n}\rb^{-1}\lb \widetilde{\mu}_{-n} - \widetilde{\mu}_{n}\rb $
                                              & $ \widetilde{\zeta}_{1,-(n+1)} = \lb\hat{c}_{-(n+1)}\rb^{-1}\lb \widetilde{\mu}_{n+1} - \widetilde{\mu}_{-(n+1)}\rb $ \\
 $-abq^{-n}(abcdq^{2n-1}| q)_1$                                                   & $= {\ds  -\frac
 {(q^{n+1}, abq^{n+1}, cdq^{n}, abcdq^{n}| q)_1}
  {q^{n+1} (abcdq^{2n+1}| q)_1 }  } $  \\ \hline                                                    
\end{tabular}
\label{default}
\end{table}%

For any signs of $n,k:$
\begin{eqnarray*}
E_{-(n+1)} &=& \hat{a}_{-(n+1)} \widetilde{U}_0\Big(E_n\Big) +   \hat{b}_{-(n+1)} E_n \\
E_{n}& =& \hat{a}_{n} \widetilde{U}_0\Big(E_{-(n+1)}\Big) +  \hat{b}_{n} E_{-(n+1)} \\
E_{-n} &=& \hat{c}_{-n} \widetilde{T}_1\Big(E_{n}\Big) +   \hat{d}_{-n} E_n \\
E_{n} &=& \hat{c}_{n} \widetilde{T}_1\Big(E_{-n}\Big) +   \hat{d}_{n} E_{-n} 
\end{eqnarray*}
\begin{eqnarray*}
\widetilde{U}_0\Big(E_{k}\Big)&=&\lb \hat{a}_{-(k+1)}\rb^{-1}    \lb E_{-(k+1)}  -  \hat{b}_{-(k+1)} E_k \rb \\
\widetilde{T}_1\Big(E_{k}\Big)&=&\lb \hat{c}_{-k}\rb^{-1} \lb E_{-k} -  \hat{d}_{-k} E_k \rb \ \  (k \neq 0)
\end{eqnarray*}
\begin{eqnarray*}
\widetilde{U}_0\Big(E_{-(k+1)}\Big)&=&\lb \hat{a}_{k}\rb^{-1}    \lb E_{k}  -   \hat{b}_{k} E_{-(k+1)} \rb \\
\widetilde{T}_1\Big(E_{-(k+1)}\Big)&=&\lb \hat{c}_{k+1}\rb^{-1} \lb E_{k+1}  -   \hat{d}_{k+1} E_{-(k+1)} \rb \ \  (k \neq -1)
\end{eqnarray*}
\begin{eqnarray*}
\widetilde{\mathcal{S'}}_{0} E_{n} = \widetilde{\zeta}_{0,-(n+1)}E_{-(n+1)}  &\hspace{10mm}  & \widetilde{\mathcal{S'}}_{0} E_{-(n+1)} = \widetilde{\zeta}_{0,n}E_{n} \\
\widetilde{\mathcal{S}}_{1} E_{n}=\widetilde{\zeta}_{1, - n}E_{-n} &\hspace{10mm} & \widetilde{\mathcal{S}}_{1} E_{-n} =\widetilde{\zeta}_{1n}E_n  
\end{eqnarray*}

Since the multiplication operator
$$ X= t_0 \widetilde{T}_1^{-1}  \widetilde{ Y}  \widetilde{U }_0,$$
the above imply general $4-$dimensional invariant subspaces for $X$ and $X^{-1}.$ Simplifications of the matrix entries of these rel the $\{E_r\}$ arise. 
The simplified expressions are mostly products and quotients of q-Pochhammer symbols, with an occasional  monomial or sum of $2$ symbols factor.

\ifLONG

\section{Appendix C1: Introduction to the Shift-c Proof \label{SHIFTC1}}

The nonsymmetric Askey-Wilson polynomials $E_n(z;a,b,c,d|q)$ are symmetric under the interchange of parameters $a$ and $b,$
as well as interchange of $c$ and $d.$ 


In the main body of this paper, we have presented a full proof of Theorem \ref{ETC}. About the corresponding `shift-c' version Theorem  \ref{ETC2},
we have thus far mentioned only that the `same method of proof' can be used.
The purpose of this addendum is to write out the details of the proof of this shift-c result. 

First, note that the connection coefficients $c_{m,n}(a, e ; b, c, d | q)$ of Theorem \ref{PTCFRM} are symmetric in the parameters
$c$ and $d.$  
So $c_{m,n}(a, e ; b, c, d | q)$is not natural to relate to the nonsymmetric connection coefficients  arising from changing the basis
$\cF=\{E_r(z;a,b,c,d|q)\}=\{E_r\}$  to $\cF''= \{E_r(z;a,b,g,d|q)\}=\{E_r''\}$ as we are doing in Theorem \ref{ETC2}.

So we start with a rewriting of equations (\ref{AWCONSUM}) and (\ref{AWCON}).

Theorem 
\ref{PTCFRM} here,
a result (with different normalization)  of Askey and Wilson in  \cite{MEMOIRS}, says
$$ P_n(z; a,b,c,d | q) =\sum_{m \le n}  c_{m,n}(a, e ; b, c, d | q)  P_m(z; e,b,c,d | q) $$
where
$$c_{m,n}(a, e ; b, c, d | q) =  \frac{ (q^{n-m+1} | q )_m(bc q^m, bd q^m, cd q^m, ae^{-1}  | q )_{n-m}   }  { (q  | q )_m (abcdq^{n+m-1}, bcdeq^{2m}  | q )_{n-m} } e^{n-m}. $$

Replacing $a \to c, e \to g, c \to a,$ gives
\begin{equation} P_n(z; c,b,a,d | q) =\sum_{m \le n}  c_{m,n}(c, g ; b, a, d | q)  P_m(z; g,b,a,d | q) \end{equation}
where
\begin{equation} \label{CRSSHIFTC} c_{m,n}( c, g ; b, a, d | q) =  \frac{g^{n-m} (q^{n-m+1} | q )_m(ab q^m, bd q^m, ad q^m, cg^{-1}  | q )_{n-m}   }  { (q  | q )_m (abcdq^{n+m-1}, abdgq^{2m}  | q )_{n-m} }  \end{equation}
which we observe is symmetric in $a$ and $b.$

Keeping in mind that the $P_n(z;a,b,c,d|q)$ are symmetric in all four parameters $a,b,c, \text{ and }d,$ one form of the shift-c result for
symmetric Askey-Wilson polynomials is
\begin{equation} P_n(z; a,b,c,d | q) =\sum_{m \le n}  c_{m,n}(c, g ; a, b, d | q)  P_m(z; a,b,g,d | q) \end{equation}
with equation (\ref{CRSSHIFTC}) defining $c_{m,n}(c, g ; a, b, d | q)= c_{m,n}(c, g ; b, a, d | q).$

So the result of Theorem \ref{ETC2} may also be expressed in the following style:

The expansion formula $ \ds E_s =\sum_{r \preceq s}\lb  d^c_{r,s}(c, g ; a, b, d | q) \rb \lb c_{|r|,|s|}(c, g ; a, b, d | q) \rb E_r''  $ holds, where
\[ E_r = E_r(z; a,b,c,d|q), \quad E_r''=E_r(z; a,b,g,d|q)\]
and
\[ \quad d^c_{r,s}= \begin{cases}  
\ds  \frac{ ( abq^s, abcdq^{r+s-1} | q)_1  }  { ( abq^r, abcdq^{2s-1} | q)_1 }  & \text{ if }   r \ge 0, s\ge 0 \\ 
\ds  -\frac{  abq^r ( q^{-r -s}  | q)_1 }  { ( q^{-s}, abq^r | q)_1  }   & \text{ if }   r \ge 0, s< 0 \\ 
\ds  -  \frac{ dq^{-r-1}g (q^{-r},  abq^s, cg^{-1} q^{s+r}  | q)_1 }  { ( abcd q^{2s-1}, abdgq^{-2r-1} | q)_1 }  & \text{ if }   r < 0, s\ge 0 \\  
\ds   \frac{ ( q^{-r} , abdg q^{ -r-s-1} | q)_1  }  { ( q^{-s} , abdg q^{-2r-1}  | q)_1 }  & \text{ if }   r < 0, s< 0  
\end{cases} \]

Thus our basic working notation for the shift-c proof becomes:

For  $n \geq 0$ and any parameter values $c,g,a,b,d,\text{ and } q:$
\begin{multline*}
E_n(a,b,c,d|q) = \sum_{m=0}^n \lb \tau^c_{m,n}(c, g ; a, b, d | q) \rb E_m(a,b,g,d|q)\\  + \sum_{m=0}^n \lb \sigma^c_{-(m+1),n}(c, g ; a, b, d | q) \rb E_{-(m+1)}(a,b,g,d|q) \\
E_{-(n+1)}(a,b,c,d|q) = \sum_{m=0}^n \lb \tau^c_{-(m+1),-(n+1)}(c, g ; a, b, d | q) \rb E_{-(m+1)}(a,b,g,d|q) \\+ \sum_{m=0}^n \lb \sigma^c_{m,-(n+1)}(c, g ; a, b, d | q) \rb E_{m}(a,b,g,d|q) 
\end{multline*}
where for $k,n \geq 0:$

\begin{multline}
\tau^c_{k,n}(c, g ; a, b, d | q) =  \frac{ (q^{n-k+1} | q)_k g^{n-k}(abq^{k+1},ad q^k, bd q^k, cg^{-1}  | q )_{n-k}   }  { (q | q)_k (abcdq^{n+k}, abdgq^{2k}  | q)_{n-k} } \\
\sigma^c_{k,-(n+1)}(c, g ; a, b, d | q) = -\frac{ (abq^{k+1} | q)_{n-k} (ad q^k, bd  q^k,cg^{-1}  | q )_{n-k+1}     }  { (q | q)_k (abcdq^{n+k} , abdgq^{2k}  | q)_{n-k+1} } \\
\times abq^k g^{n-k+1} (q^{n-k+1} | q)_k \\
\sigma^c_{-(k+1),n}(c, g ; a, b, d | q) = -\frac{(abq^{k+1},cg^{-1} | q)_{n-k} (ad q^{k+1}, bd  q^{k+1}  | q )_{n-k-1}  }  { (q | q)_k (abcdq^{n+k} , abdgq^{2k+1}  | q)_{n-k} } \\
 \times dq^k g^{n-k} (q^{n-k} | q)_{k+1}   \\
\tau^c_{-(k+1),-(n+1)}(c, g ; a, b, d | q) = \frac{(q^{n-k+1} | q)_k  g^{n-k}(abq^{k+1},adq^{k+1}, bd q^{k+1}, cg^{-1} | q )_{n-k}  }  { (q | q)_k  (abcdq^{n+k+1}, abdgq^{2k+1}  | q)_{n-k} } .
\end{multline}

\ifEXTRAPROOFS

\bc \bf (Verification below  that the above $\tau^c,\sigma^c$ using the $c_{m,n}$ of equation (\ref{CRSSHIFTC}) lead to the $d^c_{r,s}$ given in Theorem \ref{ETC2}. \ec

The basic relation is 
\begin{align*}
\tau^c_{r,s}&=&d^c_{r,s}c_{|r|,|s|} & \hspace{10mm} & \text{if } (r \ge 0 \text{ and } s \ge 0)  \text{ or } (r <  0 \text{ and } s < 0)\\
\sigma^c_{r,s}&= &d^c_{r,s}c_{|r|,|s|} & \hspace{10mm} &  \text{if } (r \ge 0 \text{ and } s < 0)  \text{ or } (r <  0 \text{ and } s > 0). 
\end{align*}

\begin{proof} For $k,n \ge 0,$
\begin{enumerate}
\item 
\bmlst
d^c_{k,n} = \frac{\tau^c_{k,n}}{c_{k,n}} \\=
\left\{ \frac{ (q^{n-k+1};q)_k g^{n-k}(abq^{k+1},ad q^k, bd q^k, cg^{-1}  | q )_{n-k}   }  { (q | q)_k (abcdq^{n+k}, abdgq^{2k}  | q)_{n-k} } \right\} \cdot \\
\left\{ 
\frac
{ (q  | q )_k (abcdq^{n+k-1}, abdgq^{2k}  | q )_{n-k} }
{g^{n-k} (q^{n-k+1} | q )_k(ab q^k, bd q^k, ad q^k, cg^{-1}  | q )_{n-k}   }  
\right\} \\
=  \frac{ ( abq^n, abcdq^{n+k-1} | q)_1  }  { ( abq^k, abcdq^{2n-1} | q)_1 } 
\end{multline*}
\item 
\bmlst
d^c_{k,-(n+1)} = \frac{\sigma^c_{k,-(n+1)}}{c_{k,n+1}} \\=
\left\{  -\frac{ (q^{n-k+1} | q)_k abq^k g^{n-k+1}(abq^{k+1} | q)_{n-k} (ad q^k, bd  q^k,cg^{-1}  | q )_{n-k+1}     }  { (q | q)_k (abcdq^{n+k} , abdgq^{2k}  | q)_{n-k+1} }  \right\} \cdot \\
\left\{ 
\frac
{ (q  | q )_k (abcdq^{n+k}, abdgq^{2k}  | q )_{n-k+1} }
{g^{n-k+1} (q^{n-k+2} | q )_k(ab q^k, bd q^k, ad q^k, cg^{-1}  | q )_{n-k+1}   }  
\right\} \\
=  -\frac{  abq^k ( q^{n-k+1}  | q)_1 }  { ( q^{n+1}, abq^k | q)_1  } 
\end{multline*}
\item 
\bmlst
d^c_{-(k+1),n} = \frac{\sigma^c_{-(k+1),n}}{c_{k+1,n}} \\=
\left\{- \frac{(q^{n-k} | q)_{k+1} dq^k g^{n-k} (abq^{k+1},cg^{-1} | q)_{n-k} (ad q^{k+1}, bd  q^{k+1}  | q )_{n-k-1}  }  { (q | q)_k (abcdq^{n+k} , abdgq^{2k+1}  | q)_{n-k} } \right\} \cdot \\
\left\{ 
\frac
{ (q  | q )_{k+1} (abcdq^{n+k}, abdgq^{2k+2}  | q )_{n-k-1} }
{g^{n-k-1} (q^{n-k} | q )_{k+1}(ab q^{k+1}, bd q^{k+1}, ad q^{k+1}, cg^{-1}  | q )_{n-k-1}   }  
\right\} \\
= -  \frac{ dq^kg (q^{k+1},  abq^n, cg^{-1} q^{n-k-1}  | q)_1 }  { ( abcd q^{2n-1}, abdgq^{2k+1} | q)_1 } 
\end{multline*}
\item 
\bmlst
d^c_{-(k+1),-(n+1)} = \frac{\tau^c_{-(k+1),-(n+1)}}{c_{k+1,n+1}} \\=
\left\{  \frac{(q^{n-k+1} | q)_k  g^{n-k}(abq^{k+1},adq^{k+1}, bd q^{k+1}, cg^{-1} | q )_{n-k}  }  { (q | q)_k  (abcdq^{n+k+1}, abdgq^{2k+1} | q)_{n-k} }  \right\} \cdot \\
\left\{ 
\frac
{ (q  | q )_{k+1} (abcdq^{n+k+1}, abdgq^{2k+2}  | q )_{n-k} }
{g^{n-k} (q^{n-k+1} | q )_{k+1}(ab q^{k+1}, bd q^{k+1}, ad q^{k+1}, cg^{-1}  | q )_{n-k}   }  
\right\} \\
=  \frac{ ( q^{k+1} , abdg q^{n+k+1} | q)_1  }  { ( q^{n+1} , abdg q^{2k+1}  | q)_1 } 
\end{multline*}
\end{enumerate}
\end{proof}
\fi 

The block structure of the true transition function $\boldsymbol{\mathcal{T}}^c(c, g ; a, b, d | q)$ and  the claimed answer $T^c(c, g ; a, b, d | q),$
in the theorem above are the same as in the shift-a case. We systematically use the superscript $c$ to denote the shift-c version
of the corresponding shift-a one. And, since the parameters $a,b,d,\text{ and } q$ are often the same, we usually omit them from 
argument lists that mostly depend on $c,g, \text{ or both.}$

Our  proof of Theorem \ref{ETC2}  has three steps which we refer to as
\begin{equation} 
\label{PLANC}
\text{ \em (PROOF PLAN C)  }
\end{equation}

\begin{enumerate}
\item Show that the entries of both $T^c(c,g)$ and  $\boldsymbol{\mathcal{T}}^c(c,g) $ are rational functions of $g$ with coefficients in the filed $\mathbb{Q}(a,b,c,d,q).$ {\em (The proof is the same as in the shift-a case.)}
\item Show $\boldsymbol{\mathcal{T}}^c(c,cq)=T^c(c,cq).$
\item Show $T^c$ also satisfies the discrete co-cycle condition $$T^c(c,cq^{p+1}) = T^c(cq^p,cq^{p+1})T(c,cq^p)$$ for any $p \in \mathbb{N}.$
\end{enumerate}
Since both $T^c(c,c)$ and $\boldsymbol{\mathcal{T}}^c(c,c)$ are the identity, and  equation 
(\ref{COCYC}) 
says $\boldsymbol{\mathcal{T}}^c$ satisfies the discrete
co-cycle condition, it is immediate from 
parts $2$ and $3$ above 
that  $\boldsymbol{\mathcal{T}}^c(c,g)$ and $T^c(c,g)$ agree whenever $g=cq^p$ for
$p$ a non-negative integer. Now using 
part $1$ of Proof Plan C, 
we see that each entry of the two matrices is a rational function of $g$ agreeing with the other at infinitely many points. So they must agree
$\big($as rational functions with coefficients in $\mathbb{Q}(a,b,c,d,q) \big)$ for all $g.$ 

\section{Appendix C2: The True $\boldsymbol{\mathcal{T}}^c(c,cq)$ Matches the $\boldsymbol{T}^c(c,cq) $ of Theorem 
\ref{ETC2}
\label{SHIFTC2}
}

\ifJOLT \begin{Proposition} \else \begin{props} \fi
\label{LCDGTRC}
\begin{eqnarray}
\boldsymbol{\mathcal{T}}^{c00}_{n-1, n}(c,g) &=&\tau^c_{n-1,n}= -\Delta_{cg} \lambda_{n-1,n} +\lb \mu_{n-1,-n}(g) \rb \Delta_{ae} \mu_{-n,n} \label{DELTA00C} \\
 \boldsymbol{\mathcal{T}}^{c01}_{n, n}(c,g) &=&\sigma^c_{n,-(n+1)}=-\Delta_{cg}\mu_{n,-(n+1)}  \label{DELTA01C}  \\
\boldsymbol{\mathcal{T}}^{c10}_{n, n+1}(c,g) &=&\sigma^c_{-(n+1),n+1}= -\Delta_{cg}\mu_{-(n+1),n+1} \label{DELTA10C}  \\
\boldsymbol{\mathcal{T}}^{c11}_{n-1,n}(c,g)&=&\tau^c_{-n,-(n+1)} =-\Delta_{cg} \lambda_{-n, -(n+1)} +\lb \mu_{-n,n}(g) \rb \Delta_{cg} \mu_{n,-(n+1)} 
\end{eqnarray}
\ifJOLT \end{Proposition} \else \end{props} \fi
\begin{proof}
We write down the proof of (\ref{DELTA00C}) and  (\ref{DELTA10C}), the other two being similar.

Recall $\lambda_{r,r}=1$ for all parameters and any sign of $r.$ Then  modulo $\mathcal{R}_{-(n-1)}:$
\begin{eqnarray*}
E_{n-1}(g) =\lb  \lambda_{n-1,n-1}(g)\rb f_{n-1} & \Rightarrow & f_{n-1}=E_{n-1}(g)\\
E_{-n}(g)= \lb \lambda_{-n,-n}(g) \rb h_n + \lb \mu_{n-1,-n}(g) f_{n-1}(g) \rb & \Rightarrow \\
  h_{n}=E_{-n}(g) - \lb \mu_{n-1,-n}(g) \rb E_{n-1}(g) \\
f_n = E_n(g) -  \lb \mu_{-n,n}(g) \rb h_n - \lb \lambda_{n-1,n}(g) \rb f_{n-1} & 
\end{eqnarray*}
$$\hphantom{f_n}    =E_n(g)   -  \lb \mu_{-n,n}(g) \rb \big( E_{-n}(g) - \lb \mu_{n-1,-n}(g) \rb E_{n-1}(g)  \big) - \lb  \lambda_{n-1,n}(g)\rb E_{n-1}(g) $$
So
\begin{multline*}
E_n(c) = f_n+ \lb \mu_{-n,n}(c) \rb h_n + \lb  \lambda_{n-1,n}(c)\rb f_{n-1} \\
= E_n(g)   -  \lb \mu_{-n,n}(g) \rb \big( E_{-n}(g) - \lb \mu_{n-1,-n}(g) \rb E_{n-1}(g)  \big) - \lb  \lambda_{n-1,n}(g)\rb E_{n-1}(g) \\
+ \lb \mu_{-n,n}(c) \rb \big( E_{-n}(g) - \lb \mu_{n-1,-n}(g) \rb E_{n-1}(g) \big) +  \lb  \lambda_{n-1,n}(c)\rb E_{n-1}(g) .
\end{multline*}
Combining terms and comparing coefficients of $E_{-n}(g)$ (for  (\ref{DELTA10C}))  and $E_{n-1}(g)$ (for  (\ref{DELTA00C}))with the definition
$$
E_n(c)= \tau^c_{n,n}E_n(g) + \sigma^c_{-n,n} E_{-n}(g) + \tau^c_{n-1,n}E_{n-1}(g) \ \ \  \text{ mod }\mathcal{R}_{-(n-1)}
$$
gives the asserted formulas.
\end{proof}

The zig-zag co-degree $1$ formulas (\ref{DELTA01C}) and (\ref{DELTA10C}) above have an obvious linearity based on 
$$
 \Delta_{uw} f= \Delta_{uv} f + \Delta_{vw} f \text{ since } \lp f(w) - f(u) \rp = \lp f(v) - f(u) \rp + \lp f(w) - f(v) \rp.
$$
This immediately implies that the zig-zag co-degree $1$ matrix entries of $T^c$ satisfy what is required by the discrete co-cycle condition. We shall need
these special cases in the next section, so we record them in the corollary below.
\ifJOLT \begin{Corollary} \else \begin{cors} \fi
\label{CDG1LINC}
For any $c,g,\text{ and } f,$
\begin{eqnarray*}
\sigma^c_{n,-(n+1)}(c,g) &= &\sigma^c_{n,-(n+1)}(f,g) +  \sigma^c_{n,-(n+1)}(c,f) \\
\sigma^c_{-(n+1), n+1}(c,g) &= & \sigma^c_{-(n+1), n+1}(f,g) +  \sigma^c_{-(n+1), n+1}(c,f) .
\end{eqnarray*}
\ifJOLT \end{Corollary} \else \end{cors} \fi

\ifJOLT \begin{Proposition} \else \begin{props} \fi
\label{CDG1CGQ}

In the special case of $g=cq,$ the zig-zag co-degree $1$ transition functions satisfy 
\begin{eqnarray*}
\mathcal{T}^{c01}_{n,n}(c,cq) &=&\sigma^c_{n, -(n+1)}(c,cq) = -\frac{  abcq^{n+1} (ad q^n, bd  q^n, q^{-1} | q )_{1}   }  {  (abcdq^{2n}, abcdq^{2n+1} | q)_{1}   } \\
\mathcal{T}^{c10}_{n, n+1}(c,cq) &=& \sigma^c_{-(n+1), n+1}(c,cq) = 
 -\frac{cd    q^{n+1}   (q^{n+1},abq^{n+1},  q^{-1} | q )_{1}  }  { (abcdq^{2n+1}, abcdq^{2(n+1)} | q)_{1}  } 
\end{eqnarray*}
\ifJOLT \end{Proposition} \else \end{props} \fi
\begin{proof}
\begin{enumerate}
\item 
$$  \mathcal{T}^{c01}_{n, n}(c,cq)  =\sigma^c_{n, -(n+1)}(c,q) =-\Delta_{cg}\mu_{n,-(n+1)} $$
So 
\begin{multline*}
\mathcal{T}^{c01}_{n,n}(c,cq) = - \Big\{ \Big\{ - \frac{ \lb ( a bq^{n} | q)_1 -1 \rb (cq+d) + ( a +b)  }{  ( a b\lb cq \rb dq^{2n} | q)_1 } \Big\} \\
- \Big\{  - \frac{ \lb (abq^{n} | q)_1 -1 \rb (c+d) + (a+b)  }{  (abcdq^{2n} | q)_1 }\Big\} \Big\} \\
=\lp \frac{1} {(abcdq^{2n} | q)_1 (abcdq^{2n+1} | q)_1  } \rp \cdot \\
\Big\{ (abcdq^{2n} | q)_1\lb  - abq^{n} (cq+d) +(a+b)\rb \\
 - (abcdq^{2n+1} | q)_1\lb  - abq^{n} (c+d) +(a+b)\rb \Big\}  \\
=  \lp \frac{ 1 }{ (abcdq^{2n}, abcdq^{2n+1} | q)_1 } \rp \Big\{  - abcq^{n+1} \\
+abcq^n -abcdq^{2n}(-abdq^n+a+b)
+abcdq^{2n+1}(-abdq^n+a+b)\Big\} \\
=   \frac{ abcq^n(1-q) }{ (abcdq^{2n}, abcdq^{2n+1} | q)_1 }\Big\{1 +abd^2q^{2n} -d(a+b)q^n\Big\} \\
=   \frac{ abcq^n \lb -q(q^{-1} | q)_1 \rb}{ (abcdq^{2n}, abcdq^{2n+1} | q)_1 }\Big\{(1 -adq^n)(1 -bdq^n)\Big\}. 
\end{multline*}
\item
 $$ \mathcal{T}^{c10}_{n, n+1}(c,cq) = \sigma^c_{-(n+1), n+1}(c,cq)= -\Delta_{cg}\mu_{-(n+1),n+1} $$
 So
\bmlst
\mathcal{T}^{c10}_{n,n+1}(c,cq) = - \Big\{ \Big\{   - \frac{ (\lb cq \rb dq^{n } | q)_1 ( q^{n+1} | q  )_1 } {  ( a b \lb cq \rb dq^{2n+1} | q)_1 } \Big\} 
 - \Big\{ - \frac{ (cdq^{n} | q)_1 ( q^{n+1} | q  )_1 } {  (abcdq^{2n+1} | q)_1 }  \Big\} \Big\}\\
 = \frac{ ( q^{n+1} | q  )_1  }{ (abcdq^{2n+1} | q)_1 (abcdq^{2(n+1)} | q)_1 }\Big\{(abcdq^{2n+1} | q)_1 (cdq^{n+1} | q)_1 \\
 -(abcdq^{2(n+1)} | q)_1 (cdq^{n} | q)_1\Big\}\\
  = \frac{ ( q^{n+1} | q  )_1  }{ (abcdq^{2n+1} , abcdq^{2(n+1)} | q)_1 }\Big\{ - abcdq^{2n+1} - cdq^{n+1} + abcdq^{2(n+1)} + cdq^{n} \Big\}\\
    = \frac{ ( q^{n+1} | q  )_1  }{ (abcdq^{2n+1}, abcdq^{2(n+1)} | q)_1 }\Big\{ - abcdq^{2n+1}(q | q)_1  + cdq^{n}(q | q)_1 \Big\}\\
        = \frac{cdq^n (q, q^{n+1}, abq^{n+1} | q  )_1  }{ (abcdq^{2n+1} , abcdq^{2(n+1)}  | q)_1 }\\
                = \frac{cdq^n\lb-q (q^{-1} | q)_1 \rb ( q^{n+1}, abq^{n+1} | q  )_1  }{ (abcdq^{2n+1} , abcdq^{2(n+1)} | q)_1 } . 
\end{multline*}
\end{enumerate}

\end{proof}

\ifJOLT \begin{Proposition} \else \begin{props} \fi
\label{CDG2CGQ}

In the special case of $g=cq,$ the zig-zag co-degree $2$ transition functions satisfy:
\ 
\begin{eqnarray*}
\boldsymbol{\mathcal{T}}^{c00}_{n-1, n}(c,cq) =\tau^c_{n-1,n}= - \frac{   c  (q^n, abq^n,ad q^{n-1}, bd q^{n-1}| q )_{1}    }  { \lb (abcdq^{2n-1} | q)_{1}  \rb^2  } \\
\boldsymbol{\mathcal{T}}^{c11}_{n-1,n}(c,cq) =\tau^c_{-n,-(n+1)}= - \frac{   c  (q^n, abq^n,ad q^{n}, bd q^{n}| q )_{1}    }  { \lb (abcdq^{2n} | q)_{1}  \rb^2  }
\end{eqnarray*}

\ifJOLT \end{Proposition} \else \end{props} \fi

\begin{proof}
\begin{enumerate}
\item
By proposition \ref{LCDGTRC}
$$\boldsymbol{\mathcal{T}}^{c00}_{n-1, n}(c,cq) =-\Delta_{cg} \lambda_{n-1,n} +\lb \mu_{n-1,-n}(g) \rb \Delta_{cg} \mu_{-n,n} .$$
So 
\bmlst
\boldsymbol{\mathcal{T}}^{c00}_{n-1,n}(c,c
q)  = \\ - \Big\{
\Big\{-{\ds \frac{(q^n| q)_1}{(q| q)_1  (ab \lb cq \rb dq^{2n-1}| q)_1}} \left[(\lb cq
 \rb+d)(\ a bq^{n}| q)_1 + q( a +b) (\lb cq \rb dq^{n-1}| q)_1 \right] \Big\}\\
   - \Big\{ -{\ds \frac{(q^n| q)_1}{(q| q)_1  (abcdq^{2n-1}| q)_1}} \left[(c+d)(abq^{n}| q)_1 + q(a+b) (cdq^{n-1}| q)_1 \right]\Big\}\Big\} \\
   + \Big\{  - \frac{ \lb ( abq^{n-1} | q)_1 -1 \rb (cq+d) + (   a +b)  }{  ( ab\lb cq \rb dq^{2(n-1)} | q)_1 } \Big\} \cdot 
   \Big\{ \\
   \Big\{- \frac{ (\lb cq \rb dq^{n - 1} | q)_1 ( q^n| q  )_1 } {  ( ab \lb cq \rb dq^{2n-1} | q)_1 }  \Big\} 
-  \Big\{ - \frac{ (cdq^{n - 1} | q)_1 ( q^n| q  )_1 } {  (abcdq^{2n-1} | q)_1 }  \Big\}\Big\} . 
\end{multline*}

Using the notation $\big[a,b,c,\dots\big]=(a| q)_1(b| q)_1(c| q)_1\ldots,$ 
the abbreviation variables $u=q^{n-1},y=abcd$
and multiplying by 
$$
 \frac{ (q,q)_1 \lb ( yq^{2n-1} | q)_1 \rb^2 ( yq^{2n} | q)_1 } {(q^n| q)_1}=\frac{ \lb q,yu^2q , yu^2q , yu^2q^{2}\rb } {(uq| q)_1}
$$
it suffices to show the vanishing of 
\bmlst 
p_1=\Big\{c \big[q,  yu^2q^{2}, abuq, adu, bdu \big]    \Big\} \\
+ \Big\{ \lb ( yu^2q| q)_1\rb^2  \left[(cq+d)(a buq| q)_1 + q(  a  +b) (cd uq | q)_1 \right] \Big\}\\
- \Big\{   ( yu^2q | q)_1  ( yu^2q^{2} | q)_1 \left[(c+d)(abuq| q)_1 + q(a+b) (cdu| q)_1 \right] \Big\}  \\ 
+\Big\{   (q,q)_1 ( yu^2q | q)_1 \Big[ \lb (abu | q)_1 -1 \rb (cq+d) + (   a  +b) \Big]  (cduq | q)_1\Big\} \\
 -  \Big\{   (q,q)_1( yu^2q^{2} | q)_1 \Big[ \lb ( abu | q)_1 -1 \rb (cq+d) + (   a  +b) \Big] (cdu | q)_1 \Big\} \\
 = \Big\{c \big[q,  yu^2q^{2}, abuq, adu, bdu \big]    \Big\} \\
 + ( yu^2q| q)_1 \Big\{ (a buq| q)_1 \Big\{ c\big\{ ( yu^2q| q)_1  (q) -  ( yu^2q^{2} | q)_1 (1)  \big\} \\
                                                                 + d  \big\{ ( yu^2q | q)_1 -  ( yu^2q^{2} | q)_1 \big\} \Big\} \\
                                  +q(a+b)  \big\{ ( yu^2q | q)_1 (cd uq | q)_1  - ( yu^2q^2 | q)_1  (cd u | q)_1 \big\} \Big\}\\
                                  +(q|q)_1\big\{ -abu(cq+d) +(a+b) \big\} \big\{( yu^2q | q)_1 (cd uq | q)_1  - ( yu^2q^2 | q)_1  (cd u | q)_1  \big\} 
\end{multline*}

Note (keeping in mind, e.g., $(yu^2q)(q) - (yu^2q^2)(1) =0$) that
\begin{eqnarray*}
( yu^2q | q)_1  (q) -  ( yu^2q^{2} | q)_1 (1)  &= &q-1  = -  (q| q)_1 \\
( yu^2q | q)_1 -  ( yu^2q^{2} | q)_1&=&- yu^2q (q| q)_1 
\end{eqnarray*}
\bmlst
( yu^2q | q)_1 (cd uq | q)_1  - ( yu^2q^2 | q)_1  (cd u | q)_1 \\
= -yu^2q -cduq +yu^2q^2 +cdu \\
=u(q|q)_1 (cd -yuq) \\
= u(q|q)_1(cd -abcdq^n)  = cdu (q|q)_1 (abuq| q)_1 
\end{multline*} 
 So $p_1$ is also equal to 
\bmlst
p_2 = \Big\{c \big[q,  yu^2q^{2}, abuq, adu, bdu \big]    \Big\} \\
 + ( yu^2q| q)_1 \Big\{ (a buq| q)_1 \Big\{c \big\{ -  (q| q)_1 \big\} \\
                                                                 + d  \big\{ - yu^2q (q| q)_1 \big\} \Big\} \\
                                  +q(a+b)  \big\{ cdu (q|q)_1 (abuq| q)_1  \big\} \Big\}\\
                                  +(q|q)_1\big\{ -abu(cq+d) +(a+b) \big\} \big\{  cdu (q|q)_1 (abuq| q)_1  \big\}. 
\end{multline*}
Multiplying $p_2$ by $\lb (q|q)_1\  (abuq| q)_1 \rb^{-1}$
we are reduced to showing the vanishing of
\bmlst
p_3= \Big\{c \big[  yu^2q^{2}, adu, bdu \big]    \Big\} 
 + ( yu^2q| q)_1 \Big\{ \Big\{ c\big\{ -  1 \big\} 
                                                                 + d  \big\{ - yu^2q  \big\} \Big\} \\
                                  +q(a+b)  \big\{ cdu   \big\} \Big\}
                                  +(q|q)_1\big\{ -abu(cq+d) +(a+b) \big\} \big\{  cdu  \big\} \\ 
        =c\big[  abcdu^2q^{2}, adu, bdu \big] + ( abcdu^2q| q)_1\big\{ -c -abcd^2u^2q +cduq(a+b)\big\}\\
        +cdu(q|q)_1\big\{ -abu(cq+d) +(a+b) \big\} .
\end{multline*}
This may be viewed as an at most degree $2$ polynomial in $q$ with coefficients in the field $\mathbb{Q}(a,b,c,d,u).$

The coefficient of $q^2$ is
\bmlst
(-abcdu^2) c(1-adu)(1-bdu) + ( -abcdu^2)\big\{( -abcd^2u^2) \\ +cdu(a+b)\big\}  
+(cdu)(-1) (-abcu)\\
=(-abcdu^2\big\{  c(1-adu)(1-bdu) -abcd^2u^2 +cdu(a+b) -c \big\} =0.
\end{multline*}
So $p_3$ is  an at most degree $1$ polynomial in $q$ with coefficients in the field $\mathbb{Q}(a,b,c,d,u).$
Evaluating:
\begin{description}
\item[at $q=0$ ] 
$$
p_3(0) = c(1-adu)(1-bdu)  -c +cdu( -abdu + a +b) =0.
$$
\item[at $q= (abcdu^2)^{-1} $ ] When  $q= (abcdu^2)^{-1} ,$ note   $abcdu^{2}q^2=q$ and $abc^2d u^2q=c.$ 
So
\begin{multline*}
p_3\lp (abcdu^2)^{-1} \rp= c[q,adu,bdu] +(q|q)_1 \big\{ -c -abcd^2u^2 +(a+b)(cdu)\big\}\\
=c(q|q)_1\big\{(1-adu)(1-bdu) -1 -abd^2u^2 + du(a+b)\big\}=0.
\end{multline*}
\end{description}
Since $p_3,$ when viewed as a degree at most $1$ polynomial in $q$ with coefficients in a field, vanishes at $2$ points, it must be the
zero polynomial, proving the first formula.
\item By proposition \ref{LCDGTRC}
$$
\boldsymbol{\mathcal{T}}^{c11}_{n-1,n}(c,cq) =\\ -\Delta_{cg} \lambda_{-n, -(n+1)} +\lb \mu_{-n,n}(g) \rb \Delta_{cg} \mu_{n,-(n+1)} .
$$
So
\bmlst
\boldsymbol{\mathcal{T}}^{c11}_{n-1,n}(a,aq) =\\- \Big\{ 
\Big\{   -{\ds \frac{(q^n| q)_1}{(q| q)_1 ( ab\lb cq \rb dq^{2n}| q)_1}} \left[( \lb cq \rb +d)(  a bq^{n+1}| q)_1 + q( a +b) ( \lb cq \rb dq^{n-1}| q)_1 \right] \Big\}  \\
 - \Big\{  -{\ds \frac{(q^n| q)_1}{(q| q)_1 (abcdq^{2n}| q)_1}} \left[(c+d)(abq^{n+1}| q)_1 + q(a+b) (cdq^{n-1}| q)_1 \right]   \Big\} \Big\}   \\
+\Big\{  - \frac{ ( \lb cq \rb dq^{n - 1} | q)_1 ( q^n| q  )_1 } {  ( ab\lb cq \rb dq^{2n-1} | q)_1 }\Big\} \cdot  \Big\{   \\
 \Big\{    \frac{ a bq^{n}  (  \lb cq \rb +d) - (  a +b)  }{  ( ab \lb cq \rb dq^{2n} | q)_1 }\Big\}  \\
 -   \Big\{ \frac{  abq^{n} (c+d) - (a+b)  }{  (abcdq^{2n} | q)_1 } \Big\}\Big\}  .
\end{multline*}  
Using the notation $\big[a,b,c,\dots\big]=(a| q)_1(b| q)_1(c| q)_1\ldots,$ 
the abbreviation variables $u=q^{n-1},y=abcd$
and multiplying by 
$$
 \frac{ (q,q)_1  \lb (yu^2q^2| q)_1 \rb^2 (yu^2q^3| q)_1 } {(u q | q)_1  }=\frac{\lb q,yu^2q^2, yu^2q^2, yu^2q^3 \rb } { (u q | q)_1},
$$
it suffices to show the vanishing of
\bmlst 
p_1=\Big\{ c\lb q, abuq , aduq, bduq, yu^2q^3 \rb     \Big\} \\
+ \Big\{ \lb (yu^{2}q^2, yu^2q^2 \rb   \left\{ (cq+d)(a buq^{2}| q)_1 + q(a +b) (cd uq | q )_{1} \right\} \Big\} \\
- \Big\{   \lb yu^2q^2, yu^2q^3 \rb \left\{(c+d)(abuq^2| q)_1 + q(a+b)(cd u | q )_{1}  \right\}\Big\}  \\
+ \Big\{  \lb q, yu^2q^2, cduq \rb \Big\{    -abuq (cq+d) +  (a+b)\Big\} \\
 -  \lb q, yu^2q^3, cduq \rb \Big\{   - abu q (c+d) + (a+b)  \Big\}  \Big\} \\
 = \Big\{ c\lb q, abuq , aduq, bduq, yu^2q^3 \rb     \Big\} \\
 + (yu^{2}q^2 | q)_1 \Big\{  (a buq^{2}| q)_1 \Big\{c \big\{ ( yu^2q^2 | q)_1  (q) -  ( yu^2q^{3} | q)_1 (1) \big\} \\
                      + d \big\{ ( yu^2q^2 | q)_1 -  ( yu^2q^{3} | q)_1 \big\} \Big\} \\
                      + q(a+b)\big\{ ( yu^2q^2 | q)_1 (cd uq | q)_1  - ( yu^2q^3 | q)_1  (cd u | q)_1 \big\} \Big\} \\
                     + [q, cduq] \Big\{ (-abuq) \Big\{ c \big\{ ( yu^2q^2 | q)_1  (q) -  ( yu^2q^{3} | q)_1 (1) \big\}  \\
                     + d \big\{  ( yu^2q^2 | q)_1 -  ( yu^2q^{3} | q)_1 \big\} \Big\} \\
                              + (a+b) \big\{  ( yu^2q^2 | q)_1 -  ( yu^2q^{3} | q)_1 \big\}       \Big\} .
\end{multline*}
Note (keeping in mind, e.g., $(yu^2q^2)(q) - (yu^2q^3)(1) =0$) that
\begin{eqnarray*}
( yu^2q^2 | q)_1  (q) -  ( yu^2q^{3} | q)_1 (1)  &=& q-1  = -  (q| q)_1 \\
( yu^2q^2 | q)_1 -  ( yu^2q^{3} | q)_1&=& - yu^2q^2 (q| q)_1 
\end{eqnarray*}
\begin{multline*}
( yu^2q^2 | q)_1 (cd uq | q)_1  - ( yu^2q^3 | q)_1  (cd u | q)_1 \\
= -yu^2q^2 -cduq +yu^2q^3 +cdu \\
=u(q|q)_1 (cd -yuq^2) \\
= u(q|q)_1(cd -abcdq^{n+1}) \\ = cdu (q|q)_1 (abuq^2 | q)_1 .
\end{multline*} 
So $p_1$ is also equal to 
\bmlst
p_2 =  \Big\{ c\lb q, abuq , aduq, bduq, yu^2q^3 \rb     \Big\} 
 + (yu^{2}q^2 | q)_1 \Big\{  (a buq^{2}| q)_1 \Big\{c \big\{  -  (q| q)_1 \big\} \\
                      + d \big\{  - yu^2q^2 (q| q)_1 \big\} \Big\} 
                      + q(a+b)\big\{  cdu (q|q)_1 (abuq^2 | q)_1  \big\} \Big\} \\
                     + [q, cduq] \Big\{ (-abuq) \Big\{ c \big\{  -  (q| q)_1 \big\}                       
                     + d \big\{  - yu^2q^2 (q| q)_1  \big\} \Big\} \\
                              + (a+b) \big\{  - yu^2q^2 (q| q)_1 \big\}       \Big\} .                     
\end{multline*}
Multiplying $p_2$ by $\lb (q|q)_1\   \rb^{-1}$
we are reduced to showing the vanishing of
\bmlst
p_3=  \Big\{ c\lb abuq , aduq, bduq, yu^2q^3 \rb     \Big\} \\
 + (yu^{2}q^2 | q)_1 \Big\{  (a buq^{2}| q)_1 \Big\{c \big\{  -  1 \big\} 
                      + d \big\{  - yu^2q^2 \big\} \Big\} \\
                      + q(a+b)\big\{  cdu (abuq^2 | q)_1  \big\} \Big\} \\
                     + [q, cduq ] \Big\{ (-abuq) \Big\{ c \big\{  -  1\big\}                        
                     + d \big\{  - yu^2q^2  \big\} \Big\} 
                              + (a+b) \big\{  - yu^2q^2  \big\}       \Big\} . 
\end{multline*}
In $p_3$ above, every appearance of $u$ has at least one power of $q$ next to it, so it
is natural to switch to
$$
\widetilde{u} =uq=q^n.
$$
Also remember $y=abcd.$ So it is sufficient to show the vanishing of the polynomial
\bmlst
p_4=   c\lb ab\widetilde{u} , ad\widetilde{u}, bd\widetilde{u}, abcd\widetilde{u}^2q \rb     \\
 + (abcd\widetilde{u}^{2} | q)_1 \Big\{ - (a b\widetilde{u}q | q)_1 \Big\{c 
                      +    abcd^2\widetilde{u}^2  \Big\} 
                      +   cd(a+b)\widetilde{u} (ab\widetilde{u}q | q)_1   \Big\} \\
                     + [q, cd\widetilde{u} ] \Big\{ ab\widetilde{u} \Big\{ c                        
                     +    abcd^2\widetilde{u}^2  \Big\} 
                              -     abcd\widetilde{u}^2   (a+b)     \Big\} 
\end{multline*}
as an at most degree $1$ polynomial in $q$ with coefficients in the field \\ $\mathbb{Q}(a,b,c,d,\widetilde{u}).$ 
Evaluating at the two points
\begin{description}
\item[at $q=1$]
$$p_4(1 )= \lb abcd\widetilde{u}^2, ab \widetilde{u} \rb \big\{  c \lb ad \widetilde{u}, bd \widetilde{u} \rb  -(c+abcd^2\widetilde{u}^2) +cd\widetilde{u}(a+b)\big\} =0 . $$
\item[at $q= (abcd\widetilde{u}^2)^{-1} $ ] When  $q= (abcd\widetilde{u}^2)^{-1} ,$ note   $ abcd\widetilde{u}^{2} =q^{-1}, $ \  $ab\widetilde{u}q=\lb cd\widetilde{u} \rb^{-1} $ and $\lb ab\widetilde{u} \rb^{-1} = cd\widetilde{u} q .$ 
So
\begin{multline*}
p_4\lp (abcd\widetilde{u}^2)^{-1} \rp=  (abcd\widetilde{u}^{2} | q)_1 \Big\{ - (a b\widetilde{u}q | q)_1 \Big\{c 
                      +    abcd^2\widetilde{u}^2  \Big\} \\
                      +   cd(a+b)\widetilde{u} (ab\widetilde{u}q | q)_1   \Big\} \\
                     + [q, cd\widetilde{u} ] \Big\{ ab\widetilde{u} \Big\{ c                        
                     +    abcd^2\widetilde{u}^2  \Big\} 
                              -     abcd\widetilde{u}^2   (a+b)     \Big\} \\
=  (q^{-1}| q)_1 (a b\widetilde{u}q | q)_1 \Big\{ -  \lp c 
                      +    dq^{-1} \rp 
                      +   cd(a+b)\widetilde{u}  \Big\} \\
                     + [q, ( ab\widetilde{u}q )^{-1} ] \Big\{ ab\widetilde{u} \Big\{ c                        
                     +    d q^{-1} \Big\} 
                              -    q^{-1}  (a+b)     \Big\} \\ 
  =  (q^{-1}| q)_1 (a b\widetilde{u}q | q)_1 \Big\{ -  \lp c 
                      +   d q^{-1} \rp
                      +   cd(a+b)\widetilde{u}  \Big\} \\
                     + (q^{-1})(ab\widetilde{u} q) [q, ( ab\widetilde{u}q )^{-1} ] \Big\{  \lp  c                        
                     +  d  q^{-1} \rp
                              -    q^{-1} \lb ab\widetilde{u} \rb^{-1} (a+b)     \Big\} \\        
  =  (q^{-1}| q)_1 (a b\widetilde{u}q | q)_1 \Big\{ \Big\{ -  \lp c 
                      +  d  q^{-1} \rp
                      +   cd(a+b)\widetilde{u}  \Big\} \\
                     +  \Big\{  \lp c                        
                     +  d  q^{-1} \rp 
                              -    q^{-1} \lp cd\widetilde{u} q \rp(a+b)     \Big\} \Big\}=0.                                                                                                                     
\end{multline*}
\end{description}
Thus $p_4$ being of degree at most $1$ and vanishing at $2$ points implies $p_4$ is identically $0$ and the proposition proof is complete.
\end{enumerate}

\end{proof}

Propositions \ref{CDG1CGQ} and \ref{CDG2CGQ} complete the proof of the second step of (\ref{PLANC}), Proof Plan C.

\section{Appendix C3: Proof of the Discrete Co-Cycle Identity for $\boldsymbol{T}^c$  \label{COCYCPRFC}}

The block form of the equation 
$$T^c(cq^p,cq^{p+1}) T^c(c,cq^{p}) = T^c(c,cq^{p+1})$$

is
\bmlst
\begin{bmatrix}
T^{c00} (cq^p,cq^{p+1})& T^{c01}(cq^p,cq^{p+1} \\
T^{c10}(cq^p,cq^{p+1}) & T^{c11}(cq^p,cq^{p+1} \\
\end{bmatrix}
\begin{bmatrix}
T^{c00}(c,cq^p) & T^{c01}(c,cq^p) \\
T^{c10}(c,cq^p) & T^{c11}(c,cq^p) \\
\end{bmatrix} \\
=
\begin{bmatrix}
T^{c00}(c,cq^{p+1}) & T^{c01}(c,cq^{p+1}) \\
T^{c10}(c,cq^{p+1}) & T^{c11}(c,cq^{p+1}) \\
\end{bmatrix}.
\end{multline*}
So we have the four matrix equations:
\begin{eqnarray*}
T^{c00}(c,cq^{p+1})  =  \Big[ T^{c00} (cq^p,cq^{p+1}) \Big] \Big[T^{c00}(c,cq^p) \Big] &  +  & \Big[T^{c01} (cq^p,cq^{p+1}) \Big] \Big[T^{c10}(c,cq^p) \Big] \\
T^{c01}(c,cq^{p+1})  =  \Big[T^{c00} (cq^p,cq^{p+1}) \Big] \Big[T^{c01}(c,cq^p) \Big] & +  & \Big[T^{c01} (cq^p,cq^{p+1})  \Big]\Big[T^{c11}(c,cq^p) \Big] \\
T^{c10}(c,cq^{p+1})  =  \Big[T^{c10} (cq^p,cq^{p+1}) \Big] \Big[T^{c00}(c,cq^p) \Big] &+  & \Big[T^{c11} (cq^p,cq^{p+1}) \Big] \Big[T^{c10}(c,cq^p) \Big]\\
T^{c11}(c,cq^{p+1})  =  \Big[T^{c10} (cq^p,cq^{p+1}) \Big]\Big[T^{c01}(c,cq^p) \Big] & +&  \Big[T^{c11} (cq^p,cq^{p+1} ) \Big] \Big[T^{c11}(c,cq^p) \Big].
\end{eqnarray*}
Keep in mind, as described in the shift-a case 
at the end of section 
\ref{ALMSYM},
that the entries of each $T^{cij}(cq^p,cq^{p+1})$ are zero except (possibly) on the diagonal and the superdiagonal.

In terms of $\tau^c \text{ and } \sigma^c,$ these equations are:
\bml
\label{COCYC00DSPC}
\tau^c_{k,n}(c,cq^{p+1}) = \big[ \tau^c_{k,k}(cq^p,cq^{p+1})\big] \big[ \tau^c_{k,n}(c,cq^p) \big]+ \big[ \tau^c_{k,k+1}(cq^p,cq^{p+1})\big] \cdot \big\{ \\ \big[ \tau^c_{k+1, n}(c,cq^p)\big] \big\}
   +\big[ \sigma^c_{k,-(k+1)}(cq^p,cq^{p+1})\big] \big[ \sigma^c_{-(k+1),n}(c,cq^p) \big] 
   \end{multline}  
   \bml
   \label{COCYC01DSPC}
\sigma^c_{k,-(n+1)}(c,cq^{p+1}) = \big[ \tau^c_{k,k}(cq^p,cq^{p+1})\big] \big[ \sigma^c_{k,-(n+1)}(c,cq^p) \big]+ \big[ \tau^c_{k,k+1}(cq^p,cq^{p+1})\big] \cdot \big\{ \\\big[ \sigma^c_{k+1, -(n+1))}(c,cq^p)\big] \big\}
    +\big[ \sigma^c_{k,-(k+1)}(cq^p,cq^{p+1})\big] \big[ \tau^c_{-(k+1),-(n+1)}(c,cq^p) \big]
       \end{multline}  
   \bml
   \label{COCYC10DSPC}
  \sigma^c_{-(k+1), n}(c,cq^{p+1}) = 
  \big[ \sigma^c_{-(k+1),k+1}(cq^p,cq^{p+1})\big] \big[ \tau^c_{k+1, n}(c,cq^p)\big] \\
    +\big[ \tau^c_{-(k+1),-(k+1)}(cq^p,cq^{p+1})\big]
    \big\{ 
     \big[ \sigma^c_{-(k+1),n}(c,cq^p) \big] \big\} \\
    + \big[ \tau^c_{-(k+1),-(k+2)}(cq^p,cq^{p+1})\big] \big[ \sigma^c_{-(k+2), n}(c,cq^p)\big]  
       \end{multline}  
   \bml
   \label{COCYC11DSPC}
 \tau^c_{-(k+1),-(n+1)}(c,cq^{p+1}) = 
    \big[ \sigma^c_{-(k+1),k+1}(cq^p,cq^{p+1})\big] \big[ \sigma^c_{k+1, -(n+1)}(c,cq^p)\big] \\
  + \big[ \tau^c_{-(k+1),-(k+1)}(cq^p,cq^{p+1})\big] 
  \big[ \tau^c_{-(k+1),-(n+1)}(c,cq^p) \big] \\
  + \big[ \tau^c_{-(k+1),-(k+2)}(cq^p,cq^{p+1})\big] \big[ \tau^c_{-(k+2), -(n+1)}(c,cq^p)\big] 
  \end{multline}
\ifJOLT \begin{Lemma} \else \begin{lemmas} \fi For $k,n \ge 0$
\label{URTC}
\ 

\begin{enumerate}
\item $\ds \frac{ \tau^c_{k,n}(c,cq^{p+1}) } { \sigma^c_{-(k+1),n}(c,cq^{p}) }= - \frac{ q^{n-2k} (adq^k, bdq^k, q^{-(p+1)} | q)_1} {d(q^{n-k}, q^{n-k-p-1} | q )_1 }. $

\item $\ds  \frac{ \sigma^c_{k,-(n+1)}(c,cq^{p+1}) } { \tau^c_{-(k+1),-(n+1)}(c,cq^{p}) } =  - \frac{abcq
^{n+p+1} (adq^k, bdq^k, q^{-(p+1)}  | q )_1 } {  (abcdq^{n+k}, abcdq^{n+k+p+1} | q )_1 }.  $

\item $\ds  \frac{ \sigma^c_{-(k+1),n}(c,cq^{p+1}) } { \tau^c_{k+1,n}(c,cq^{p}) } = - \frac{ cdq^{n+p} ( q^{k+1}, abq^{k+1}, q^{-(p+1)}  | q )_1 } {  ( abcdq^{n+k}, abcdq^{n+k+p+1} | q )_1 } .$

\item $\ds \frac{ \tau^c_{-(k+1),-(n+1)}(c,cq^{p+1}) } { \sigma^c_{k+1,-(n+1)}(c,cq^{p}) } =  - \frac{q^{n-2k-1} ( q^{k+1}, abq^{k+1}, q^{-(p+1)} | q )_1 } { ab ( q^{n-k}, q^{n-k-p-1} | q )_1 } . $ 
\end{enumerate}
\ifJOLT \end{Lemma} \else \end{lemmas} \fi

\ifEXTRAPROOFS
\begin{proof}
\begin{enumerate}
\item 
\bmlst
\lp 
\frac{ (q^{n-k+1} | q)_k \lb cq^{p+1} \rb^{n-k}(abq^{k+1},ad q^k, bd q^k, q^{-(p+1)}  | q )_{n-k}   }  { (q | q)_k (abcdq^{n+k}, abd \lb cq^{p+1} \rb q^{2k}  | q)_{n-k} } 
\rp \cdot \\
\lp 
 - \frac
 { (q | q)_k (abcdq^{n+k} , abd \lb cq^{p} \rb q^{2k+1}  | q)_{n-k} }
 {(q^{n-k} | q)_{k+1} dq^k \lb cq^{p} \rb^{n-k} (abq^{k+1}, q^{-p} | q)_{n-k} (ad q^{k+1}, bd  q^{k+1}  | q )_{n-k-1}  }  
\rp \\
= - \frac{ q^{n-2k} (adq^k, bdq^k, q^{-(p+1)} | q)_1} {d(q^{n-k}, q^{n-k-p-1} | q )_1 } .  \\
\end{multline*}
\item 
\bmlst
\lp 
-\frac{ (q^{n-k+1} | q)_k abq^k  \lb cq^{p+1} \rb ^{n-k+1}(abq^{k+1}| q)_{n-k} (ad q^k, bd  q^k,q^{-(p+1)}  | q )_{n-k+1}     }  { (q | q)_k (abcdq^{n+k} , abd \lb cq^{p+1} \rb q^{2k}  | q)_{n-k+1} }  
\rp \cdot \\
\lp  
\frac
{ (q | q)_k  (abcdq^{n+k+1}, abd \lb cq^p \rb q^{2k+1}  | q)_{n-k} }
{(q^{n-k+1} | q)_k   \lb cq^p \rb ^{n-k}(abq^{k+1},adq^{k+1}, bd q^{k+1}, q^{-p} | q )_{n-k}  }  
\rp \\
= - \frac{abcq
^{n+p+1} (adq^k, bdq^k, q^{-(p+1)}  | q )_1 } {  (abcdq^{n+k}, abcdq^{n+k+p+1} | q )_1 }. \\
\end{multline*}
\item 
\bmlst
\Big( 
- \big\{ (q^{n-k} | q)_{k+1} dq^k \lb cq^{p+1} \rb^{n-k} \big\} \cdot \\
 \frac
 {(abq^{k+1}, q^{-(p+1)}  | q)_{n-k} (ad q^{k+1}, bd  q^{k+1}  | q )_{n-k-1}  } 
 { (q | q)_k (abcdq^{n+k} , abd \lb cq^{p+1} \rb q^{2k+1}  | q)_{n-k} } 
\Big) 
\Big( 
\big\{ (q | q)_{k+1} \big\} \cdot \\
\Big\{ \frac
{  (abcdq^{n+k+1}, abd \lb cq^{p} \rb q^{2k+2}  | q)_{n-k-1} }
{ (q^{n-k} | q)_{k+1} \lb cq^{p} \rb^{n-k-1}(abq^{k+2},ad q^{k+1}, bd q^{k+1}, q^{-p}  | q )_{n-k-1}   }  \Big\}
\Big) \\ 
= - \frac{ cdq^{n+p} ( q^{k+1}, abq^{k+1}, q^{-(p+1)}  | q )_1 } {  ( abcdq^{n+k}, abcdq^{n+k+p+1} | q )_1 }. \\
\end{multline*}
\item 
\bmlst
\lp
\frac{(q^{n-k+1} | q)_k   \lb cq^{p+1} \rb^{n-k}(abq^{k+1},adq^{k+1}, bd q^{k+1}, q^{-(p+1)} | q )_{n-k}  }  { (q | q)_k  (abcdq^{n+k+1}, abd \lb cq^{p+1} \rb q^{2k+1}  | q)_{n-k} }   
\rp \cdot \\
\lp 
-\frac
{ (q | q)_{k+1} (abcdq^{n+k+1} , abd  \lb cq^{p} \rb q^{2k+2}  | q)_{n-k} }
{ (q^{n-k} | q)_{k+1} abq^{k+1}  \lb cq^{p} \rb^{n-k}(abq^{k+2} | q)_{n-k-1} (ad q^{k+1}, bd  q^{k+1} , q^{-p}  | q )_{n-k}     }  
\rp
\\
= - \frac{q^{n-2k-1} ( q^{k+1}, abq^{k+1}, q^{-(p+1)} | q )_1 } { ab ( q^{n-k}, q^{n-k-p-1} | q )_1 }.\\
\end{multline*}
\end{enumerate}

\end{proof}
\fi 

\ifJOLT \begin{Lemma} \else \begin{lemmas} \fi  For $k,n \ge 0:$
\label{VRTC}
\ 

\begin{enumerate}
\item $\ds \frac{ \tau^c_{k,n}(c,cq^{p}) } { \sigma^c_{-(k+1),n}(c,cq^{p}) } = - \frac {(adq^k,bdq^k,abcdq^{n+k+p} | q)_1 } { dq^k ( q^{n-k}, abcd q^{2k+p} | q)_1} .$

\item $\ds \frac{ \sigma^c_{k,-(n+1)}(c,cq^{p}) } { \tau^c_{-(k+1),-(n+1)}(c,cq^{p}) } = - \frac{abc q^{k+p}  ( adq^k, bdq^k, q^{n-k-p}  | q )_1 } {  (abcdq^{n+k}, abcdq^{2k+p}  | q )_1 }.$

\item  $\ds \frac{ \sigma^c_{-(k+1),n}(c,cq^{p}) } { \tau^c_{k+1,n}(c,cq^{p}) } = - \frac {cdq^{k+p}(q^{k+1},q^{n-k-p-1}, abq^{k+1} | q )_1 } { (abcdq^{n+k}, abcd q^{2k+p+1} | q )_1}.$

\item $\ds \frac{ \tau^c_{-(k+1),-(n+1)}(c,cq^{p}) } { \sigma^c_{k+1,-(n+1)}(c,cq^{p}) } =  - \frac{ ( q^{k+1}, abq^{k+1}, abcdq^{n+k+p+1}  | q )_1 } { abq^{k+1}  ( q^{n-k}, abcdq^{2k+p+1}  | q )_1 }  . $

\end{enumerate}
\ifJOLT \end{Lemma} \else \end{lemmas} \fi

\ifEXTRAPROOFS

\begin{proof}

\begin{enumerate}
\item 
\bmlst
\lp 
\frac{ (q^{n-k+1} | q)_k \lb cq^{p} \rb^{n-k}(abq^{k+1},ad q^k, bd q^k, q^{-p}  | q )_{n-k}   }  { (q | q)_k (abcdq^{n+k}, abd \lb cq^{p} \rb q^{2k}  | q)_{n-k} }
\rp \cdot \\
\lp 
- \frac
 { (q | q)_k (abcdq^{n+k} , abd \lb cq^{p} \rb q^{2k+1}  | q)_{n-k} }
 {(q^{n-k} | q)_{k+1} dq^k \lb cq^{p} \rb^{n-k} (abq^{k+1}, q^{-p}  | q)_{n-k} (ad q^{k+1}, bd  q^{k+1}  | q )_{n-k-1}  }  
\rp \\
= - \frac {(adq^k,bdq^k,abcdq^{n+k+p} | q)_1 } { dq^k ( q^{n-k}, abcd q^{2k+p} | q)_1}  .\\
\end{multline*}
\item 
\bmlst
\lp  
-\frac{ (q^{n-k+1} | q)_k abq^k  \lb cq^p \rb ^{n-k+1}(abq^{k+1} | q)_{n-k} (ad q^k, bd  q^k,q^{-p}  | q )_{n-k+1}     }  { (q | q)_k (abcdq^{n+k} , abd \lb cq^p \rb q^{2k}  | q)_{n-k+1} } 
\rp \cdot \\
\lp
\frac
{ (q | q)_k  (abcdq^{n+k+1}, abd \lb cq^p \rb q^{2k+1}  | q)_{n-k} }
{(q^{n-k+1} | q)_k   \lb cq^p \rb ^{n-k}(abq^{k+1},adq^{k+1}, bd q^{k+1}, q^{-p} | q )_{n-k}  }  
\rp \\
= - \frac{abc q^{k+p}  ( adq^k, bdq^k, q^{n-k-p}  | q )_1 } {  (abcdq^{n+k}, abcdq^{2k+p}  | q )_1 }. \\
\end{multline*}
\item 
\bmlst
\lp 
 - \frac{(q^{n-k} | q)_{k+1} dq^k \lb cq^{p} \rb ^{n-k} (abq^{k+1}, q^{-p}  | q)_{n-k} (ad q^{k+1}, bd  q^{k+1}  | q )_{n-k-1}  }  { (q | q)_k (abcdq^{n+k} , abd \lb cq^{p} \rb q^{2k+1}  | q)_{n-k} } 
\rp \\
\lp 
\frac
{ (q | q)_{k+1} (abcdq^{n+k+1}, abd \lb cq^{p} \rb q^{2k+2}  | q)_{n-k-1} } 
{ (q^{n-k} | q)_{k+1} \lb cq^{p} \rb ^{n-k-1}(abq^{k+2},ad q^{k+1}, bd q^{k+1},  q^{-p}   | q )_{n-k-1}   }  
\rp \\
= - \frac {cdq^{k+p}(q^{k+1},q^{n-k-p-1}, abq^{k+1} | q )_1 } { (abcdq^{n+k}, abcd q^{2k+p+1} | q )_1}. \\
\end{multline*}
\item 
\bmlst
\lp
\frac{(q^{n-k+1} | q)_k   \lb cq^{p} \rb ^{n-k}(abq^{k+1},adq^{k+1}, bd q^{k+1}, q^{-p} | q )_{n-k}  }  { (q | q)_k  (abcdq^{n+k+1}, abd \lb cq^{p} \rb q^{2k+1}  | q)_{n-k} }  
  \rp \cdot \\
\lp 
-\frac
{ (q | q)_{k+1} (abcdq^{n+k+1} , abd  \lb cq^{p} \rb q^{2k+2}  | q)_{n-k} }
{ (q^{n-k} | q)_{k+1} abq^{k+1}  \lb cq^{p} \rb^{n-k}(abq^{k+2} | q)_{n-k-1} (ad q^{k+1}, bd  q^{k+1} , q^{-p}  | q )_{n-k}     }  
  \rp \\
= - \frac{ ( q^{k+1}, abq^{k+1}, abcdq^{n+k+p+1}  | q )_1 } { abq^{k+1}  ( q^{n-k}, abcdq^{2k+p+1}  | q )_1 } .  \\
\end{multline*}
\end{enumerate}

\end{proof}
\fi 

\bc \bf The $T^{c00}$ Identity\label{T00C} \ec
\ifJOLT \begin{Proposition} \else \begin{props} \fi
\label{T00PROPC}
When $0 \le k \le n-1$
\bmlst
\tau^c_{k,n}(c,cq^{p+1}) = \big[ \tau^c_{k,k}(cq^p,cq^{p+1})\big] \big[ \tau^c_{k,n}(c,cq^p) \big] + \big[ \tau^c_{k,k+1}(cq^p,cq^{p+1})\big] \big[ \tau^c_{k+1, n}(c,cq^p)\big] \\
  +\big[ \sigma^c_{k,-(k+1)}(cq^p,cq^{p+1})\big] \big[ \sigma^c_{-(k+1),n}(c,cq^p) \big].
\end{multline*}
 And 
$$  \tau^c_{n,n}(c,cq^{p+1}) = \big[ \tau^c_{n,n}(cq^p,cq^{p+1})\big] \big[ \tau^c_{n,n}(c,cq^p) \big].$$
\ifJOLT \end{Proposition} \else \end{props} \fi

\begin{proof}
The second identity just says $1=1 \cdot 1.$


For the first, upon dividing by $\sigma^c_{-(k+1),n}(c,cq^p),$ we see it is sufficient to prove
\bmlst
\frac{\tau^c_{k,n}(c,cq^{p+1})}{\sigma^c_{-(k+1),n}(c,cq^p)} = \left\{ \tau^c_{k,k}(cq^p,cq^{p+1})\right\} \left\{ \frac{\tau^c_{k,n}(c,cq^p)}{\sigma^c_{-(k+1),n}(c,cq^p)} \right\} \\ + \left\{ \tau^c_{k,k+1}(cq^p,cq^{p+1})\right\} \left\{ \frac{\tau^c_{k+1, n}(c,cq^p)}{\sigma^c_{-(k+1),n}(c,cq^p)} \right\} 
  +\left\{ \sigma^c_{k,-(k+1)}(cq^p,cq^{p+1})\right\} .
\end{multline*}
Using Lemmas \ref{URTC} and \ref{VRTC}, this means we need to show
\bml
\label{EQN1ac}
\left\{ 
- \frac{ q^{n-2k} (adq^k, bdq^k, q^{-(p+1)} | q)_1} {d(q^{n-k}, q^{n-k-p-1} | q )_1 }
\right\} 
  =   \left\{ 1 \right\} \cdot
 \left\{  - \frac {(adq^k,bdq^k,abcdq^{n+k+p} | q)_1 } { dq^k ( q^{n-k}, abcd q^{2k+p} | q)_1}\right\} \\
+  
\left\{ \frac{ (q^{2} | q)_k \lb cq^{p+1} \rb (abq^{k+1},ad q^k, bd q^k, q^{-1}  | q )_{1}   }  { (q | q)_k (ab \lb cq^p \rb dq^{2k+1}, abd \lb cq^{p+1} \rb q^{2k}  | q)_{1} }  \right\} \cdot \\
\left\{ 
- \frac 
{ (abcdq^{n+k}, abcd q^{2k+p+1} | q )_1}
{cdq^{k+p}(q^{k+1},q^{n-k-p-1}, abq^{k+1} | q )_1 } 
\right\} \\
  + \left\{ 
   -\frac{ (q | q)_k abq^k  \lb cq^{p+1} \rb (ad q^k, bd  q^k,q^{-1}  | q )_{1}     }  { (q | q)_k (ab \lb cq^p \rb dq^{2k} , abd \lb cq^{p+1} \rb q^{2k}  | q)_{1} } \right\} 
\cdot 1 \\
\end{multline}

Multiplying (\ref{EQN1ac}) by
$$ 
\frac{ dq^{2k+p+1} (q, q^{n-k}, q^{n-k-p-1}, abcdq^{2k+p}, abcdq^{2k+p+1} | q )_1 } { (adq^k, bdq^k | q)_1 }
$$ 
we see it is sufficient to show the vanishing of the  polynomial $p_1$ below.

We will eventually reduce this identity to the vanishing of a 1-variable polynomial in $q$ with coefficients in the field $\mathbb{Q}(a,b,c,d,y,u,v,w)$ 
with the property that when 
$$
y=abcd \hspace{5mm} u= q^n \hspace{5mm} v= q^k\hspace{5mm} w= q^p
$$
we obtain a unit in the coefficient field times the difference between the two sides of equation (\ref{EQN1ac}) above. So, effectively, we can use the variables $y,u,v,w$
as abbreviations for these expressions.

\begin{multline}
p_1= \left\{q^{n+p+1} (yq^{2k+p} | q)_1 (yq^{2k+p+1} | q)_1  (q| q)_1       (q^{-(p+1)}| q )_{1}   \right\}   \\
  -\left\{ q^{k+p+1} (yq^{n+k+p} | q)_1  (yq^{2k+p+1} | q)_1     (q| q)_1   (q^{n-k-p-1} | q )_{1}    \right\} \\
    -\left\{ q^{k+p+2} (q^{n-k}|q)_1 (yq^{n+k} | q)_1  (yq^{2k+p} | q)_1      (q^{-1} | q )_{1}    
       \right\}  \\
    - \left\{q^{3k+2p+2} (q^{n-k}| q)_1y (q| q)_1  (q^{n-k-p-1}|q)_1   (q^{-1} | q )_{1}  
      \right\} 
\end{multline}
\begin{multline*}
\hphantom{p_1}=\left\{ q^{n+p+1} (yq^{2k+p} | q)_1 (yq^{2k+p+1} | q)_1  (q| q)_1  \lb - q^{-(p+1)}(q^{p+1}| q)_1\rb   \right\} \\
-  \left\{  q^{k+p+1} (yq^{n+k+p} | q)_1  (yq^{2k+p+1} | q)_1     (q| q)_1  \lb q^{-(k+p+1)} \lb (q^{n}| q)_1-(q^{k+p+1}| q)_1 \rb \rb    \right\} \\
  - \left\{ q^{k+p+2}  (yq^{n+k} | q)_1  (yq^{2k+p} | q)_1   \lb q^{-k} \lb (q^{n}| q)_1-(q^k| q)_1 \rb \rb       \lb - q^{-1}(q| q)_1 \rb   
       \right\} \\
 - \left\{  yq^{3k+2p+2}(q| q)_1   \lb q^{-k} \lb (q^{n}| q)_1-(q^k| q)_1 \rb \rb   \right. \cdot \\
\left.  \lb q^{-(k+p+1)} \lb (q^{n}| q)_1-(q^{k+p+1}| q)_1 \rb \rb      \lb- q^{-1}(q| q)_1 \rb   \right\} 
\end{multline*}
\begin{multline*}
\hphantom{p_1}=- \left\{ q^{n} (yq^{2k+p} | q)_1 (yq^{2k+p+1} | q)_1  (q| q)_1  \lb (q^{p+1}| q)_1\rb   \right\} \\
-  \left\{   (yq^{n+k+p} | q)_1  (yq^{2k+p+1} | q)_1     (q| q)_1   \lb (q^{n}| q)_1-(q^{k+p+1}| q)_1 \rb     \right\} \\
  + \left\{ q^{p+1}  (yq^{n+k} | q)_1  (yq^{2k+p} | q)_1     \lb (q^{n}| q)_1-(q^k| q)_1 \rb        \lb (q| q)_1 \rb   
       \right\} \\
 + \left\{  yq^{k+p}(q| q)_1     \lb (q^{n}| q)_1-(q^k| q)_1 \rb     
  \lb (q^{n}| q)_1-(q^{k+p+1}| q)_1 \rb       \lb (q| q)_1  \rb  \right\} 
\end{multline*}
Using the notation $\big[a,b,c,\dots\big]=(a| q)_1(b| q)_1(c| q)_1\ldots,$ 
and multiplying by $\lb  (q| q)_1 \rb^{-1}$ 
it suffices to show the vanishing of 
\bmlst
p_2 = -u[yv^2w, yv^2wq, wq] 
        -(vwq - u)[yuvw, yv^2wq] \\
       +wq(v-u)[yuv, yv^2w]
       +yvw(v-u)(vwq-u) (q| q)_1.
\end{multline*}

This expression may be interpreted as a one variable polynomial of degree at most $2$ in $q$ with coefficients in the field $\mathbb{Q}(y,u,v,w).$

The coefficient of $q^2$ is
\bmlst
-u (1-yv^2w) (-yv^2w) (-w) - (vw) (1-yuvw) (-yv^2w) +yvw (v-u) (vw) (-1) \\
=yv^2w^2 \big\{ -u (1-yv^2w) - (v) (1-yuvw)(-1)+(v-u)(-1)   \big\} 
=0.
\end{multline*}
So the polynomial $p_2$ is of degree at most $1$ in $q.$
Evaluating 
\begin{description}
\item[at $q=0$ ] $p_2(0)=-u(1-yv^2w)+u(1-yuvw) +yvw(v-u)(-u)=0$
\item[at $ q=u(vw)^{-1}$ ] Note at this point $yv^2wq=yuv$ and $wq=uv^{-1}.$ So
$$ p_2\big(u(vw)^{-1}\big) =- u[yv^2w, yuv](1-uv^{-1}) +uv^{-1}(v-u)[yv^2w,yuv]=0. $$
\end{description}

Thus the polynomial $p_2$ is $0$ and the proposition is proven.
\end{proof}

\bc \bf The $T^{c01}$ Identity \ec

\ifJOLT \begin{Proposition} \else \begin{props} \fi
\label{T01PROPC}
When $0 \le k \le n-1$
\bml
\sigma^c_{k,-(n+1)}(c,cq^{p+1}) = \big[ \tau^c_{k,k}(cq^p,cq^{p+1})\big] \big[ \sigma^c_{k,-(n+1)}(c,cq^{p}) \big] \\
+ \big[ \tau^c_{k,k+1}(cq^p,cq^{p+1})\big] \big[ \sigma^c_{k+1, -(n+1))}(c,cq^{p})\big] \\
  +\big[ \sigma^c_{k,-(k+1)}(cq^p,cq^{p+1})\big] \big[ \tau^c_{-(k+1),-(n+1)}(c,cq^{p}) \big].
\end{multline}
And
\bml
\sigma^c_{n,-(n+1)}(c,cq^{p+1}) = \big[ \tau^c_{n,n}(cq^p,cq^{p+1})\big] \big[ \sigma^c_{n,-(n+1)}(c,cq^{p}) \big] \\
  +\big[ \sigma^c_{n,-(n+1)}(cq^p,cq^{p+1})\big] \big[ \tau^c_{-(n+1),-(n+1)}(c,cq^{p}) \big].
\end{multline}

\ifJOLT \end{Proposition} \else \end{props} \fi

\begin{proof} The second is immediate from corollary \ref{CDG1LINC} together with the observation that $\tau^c_{r,r}=1$ for any sign of $r.$

For the first, upon dividing by $\tau^c_{-(k+1),-(n+1)}(c,cq^{p}),$ we see it is sufficient to prove
\bmlst
\frac{\sigma^c_{k,-(n+1)}(c,cq^{p+1}) } {\tau^c_{-(k+1),-(n+1)}(c,cq^{p})} 
= \left\{ \tau^c_{k,k}(cq^p,cq^{p+1})\right\} \left\{ \frac{\sigma^c_{k,-(n+1)}(c,cq^{p})} {\tau^c_{-(k+1),-(n+1)}(c,cq^{p})} \right\} \\
+ \left\{ \tau^c_{k,k+1}(cq^p,cq^{p+1})\right\} \left\{ \frac{\sigma^c_{k+1, -(n+1))}(c,cq^{p}) } {\tau^c_{-(k+1),-(n+1)}(c,cq^{p})} \right\}
  +\left\{ \sigma^c_{k,-(k+1)}(cq^p,cq^{p+1})\right\} \cdot 1.
\end{multline*}
That means we need to show
\bml
\label{EQN1bc}
\left\{
- \frac{abcq
^{n+p+1} (adq^k, bdq^k, q^{-(p+1)}  | q )_1 } {  (abcdq^{n+k}, abcdq^{n+k+p+1} | q )_1 }
 \right\}  
  =   \left\{1  \right\} 
  \left\{  
  - \frac{abc q^{k+p}  ( adq^k, bdq^k, q^{n-k-p}  | q )_1 } {  (abcdq^{n+k}, abcdq^{2k+p}  | q )_1 }
 \right\} 
   \\
  + \left\{
  \frac{ (q^{2} | q)_k \lb cq^{p+1} \rb(abq^{k+1},ad q^k, bd q^k, q^{-1}  | q )_{1}   }  { (q | q)_k (ab \lb cq^{p} \rb dq^{2k+1}, abd \lb cq^{p+1} \rb q^{2k}  | q)_{1} } 
 \right\} 
   \cdot \\
 \left\{  
 - \frac
 { abq^{k+1}  ( q^{n-k}, abcdq^{2k+p+1}  | q )_1 }
 { ( q^{k+1}, abq^{k+1}, abcdq^{n+k+p+1}  | q )_1 } 
  \right\} \\
 + \left\{ 
  -\frac{ (q | q)_k abq^k  \lb cq^{p+1} \rb  (ad q^k, bd  q^k,q^{-1}  | q )_{1}     }  { (q | q)_k (ab \lb cq^{p} \rb dq^{2k} , abd \lb cq^{p+1} \rb q^{2k}  | q)_{1} } 
 \right\} 
   \cdot   
  1
\end{multline}

We will eventually reduce this identity to the vanishing of a 1-variable polynomial in $y$ with coefficients in the field $\mathbb{Q}(a,b,c,d,u,v,w,q)$ 
with the property that when 
$$
y=abcd \hspace{5mm} u= q^n \hspace{5mm} v= q^k\hspace{5mm} w= q^p
$$
we obtain a unit in the coefficient field times the difference between the two sides of equation (\ref{EQN1bc})  above. So, effectively, we can use the variables $y,u,v,w$
as abbreviations for the above expressions.

Multiplying (\ref{EQN1bc})  by  
\begin{equation}
-\frac{ (yq^{n+k}, yq^{n+k+p+1}, yq^{2k+p}, yq^{2k+p+1}  |q)_1 } { abc (adq^k, bdq^k  | q)_1 }
\end{equation}
we see it is sufficient to show the vanishing of the  polynomial $p_1$ below. We use $y$ as an abbreviation for $abcd.$
\bc (We chose the $-1$ factor above to make $p_1$ a perfect match for the $p_1$ in the shift-a case.)\ec
\bmlst
p_1 =  -\left\{  (yq^{2k +p} | q)_{1}   (yq^{2k +p+1} | q)_{1}  
   q^{n+p+1}  (q^{-(p+1)}| q)_1  \right\}\\
  +   \left\{1  \right\}   
  \left\{    (yq^{n+k +p+1} | q)_{1}  (yq^{2k +p+1} | q)_{1}  q^{k+p}  
    (q^{n-k-p}| q)_1      \right\} \\
 + \left\{\frac{ (q^{n-k}| q)_1(yq^{n+k} | q)_{1}   (y q^{2k+p} | q)_{1}   q^{k+p+2} 
   (q^{-1}| q)_1 
  }  { (q| q)_1 
   }  \right\} \\ 
 +  \left\{ (y q^{n+k} | q)_{1} (y q^{n+k+p+1} | q)_{1}    q^{k+p+1}      (q^{-1}| q)_1   \right\} 
   \end{multline*}      
\bmlst
\hphantom{p_1 }=  -\left\{  (yq^{2k +p} | q)_{1}   (yq^{2k +p+1} | q)_{1}  q^{n+p+1} \lb   - q^{-(p+1)}(q^{p+1}| q)_1 \rb   \right\}\\
  +   \left\{1  \right\}   
  \left\{    (yq^{n+k +p+1} | q)_{1}  (yq^{2k +p+1} | q)_{1}   q^{k+p} 
 \lb q^{-(k+p)} \lb (q^n| q)_1 -(q^{k+p}| q)_1 \rb  \rb     \right\}   \\
 +\frac{ \left\{ \lb q^{-k} \lb(q^n| q)_1 -(q^k| q)_1 \rb  \rb (yq^{n+k} | q)_{1}   (y q^{2k+p} | q)_{1}  q^{k+p+2} 
  \lb -q^{-1}(q| q)_1 \rb 
  \right\}}  { (q| q)_1 } \\ 
 +  \left\{ (y q^{n+k} | q)_{1} (y q^{n+k+p+1} | q)_{1}    q^{k+p+1}     \lb -q^{-1}(q| q)_1 \rb  \right\} 
  \end{multline*}
 \bmlst
\hphantom{p_1 }=  \left\{  (yq^{2k +p} | q)_{1}   (yq^{2k +p+1} | q)_{1}   q^{n} \lb (q^{p+1}| q)_1 \rb   \right\}\\
  +   \left\{1  \right\}   
  \left\{    (yq^{n+k +p+1} | q)_{1}  (yq^{2k +p+1} | q)_{1}   
   \lb \lb (q^{n}| q)_1-(q^{k+p}| q)_1\rb    \rb
      \right\} \\
 - \left\{  q^{p+1} (yq^{n+k} | q)_{1}   (y q^{2k+p} | q)_{1}  
    \right\} 
 \lb \lb (q^{n}| q)_1-(q^{k}| q)_1 \rb \rb \\
 -  \left\{ (y q^{n+k} | q)_{1} (y q^{n+k+p+1} | q)_{1}    q^{k+p}     \lb (q| q)_1 \rb  \right\}. 
  \end{multline*} 
 
  Using the notation $\big[a,b,c,\dots\big]=(a| q)_1(b| q)_1(c| q)_1\ldots,$ 
setting $y=abcd, u=q^n,v=q^k, \text{ and } w=q^p,$ we see it is sufficient to show the vanishing of

\bmlst
p_2=u\lb yv^2w, yv^2wq, wq\rb 
+(vw-u)\lb yuvwq, yv^2wq \rb\\
-wq(v-u) \lb yuv, yv^2w\rb
-vw \lb yuv, yuvwq,q\rb. 
\end{multline*}
 
 This expression may be interpreted as a one variable polynomial of degree at most $2$ in $y$ with coefficients in the field $\mathbb{Q}(u,v,w,q).$
 Evaluating
 \begin{description}
\item [at $y=0$] $p_2(0)=u(1-wq)+(vw-u) -wq(v-u)-vw(1-q) =0.$
\item [at $y=(v^2wq)^{-1}$] When $y=(v^2wq)^{-1},$ note $yv^2w =q^{-1},$ $yuv= u(vwq)^{-1},$ and $yuvwq=uv^{-1}.$ So
\bmlst
p_2\big( (v^2wq)^{-1} \big) = -wq(v-u) \lb u(vwq)^{-1}, q^{-1} \rb -vw \lb u(vwq)^{-1}, uv^{-1}, q \rb\\
= (1-u(vwq)^{-1})\big\{  -wq(v-u)  \lp -q^{-1} (1-q) \rp -vw(1-uv^{-1}) (1-q) \big\} \\
= w(1-u(vwq)^{-1})(1-q)\big\{ v-u -v(1-uv^{-1}) \big\} =0.
\end{multline*}
\item [at $y=(uv)^{-1}$] When $y=(uv)^{-1},$ note $yv^2w= u^{-1}vw,$ $yv^2wq=u^{-1}vwq ,$ and $yuvwq=wq .$ So
\bmlst
p_2\big( (uv)^{-1} \big) = (1-u^{-1}vwq)(1-wq) \big\{u(1 -  u^{-1}vw) +vw-u\big\}=0. 
\end{multline*}
\end{description}
Thus the polynomial $p_2$ is $0$ and the proposition is proven.

\end{proof}

\bc \bf The $T^{c10}$ Identity \ec

\ifJOLT \begin{Proposition} \else \begin{props} \fi
\label{T10PROPC}
When $0 \le k \le n-2$
\bml
    \sigma^c_{-(k+1), n}(c,cq^{p+1}) =  \big[ \sigma^c_{-(k+1),k+1}(cq^p,cq^{p+1})\big] \big[ \tau^c_{k+1, n}(c,cq^{p})\big] \\
+\big[ \tau^c_{-(k+1),-(k+1)}(cq^p,cq^{p+1})\big] \big[ \sigma^c_{-(k+1),n}(c,cq^{p}) \big] \\
  + \big[ \tau^c_{-(k+1),-(k+2)}(cq^p,cq^{p+1})\big] \big[ \sigma^c_{-(k+2), n}(c,cq^{p})\big] 
\end{multline}
And
\bmlst
   \sigma^c_{-n, n}(c,cq^{p+1}) =  \big[ \sigma^c_{-n,n}(cq^p,cq^{p+1})\big] \big[ \tau^c_{n,n}(c,cq^{p})\big] \\
+\big[ \tau^c_{-n,-n}(cq^p,cq^{p+1})\big] \big[ \sigma^c_{-n,n}(c,cq^{p}) \big]. 
\end{multline*}
\ifJOLT \end{Proposition} \else \end{props} \fi

\begin{proof} The second is immediate from corollary \ref{CDG1LINC} together with the observation that $\tau^c_{rr}=1$ for any sign of $r.$

For the first, upon dividing by $\tau^c_{k+1,n}(c,cq^{p}),$ we see it is sufficient to prove
\bmlst
    \frac{ \sigma^c_{-(k+1), n}(c,cq^{p+1}) } { \tau^c_{k+1,n}(c,cq^{p})} =  \left\{ \sigma^c_{-(k+1),k+1}(cq^p,cq^{p+1}) \right\} \cdot 1 \\
+\left\{  \tau^c_{-(k+1),-(k+1)}(cq^p,cq^{p+1})  \right\} \left\{ \frac{ \sigma^c_{-(k+1),n}(c,cq^{p}) } { \tau^c_{k+1,n}(c,cq^{p}) } \right\} \\
  + \left\{  \tau^c_{-(k+1),-(k+2)}(cq^p,cq^{p+1}) \right\} \left\{ \frac { \sigma^c_{-(k+2), n}(c,cq^{p}) } { \tau^c_{k+1,n}(c,cq^{p}) } \right\} .
\end{multline*}
That means we need to show
\bml
\label{EQN1cc}\ 
\left\{
- \frac{ cdq^{n+p} ( q^{k+1}, abq^{k+1}, q^{-(p+1)}  | q )_1 } {  ( abcdq^{n+k}, abcdq^{n+k+p+1} | q )_1 }
\right\}
\\
= 
\left\{
- \frac{(q | q)_{k+1} dq^k  \lb cq^{p+1} \rb  (abq^{k+1},q^{-1} | q)_{1}   }  { (q | q)_k (ab \lb cq^{p} \rb dq^{2k+1} , abd \lb cq^{p+1} \rb q^{2k+1}  | q)_{1} } 
\right\} 
\cdot 1 \\
+  \left\{1 \right\}  \cdot 
\left\{
- \frac {cdq^{k+p}(q^{k+1},q^{n-k-p-1}, abq^{k+1} | q )_1 } { (abcdq^{n+k}, abcd q^{2k+p+1} | q )_1}
\right\}
 \\
+ 
\left\{
\frac{(q^{2} | q)_k   \lb cq^{p+1} \rb (abq^{k+1},adq^{k+1}, bd q^{k+1}, q^{-1} | q )_{1}  }  { (q | q)_k  (ab \lb cq^{p} \rb dq^{2k+2}, abd \lb cq^{p+1} \rb q^{2k+1}  | q)_{1} } 
\right\} 
\cdot \\   
\left\{  - \frac 
{ dq^{k+1} ( q^{n-k-1}, abcd q^{2k+p+2} | q)_1}
{(adq^{k+1},bdq^{k+1},abcdq^{n+k+p+1} | q)_1 } 
\right\} \\
\end{multline}

We will eventually reduce this identity to the vanishing of a 1-variable polynomial in $y$ with coefficients in the field $\mathbb{Q}(a,b,c,d,u,v,w,q)$ 
with the property that when 
$$
y=abcd \hspace{5mm} u= q^n \hspace{5mm} v= q^k\hspace{5mm} w= q^p
$$
we obtain a unit in the coefficient field times the difference between the two sides of equation (\ref{EQN1cc})  above. So, effectively, we can use the variables $y,u,v,w$
as abbreviations for the above expressions.

Multiplying (\ref{EQN1cc})  by
$$ 
- \frac{ q( q, yq^{n+k}, yq^{n+k+p+1}, yq^{2k+p+1}, yq^{2k+p+2}  | q )_1 } {cd (q^{k+1}, abq^{k+1}  | q)_1 } \\
$$
we see it is sufficient to show the vanishing of the  polynomial $p_1$ below.


\bmlst
p_1  =  
- \left\{ q^{n +p+1}  (q| q)_1  (yq^{2k+p+1}| q)_1(yq^{2k+p+2}| q)_1 
 (q^{-(p+1)} | q )_{1}    \right\}  \\
+ 
  \left\{ q^{k+p+2}  (q| q)_1  (yq^{n+k}| q)_1 (yq^{n+k+p+1}| q)_1  (q^{-1}| q)_1    \right\}  \\  
   + \left\{1 \right\}  \cdot    
\left\{ q^{k+p+1} (q| q)_1  (q^{n-k-p-1}| q)_1  (yq^{n+k+p+1}| q)_1 (yq^{2k+p+2}| q)_1   
    \right\} \\
+ \left\{ q^{k+p+3} (q^{n-k-1}| q)_1  (yq^{n+k}| q)_1 (yq^{2k+p+1}| q)_1    
(q^{-1}| q)_1    \right\} \cdot 
\end{multline*}
\bmlst 
\hphantom{p_1}  =  
-\left\{q^{n +p+1}(q| q)_1  (yq^{2k+p+1}| q)_1(yq^{2k+p+2}| q)_1  
  \lb - q^{-(p+1)}(q^{p+1}| q)_1  \rb   \right\}  \\
+  
     \left\{ q^{k+p+2}  (q| q)_1  (yq^{n+k}| q)_1 (yq^{n+k+p+1}| q)_1    \lb  - q^{-1}(q| q)_1 \rb   \right\}   \\  
  +  \left\{1 \right\}  \cdot   
\left\{ q^{k+p+1}  (q| q)_1 (yq^{n+k+p+1}| q)_1 (yq^{2k+p+2}| q)_1 
  \lb q^{-(k+p+1)}   \lb (q^{n}| q)_1-(q^{k+p+1}| q)_1 \rb \rb  \right\} \\
+  \left\{ q^{k+p+3}  (yq^{n+k}| q)_1 (yq^{2k+p+1}| q)_1 \lb  q^{-(k+1)} \lb (q^{n}| q)_1-(q^{k+1}| q)_1 \rb\ \rb  
 \lb - q^{-1}(q| q)_1 \rb   \right\}   
\end{multline*}
\bmlst
\hphantom{p_1} =  
\left\{ q^{n }   (yq^{2k+p+1}| q)_1(yq^{2k+p+2}| q)_1 (q| q)_1 
  (q^{p+1}| q)_1     \right\}\\
 -   
     \left\{  q^{k+p+1}   (yq^{n+k}| q)_1 (yq^{n+k+p+1}| q)_1   \lb (q| q)_1 \rb^2  \right\}   \hspace{50mm} \\  
  +  \left\{1 \right\}  \cdot   %
\left\{ (yq^{n+k+p+1}| q)_1 (yq^{2k+p+2}| q)_1 (q| q)_1 
   \lb  \lb (q^{n}| q)_1-(q^{k+p+1}| q)_1 \rb \rb  \right\} \\
- \left\{q^{p+1} (yq^{n+k}| q)_1 (yq^{2k+p+1}| q)_1 \lb  \lb (q^{n}| q)_1-(q^{k+1}| q)_1 \rb\ \rb 
 (q| q)_1   \right\} .
\end{multline*}
Using the notation $\big[a,b,c,\dots\big]=(a| q)_1(b| q)_1(c| q)_1\ldots,$ 
recalling our abbreviation variables $u,v,\text{ and } w,$
and multiplying by $ (q| q)^{-1}$ 
it suffices to show the vanishing of 
\bmlst
p_2 =  u[ yv^2wq,yv^2wq^2,wq]
-vwq[yuv, yuvwq,q] \\
 +(vwq-u ) [yuvwq, yv^2wq^2 ] 
-wq(vq-u)[yuv, yv^2wq] .
\end{multline*}

This expression may be interpreted as a one variable polynomial of degree at most $2$ in $y$ with coefficients in the field $\mathbb{Q}(u,v,w,q).$
Evaluating
\begin{description}
\item[at $y=0$ ] 
$$
p_2(0)= u(1-wq) -vwq(1-q) +(vwq-u) -wq(vq-u)=0.
$$
\item[at $y=(v^2wq)^{-1}$ ] Note at this point $yv^2wq^2=q, \ yuvwq=uv^{-1}, \text{ and }yuv= u(vwq)^{-1}.\ $ So
\bmlst
p_2\big((v^2wq)^{-1} \big) =(1-uv^{-1}) \big\{  -vwq(1- u(vwq)^{-1}) (1-q) +(vwq-u)(1-q)\big\} \\
                                         =(1-uv^{-1})(1-q)(vwq-u)  \big\{ -1+1 \big\} =0.
\end{multline*}
\item[at $y=(uvwq)^{-1}$ ] Note at this point $yv^2wq=u^{-1}v, \ yv^2wq^2= u^{-1}v q\  \text{ and }yuv=(wq)^{-1} $ So
\bmlst
p_2\big((uvwq)^{-1} \big)= (1-u^{-1}v)\big\{ u(1-u^{-1}vq)(1-wq)-wq(vq-u)(1-(wq)^{-1} \big\} \\
                                       = (1-u^{-1}v)(vq-u)(1-wq)\big\{-1 +1\big\} =0.
\end{multline*}
\end{description}

Thus $p_2$ being of degree at most  $2$ and vanishing at $3$ points implies $p_2$ is identically $0$ and the proposition proof is complete.

\end{proof}

\bc \bf The $T^{c11}$ Identity \ec
\ifJOLT \begin{Proposition} \else \begin{props} \fi
\label{T11PROPC}
When $0 \le k \le n-1$
\bml
  \tau^c_{-(k+1),-(n+1)}(c,cq^{p+1}) = 
   \big[ \sigma^c_{-(k+1),k+1}(cq^p,cq^{p+1})\big] \big[ \sigma^c_{k+1, -(n+1)}(c,cq^{p})\big] \\
+ \big[ \tau^c_{-(k+1),-(k+1)}(cq^p,cq^{p+1})\big] \big[ \tau^c_{-(k+1),-(n+1)}(c,cq^{p}) \big] \\
+ \big[ \tau^c_{-(k+1),-(k+2)}(cq^p,cq^{p+1})\big] \big[ \tau^c_{-(k+2), -(n+1)}(c,cq^{p})\big].
\end{multline}
And 
$$  \tau^c_{-(n+1),-(n+1)}(c,cq^{p+1}) = \big[ \tau^c_{-(n+1),-(n+1)}(cq^p,cq^{p+1})\big] \big[ \tau^c_{-(n+1),-(n+1)}(c,cq^{p}) \big].$$
\ifJOLT \end{Proposition} \else \end{props} \fi

\begin{proof}
The second identity just says $1=1 \cdot 1.$

For the first, upon dividing by $ \sigma^c_{k+1,-(n+1)}(c,cq^{p}) ,$ we see it is sufficient to prove
\bmlst
  \frac{ \tau^c_{-(k+1),-(n+1)}(c,cq^{p+1}) } { \sigma^c_{k+1,-(n+1)}(c,cq^{p}) }  = 
   \left\{ \sigma^c_{-(k+1),k+1}(cq^p,cq^{p+1}) \right\} \cdot 1  \\
+ \left\{ \tau^c_{-(k+1),-(k+1)}(cq^p,cq^{p+1})\right\} \left\{  \frac{ \tau^c_{-(k+1),-(n+1)}(c,cq^{p}) } { \sigma^c_{k+1,-(n+1)}(c,cq^{p}) } \right\} \\
+ \left\{ \tau^c_{-(k+1),-(k+2)}(cq^p,cq^{p+1})\right\} \left\{  \frac{ \tau^c_{-(k+2), -(n+1)}(c,cq^{p}) }{ \sigma^c_{k+1,-(n+1)}(c,cq^{p}) } \right\}.
\end{multline*}
That means we need to show
\bml
\label{EQN1dc}
\left\{
- \frac{q^{n-2k-1} ( q^{k+1}, abq^{k+1}, q^{-(p+1)} | q )_1 } { ab ( q^{n-k}, q^{n-k-p-1} | q )_1 } 
\right\}
\\
= 
\left\{
- \frac{(q | q)_{k+1} dq^k  \lb cq^{p+1} \rb (abq^{k+1}, q^{-1} | q)_{1}  }  { (q | q)_k (ab \lb cq^{p} \rb dq^{2k+1} , abd \lb cq^{p+1} \rb q^{2k+1}  | q)_{1} } 
\right\} 
  \cdot 1  \\
+\left\{ 1
\right\}
\cdot
\left\{
- \frac{ ( q^{k+1}, abq^{k+1}, abcdq^{n+k+p+1}  | q )_1 } { abq^{k+1}  ( q^{n-k}, abcdq^{2k+p+1}  | q )_1 }  
\right\}
\\
+\left\{
\frac{(q^{2} | q)_k   \lb cq^{p+1} \rb (abq^{k+1},adq^{k+1}, bd q^{k+1}, q^{-1} | q )_{1}  }  { (q | q)_k  (ab \lb cq^{p} \rb dq^{2k+2}, abd \lb cq^{p+1} \rb q^{2k+1}  | q)_{1} } 
\right\} 
\cdot \\
    \left\{
    - \frac
    {  (abcdq^{n+k+1}, abcdq^{2k+p+2}  | q )_1 }
    {abc q^{k+p+1}  ( adq^{k+1}, bdq^{k+1}, q^{n-k-p-1}  | q )_1 } 
    \right\} 
\end{multline}

Multiplying (\ref{EQN1dc}) by

$$
- \frac{abq^{2k+p+2} (q, q^{n-k}, q^{n-k-p-1}, yq^{2k+p+1}, yq^{2k+p+2}  | q)_1} { ( q^{k+1}, abq^{k+1} | q)_1} \\
$$
we see it is sufficient to show the vanishing of the  polynomial $p_1$ below.

We will eventually reduce this identity to the vanishing of a 1-variable polynomial in $y$ with coefficients in the field $\mathbb{Q}(a,b,c,d,u,v,w,q)$ 
with the property that when 
$$
y=abcd \hspace{5mm} u= q^n \hspace{5mm} v= q^k\hspace{5mm} w= q^p
$$
we obtain a unit in the coefficient field times the difference between the two sides of equation (\ref{EQN1dc}) above. So, effectively, we can use the variables $y,u,v,w$
as abbreviations for the above expressions.

\bmlst
p_1=
-\left\{ q^{n+p+1} (q| q)_1   (yq^{2k+p+1} | q)_1 (yq^{2k+p+2} | q)_1        (q^{-(p+1)} | q )_{1}    \right\} \\
 +
   \left\{ yq ^{3k+2p+3}  (q| q)_1     (q^{n-k}| q)_{1}    (q^{n-k-p-1} | q )_{1}      (q^{-1} | q )_{1}     
   \right\} \\
  +  \left\{1 \right\} 
   \left\{ q^{k+p+1} (q| q)_1    (q^{n-k-p-1} | q )_{1}  (yq^{n+k+p+1} | q)_1  (yq^{2k+p+2} | q)_1   
   \right\} \\
 +  \left\{ q^{k+p+2}  (q^{n-k}| q)_1  (yq^{n+k+1} | q)_1  (yq^{2k+p+1} | q)_1   (q^{-1} | q )_{1}        
  \right\}  
 \end{multline*} 
\bmlst
\hphantom{p_1}=
-\left\{q^{n+p+1}  (q| q)_1  (yq^{2k+p+1} | q)_1 (yq^{2k+p+2} | q)_1       \lb  - q^{-(p+1)}(q^{p+1}| q)_1 \rb    \right.   \\
 +
   \Big\{   y  q^{3k+2p+3}    (q| q)_1  \lb  q^{-k} \lb (q^{n}| q)_1-(q^k| q)_1 \rb\rb    \cdot \\
   \lb q^{-k-p-1} \lb (q^n| q)_1 - (q^{k+p+1}| q)_1 \rb  \rb   \lb  - q^{-1}(q| q)_1  \rb      \Big\} \\
  +  \left\{1 \right\} 
   \left\{   q^{k+p+1}  (q| q)_1   \lb q^{-k-p-1} \lb (q^n| q)_1 - (q^{k+p+1}| q)_1 \rb  \rb   (yq^{n+k+p+1} | q)_1   (yq^{2k+p+2} | q)_1        \right\} \\
+  \left\{  q^{k+p+2}  \lb q^{-k} \lb (q^{n}| q)_1-(q^k| q)_1 \rb \rb    (yq^{n+k+1} | q)_1     (yq^{2k+p+1} | q)_1   \lb  - q^{-1}(q| q)_1  \rb    
  \right\}  
 \end{multline*}
\bmlst
\hphantom{p_1}=
\left\{ q^{n}  (q| q)_1   (yq^{2k+p+1} | q)_1 (yq^{2k+p+2} | q)_1        (q^{p+1}| q)_1      \right\} \\
 - \left\{   y  q^{k+p+1}   \lb (q| q)_1  \rb^2    \lb (q^{n}| q)_1-(q^k| q)_1 \rb     
   \lb (q^n| q)_1 - (q^{k+p+1}| q)_1 \rb        \right\} \\
+ \left\{1 \right\} 
   \left\{  (q| q)_1    \lb (q^n| q)_1 - (q^{k+p+1}| q)_1 \rb          (yq^{n+k+p+1} | q)_1  (yq^{2k+p+2} | q)_1    \right\} \\
 -  \left\{ q^{p+1}   (q| q)_1  \lb (q^{n}| q)_1-(q^k| q)_1 \rb   (yq^{n+k+1} | q)_1  (yq^{2k+p+1} | q)_1      
  \right\}  
 \end{multline*}
Using the notation $\big[a,b,c,\dots\big]=(a| q)_1(b| q)_1(c| q)_1\ldots,$ 
recalling our abbreviation variables $u,v,\text{ and } w,$
and multiplying by $\lb (q| q)_1 \rb^{-1}$ 
it suffices to show the vanishing of 
\bmlst
p_2 = u[yv^2wq, yv^2wq^2, wq ] 
- yvwq(v-u)(vwq-u)(q| q)_1\\
+ (vwq-u) [yuvwq,yv^2wq^2] 
 -wq(v-u) [yuvq,yv^2wq] 
\end{multline*}
This expression may be interpreted as a one variable polynomial of degree at most $2$ in $y$ with coefficients in the field $\mathbb{Q}(u,v,w,q).$
Evaluating 
\begin{description}
\item[at $y=0$ ] 
$$
p_2(0)= u(1-wq) + (vwq-u)  -wq(v-u)=0.
$$
\item[at $y=(v^2wq)^{-1}$ ] Note at this point $yv^2wq^2=q, \ yvwq=v^{-1}, \text{ and } yuvwq=uv^{-1}.$ So
\bmlst
p_2\big((v^2wq)^{-1} \big)=(vwq-u) \big\{ -v^{-1}(v-u)(1-q) + (1-uv^{-1})(1-q)  \big\} \\
                                         =(vwq-u) (1-uv^{-1})(1-q)  \big\{ -1+1 \big\} =0.
\end{multline*}
\item[at $y=(v^2wq^2)^{-1}$ ] Note at this point $yv^2wq=q^{-1}, \ yvwq=(vq)^{-1}, \text{ and } yuvq=u(vwq)^{-1}.$ So
\bmlst
p_2\big((v^2wq^2)^{-1} \big)=(v-u)  \big\{-(vq)^{-1}(vwq-u)(1-q) - wq(1- u(vwq)^{-1})(1-q^{-1}) \big\} \\
                                       =(v-u) (1-q)\big\{ -( w -u (vq)^{-1}  ) + (w- u(vq)^{-1})\big\}=0. 
\end{multline*}
\end{description}
Thus $p_2$ being of degree $2$ and vanishing at $3$ points implies $p_2$ is identically $0$ and the proposition proof is complete.

\end{proof}

The combination of Propositions \ref{PLANA1} (the $\boldsymbol{\mathcal{T}}(c, g ; a, b, d | q) $ version),  \ref{CDG1CGQ}, \ref{CDG2CGQ}, \ref{T00PROPC}, \ref{T01PROPC}, \ref{T10PROPC}, and \ref{T11PROPC}   finishes the proof of all $3$ steps of (\ref{PLANC}), Proof Plan C, and so completes the proof of Theorem \ref{ETC2}.

\fi 


\end{document}